\documentclass[letterpaper,11pt,reqno]{amsart} 
\RequirePackage[OT1]{fontenc}
\usepackage[portrait,margin=3cm]{geometry}  
\usepackage{mathrsfs}
\usepackage[colorlinks,citecolor=blue,urlcolor=blue,linkcolor=blue,pagecolor=blue,linktocpage=true,backref=true]%
{hyperref}
\usepackage[foot]{amsaddr}
\usepackage{amssymb,amsthm,amsfonts,amsbsy,latexsym,dsfont}
\usepackage{graphicx}
\usepackage[numeric,initials,nobysame]{amsrefs}
\usepackage{upref,setspace}

\usepackage{enumerate}

\newenvironment{enumeratea}{\begin{enumerate}[\upshape (a)]}{\end{enumerate}}
\newenvironment{enumeraten}{\begin{enumerate}[\upshape 1.]}{\end{enumerate}}

\numberwithin{equation}{section}
\numberwithin{figure}{section}
\numberwithin{table}{section}

\usepackage{color}              % Need the color package
\usepackage{graphicx}
\usepackage[]{amsmath}
\usepackage{amssymb}
\usepackage{hyperref}
\usepackage{amsfonts}
\definecolor{MyDarkBlue}{rgb}{0,0.08,0.50}
\definecolor{BrickRed}{rgb}{0.65,0.08,0}
\usepackage[textsize=small]{todonotes}       % includes TO DO LIST
\usepackage{amsmath, amsthm, amssymb}
\newcommand{\BS}{{\bf{BSR}}}

\newcommand{\DDD}{{\mathbb{D}}}

\newcommand{\bfm}{\boldsymbol{M}}
\newcommand{\bfv}{\boldsymbol{v}}

\newcommand{\BF}{{\bf {BF}}}

\newcommand{\BM}{{\mbox{BM}}}

\newcommand{\RR}{\mathcal{R}}
\newcommand{\DD}{\mathcal{D}}
\newcommand{\CC}{\mathcal{C}}

\newcommand{\ompar}{\varpi}

\newcommand{\clf}{\mathcal{F}}
\newcommand{\clj}{\mathcal{J}}
\newcommand{\clc}{\mathcal{C}}

\newcommand{\clt}{\mathcal{T}}

\newcommand{\clp}{\mathcal{P}}
\newcommand{\clb}{\mathcal{B}}
\newcommand{\VV}{\mathcal{V}}
\newcommand{\NNN}{\mathbb{N}}
\newcommand{\ZZZ}{\mathbb{Z}}
\newcommand{\III}{\mathbb{I}}
\newcommand{\RRR}{\mathbb{R}}
\newcommand{\TTT}{\mathbb{T}}
\newcommand{\cals}{\mathcal{S}}

\newcommand{\XX}{\mathcal{X}}
\newcommand{\bfx}{{\bf{x}}}

\newcommand{\bfG}{{\bf{G}}}

\newcommand{\be}{\begin{equation}}
\newcommand{\ee}{\end{equation}}
\newcommand{\beq}{\begin{eqnarray*}}
\newcommand{\eeq}{\end{eqnarray*}}

\newcommand{\beqn}{\begin{eqnarray}}
\newcommand{\eeqn}{\end{eqnarray}}

\newcommand{\ba}{\begin{aligned}}
\newcommand{\ea}{\end{aligned}}
\newcommand{\bes}{\begin{equation*}}
\newcommand{\ees}{\end{equation*}}
\newtheorem{Lemma}{Lemma}[section]
\newtheorem{Proposition}[Lemma]{Proposition}

\newtheorem{Theorem}[Lemma]{Theorem}

\theoremstyle{definition}
\newtheorem{Remark}[Lemma]{Remark}

\newtheorem{Corollary}[Lemma]{Corollary}

\theoremstyle{definition}

\newcommand{\EE}{\mathcal{E}}

\newcommand{\prob}{\mathbb{P}}
\newcommand{\PP}{\mathcal{P}}
\newcommand{\E}{\mathbb{E}}

\newcommand{\FF}{\mathcal{F}}
\newcommand{\GG}{\mathcal{G}}
\newcommand{\NN}{\mathcal{N}}

\newcommand{\set}[1]{\left\{#1\right\}}

\newcommand{\Rbold}{{\mathbb{R}}}

\newcommand{\ind}[2]{1_{(e \in \pi(#1,#2))}}

\newcommand{\bfd}{{\bf d}}

\def\ind{{\rm 1\hspace{-0.90ex}1}}
\newcommand{\bfC}{\boldsymbol{C}}

\newcommand{\bfX}{\boldsymbol{X}}
\newcommand{\bfZ}{\boldsymbol{Z}}
\newcommand{\bfY}{\boldsymbol{Y}}

\newcommand{\ldown}{l^2_{\downarrow}}
\newcommand{\udown}{\mathbb{U}_{\downarrow}}

\newcommand{\calS}{\mathcal{S}}
\newcommand{\bars}{\bar{s}}
\newcommand{\sss}{\scriptscriptstyle}

\newcommand{\barx}{\bar{x}}
\def\1{{\mathchoice {1\mskip-4mu\mathrm l}      % Blackboard bold 1
{1\mskip-4mu\mathrm l}
{1\mskip-4.5mu\mathrm l} {1\mskip-5mu\mathrm l}}}

\newcommand {\convd}{\stackrel{d}{\longrightarrow}}
\newcommand {\convp}{\stackrel{\sss {\mathbb P}}{\longrightarrow}}

\newcommand{\erdos}{Erd\H{o}s-R\'enyi }

\newcommand{\spls}{{\bf spls}}
\newcommand{\ER}{{\bf ER}}
\newcommand{\vol}{{\bf vol}}

\numberwithin{equation}{section}
\begin{document}
	
	\title[Augmented multiplicative coalescent]{The augmented multiplicative coalescent and critical dynamic random graph models}

	\date{\today}
	\subjclass[2000]{Primary: 60C05, 05C80, 90B15}
	\keywords{bounded-size rules, surplus, critical random graphs, scaling window, multiplicative coalescent, entrance boundary, giant component, branching processes,
	inhomogeneous random graphs, differential equation method, dynamic random graph models.
	}

	\author[Bhamidi]{Shankar Bhamidi}
	\address{Department of Statistics and Operations Research, 304 Hanes Hall CB \#3260, University of North Carolina, Chapel Hill, NC 27599}
	\author[Budhiraja]{Amarjit Budhiraja}
	\author[Wang]{Xuan Wang}
	\email{bhamidi@email.unc.edu, budhiraj@email.unc.edu, wangxuan@email.unc.edu}

\begin{abstract}
	Random graph models with limited choice have been studied extensively with the goal of understanding the mechanism of the emergence of the giant component.
	 One of the standard models are the Achlioptas random graph processes on a fixed set of $n$ vertices.  Here at each step, one chooses two edges uniformly at random and then decides which one to add to the existing configuration according to some criterion. % In the context of famous rules such as the product rule (originally suggested by Bollobas as the rule most likely to delay the emergence of the giant), the processes seems to exhibit phase transitions fundamentally different from the classical Erd\H{o}s-R\'{e}nyi random graph evolution (\cite{achlioptas2009explosive}). 
	An important class of such rules are the \emph{bounded-size rules} where for a fixed $K\geq 1$, all components of size greater than $K$ are treated equally. While a great deal of work has gone into analyzing the subcritical and supercritical regimes, the nature of the critical scaling window,  the size and complexity (deviation from trees) of the  components in the critical regime and nature of the merging dynamics has not been well understood. In this work we study such questions for general bounded-size rules. 
	Our first main contribution is the construction of an extension of Aldous's standard multiplicative coalescent process which describes the asymptotic evolution of 
	the vector of sizes and surplus of all components.  We show that this process, referred to as the {\em standard augmented multiplicative coalescent} (AMC)
	is `nearly' Feller with a suitable topology on the state space. 
	Our second main result proves the convergence of suitably scaled component size and surplus vector, for any bounded-size rule, to the standard AMC.
	This result is new even for the classical \erdos setting. The key ingredients here are a precise analysis of the asymptotic behavior of various susceptibility functions near criticality and certain bounds from \cite{bsr-2012}, on the size of the largest component in the barely subcritical regime.

	 % and further there exists a {\bf standard augmented} process.  We then analyze detailed asymptotics of the susceptibility functions and rate of singularity as $t\to t_c$. % Using a careful analysis of the operator norm of an infinite dimensional mutitype branching process approximation first constructed in \cite{bhamidi-budhiraja-wang2011}, we prove asymptotically tight bounds on the size of the largest component in the barely subcritical regime.
	 % Using the singularity of the susceptibility functions as well as bounds on the largest component in the barely subcritical regime derived in \cite{bsr-2012} we are able to identify the critical scaling window and prove that  the largest components all scale like $\Theta(n^{2/3})$ and further properly rescaled, these sizes along with the surplus converge to the standard augmented multiplicative coalescent in the critical scaling window. 
	 % We finally provide conditions under which the diameter of the scaling window for truncations of general rules to bounded-size rules of order $K$ (including the truncations of the product rule) converge to $0$ as the level of truncation $K\to\infty$. 
\end{abstract}

\maketitle

\section{Introduction}
\label{sec:main-intro}
Profusion of  empirical data on real world networks has given impetus to research in mathematical models for such systems that explain the various observed statistics such as scale free degree distribution, small world properties and clustering. A range of mathematical models have been proposed both static as well as dynamic to understand the structural properties of such real world networks and their evolution over time. 
% This has resulted in an intense activity in a wide range of disciplines ranging from combinatorics to computer science in deriving rigorous results on various aspects of such models in the large network limit. In the context of mathematical probability, deep connections have arisen between these questions and classical probabilistic constructs such as Jagers-Nerman stable age distribution theory, coagulation processes and near critical branching processes. 
% \begin{center}
% 	\begin{figure}[htbp]
% 		\centering
% 			\includegraphics[scale=.35]{dsouza-graph.pdf}
% 		\caption{Taken from \cite{achlioptas2009explosive}, displaying the largest component $\calC_n^{\sss(1)}(r)/n$ for the three processes. Note that the critical time for the BF process is $t_c(BF)>1/2$. Also note the abrupt emergence of the giant component in the product rule.}
% 		\label{fig:ads-graph}
% 	\end{figure}
% 	\vspace{-.4in}
% \end{center}
One particular direction of significant research is focused on  understanding the effect of choice in the evolution of random network models (see \cite{spencer2007birth} and references therein).  
%In 2000, Achlioptas asked the following question: what is the effect of limited choice in delaying or accelerating the emergence of a giant component in network models?  
More precisely, suppose that at time $t=0$ we start with the empty configuration on $[n]:=\set{1,2,\ldots, n}$ vertices. At each discrete step $k=0,1,2,\ldots$,  we choose two edges $(e_1(k),e_2(k))$ uniformly at random amongst all ${n\choose 2}$ edges and decide whether the graph at instant $(k+1)$, denoted as  $\bfG_n(k+1)$, is $\bfG_n(k)\cup e_1(k)$ or $\bfG_n(k)\cup e_2(k)$ according to some pre-specified rule that takes into account suitable properties of the chosen edges with respect to the present configuration $\bfG_n(k)$. Speeding up time by a factor of $n$ and abusing notation, for $t\geq 0$ write,  $\bfG_n(t) = \bfG_n(\lfloor nt/2 \rfloor)$.  
Then the basic goal is to understand the effect of the rule governing the edge formations in the evolution of various characteristics of the network such as, the size of the largest component, the vector of sizes of all components, component complexities, etc.
Three prototypical examples to keep in mind are as follows: 
\\(a) {\bf Erd\H{o}s-R\'{e}nyi random graph:} At each stage include edge $e_1$ and ignore $e_2$. This results in the classical Erd\H{o}s-R\'{e}nyi random graph evolution.
For a component $\CC$, define $|\CC|$ for the size (number of vertices) of the component.  Well known results \cite{er-1,bollobas-rg-book} say that the critical time for the emergence of a giant component for this model is $1$, namely for $t<1$ the size of the largest component $|\CC^{\sss \bf ER}_1(t)| = O(\log{n})$ while for $t> 1$, the size of the largest component $|\CC^{\sss \bf ER}_1(t)| = \Theta(n)$.
Here  $O,\Theta$ are defined in the usual manner. More precisely,   given a sequence of random variables $\{\xi_n\}_{n\ge 1}$ and a function $f(n)$, we say $\xi_n = O (f)$ if there is a constant $C$ such that $\xi_n \le C f(n)$ with high probability (whp), and we say $\xi_n = \Omega(f)$ if there is a constant $C$ such that $\xi_n \ge Cf(n)$ whp. Say that $\xi_n = \Theta(f)$ if $\xi_n = O(f)$ {\bf and} $\xi_n = \Omega(f)$. 
\\
(b) {\bf Bohman-Frieze (BF) process: } This was the first rigorously analyzed example of a rule that delayed the emergence of the giant  component through limited choice \cite{bohman2001avoiding}. Here the rule is to use the first edge if it connects two singletons (vertices which have no connections at the present time), otherwise use the second edge. It has been shown (\cite{spencer2007birth, janson2010phase}) that there is a critical time $t_c^{\sss \BF} \approx 1.176$ when the largest component transitions from $O(\log{n})$ to $\Theta(n)$.\\
%the critical time for the Bohman Frieze model is $t \approx 1.176$.
% \footnote{Delete: Bohman and Frieze were able to show that for some time $t> 1/2$, $\CC_1(t) = o_P(n)$.}\\
(c) {\bf General bounded-size rules (BSR):} The BF process corresponds to a 
choice rule which treats all components with size greater than one in an identical fashion. 
It is a special case of the general family of models, referred to as bounded-size rules. Here one fixes $K\geq 1$ and then the rule for attachment is invariant on components of size greater than $K$. We postpone a precise description to Section \ref{sec:bsr}. General bounded-size rules were analyzed in \cite{spencer2007birth} where it was shown that there exists a (rule dependent) critical time $t_c$ such that for $t< t_c$, the largest component $|\CC^{\sss \BS}_1(t)| = O(\log{n})$ when $t< t_c$ while $|\CC^{\sss \BS}_1(t)| = \Theta(n)$ for $t> t_c$. 

Thus as time transitions from below to above $t_c$, a giant component (of the same order as the network) emerges. Motivated by recent results on the \erdos random graph components at criticality \cite{addario2009continuum} as well as general rules such as the (unbounded-size) product rule \cite{achlioptas2009explosive}, there has been a renewed interest in understanding the precise nature of the emergence of the giant component as well as structural properties of components near $t_c$ for classes of rules which incorporate limited choice in their evolution. 
%Write $\CC_i^{\sss \BS}(t)$ for the $i^{th}$ largest component in $\BS(t)$. 
Define the surplus or complexity of a component $\spls(\CC)$ as 
\begin{equation}
\label{eqn:surp-def}
	 \spls(\CC) = \mbox{number of edges} - (|\CC|-1).
\end{equation}
If a component were a tree, its surplus would be zero, thus this is a measure of the deviation of the component from a tree. Write $\CC_i(t)$ for the $i$-th largest component and $\xi_i(t) := \spls (\CC_i(t))$ for the surplus of the component $\CC_i(t)$.  For any of the rules above and a fixed $t\geq 0$, consider the vector of component sizes and associated surplus $(|\CC_i(t)|, \xi_i(t): i\geq 1)$.  In the context of the \erdos random graph process, precise fine-scale results are known about the nature of the emergence of the giant component as time $t$ transitions through the scaling window around $t_c^{\sss \bf ER}=1$. More precisely, for fixed $\lambda \in \Rbold$ write
\[\bar{\bfC}^{\sss \ER}(\lambda) := \left(\frac{1}{n^{2/3}}\left|\CC_i^{\sss \ER}\left(1+\frac{1}{n^{1/3}}\lambda\right)\right| :i\geq 1\right) 
\mbox{ and }  \bar{\bfY}^{\sss \ER}(\lambda) = \left (\xi_i^{\sss \ER}\left(1+\frac{1}{n^{1/3}}\lambda\right) :i\geq 1\right).
\]
Then Aldous in \cite{aldous1997brownian} showed:
\begin{enumeratea}
	\item The process $(\bar{\bfC}^{\sss \ER}(\lambda): -\infty < \lambda < \infty)$ converges to a Markov process called the standard multiplicative coalescent.
	\item For fixed $\lambda \in \Rbold$, the rescaled component sizes and the corresponding surplus $(\bar{\bfC}^{\sss \ER}(\lambda), \bar{\bfY}^{\sss \ER}(\lambda))$ converge jointly to a limiting random process described by excursions from zero of an inhomogeneous reflected Brownian motion $\hat{W}_\lambda$ and a counting process $\hat N_{\lambda}$ with intensity function $\hat{W}_\lambda(\cdot)$.  
\end{enumeratea}
We give a precise description of these results in Section \ref{sec:mult}. 
Obtaining similar results on critical asymptotics
for general inhomogeneous Markovian models such as the bounded-size rules requires new ideas. These rules lack a simple description for the dependence between edges making the direct use of the component exploration and associated random walk construction, the major workhorse in understanding random graph models at criticality (\cite{aldous1997brownian,bhamidi2009novel,bhamidi-hofstad-van,riordan2012phase,joseph2010component}), intractable. Thus it is nontrivial to identify the critical scaling window for such processes, let alone distributional asympototics for the component sizes and surplus. In the current work we develop a different machinery that allows us to identify the critical scaling window for all bounded-size rules.  Furthermore, denoting the suitably scaled component sizes and surplus processes as $(\bar{\bfC}^{\sss (n)}(\lambda), \bar{\bfY}^{\sss (n)}(\lambda))$, our
results describe the joint asymptotic behavior of
\[
	\left((\bar{\bfC}^{\sss (n)}(\lambda_1), \bar{\bfY}^{\sss (n)}(\lambda_1)), \ldots, (\bar{\bfC}^{\sss (n)}(\lambda_m), \bar{\bfY}^{\sss (n)}(\lambda_m))\right)\]
	for $-\infty < \lambda_1< \lambda_2< ... < \lambda_m < \infty$.
Starting point of our work is the construction of a new Markov process that is associated with the inhomogeneous reflected Brownian motion $\{\hat W_{\lambda}\}_{\lambda \in \Rbold}$
and the associated counting process $\{\hat N_{\lambda}\}_{\lambda \in \Rbold}$ which we refer to as the augmented multiplicative coalescent (AMC).  The main result of this work
shows that AMC is the characterizing process for a new universality class that includes, in addition to critically scaled components and surplus vectors for
\erdos graphs, analogous processes for all bounded-size rules. More precisely, our  contributions are as follows.
% 
% 
% 
% Since such exploratory techniques fails, we directly work at the \emph{process level}, first constructing a variant of Aldous's multiplicative coalescent that is able to track the dynamics of both component sizes and surplus. We then identify the critical scaling window as well as process level convergence for these functionals for {\bf all bounded-size rules} as the bounded-size rule process evolves through the critical scaling window. We now informally describe our contributions.  
\begin{enumeratea}
	\item In Theorem \ref{thm:smc-surplus}  we show the existence and ``near'' Feller property of a Markov process ${\boldsymbol Z}(\lambda)$, $-\infty < \lambda < \infty$, called the augmented multiplicative coalescent, which  tracks the evolution of both component sizes and surplus edges over the critical window. Aldous's standard multiplicative coalescent corresponds to the first coordinate of   this process. Identifying the correct state space and topology that is suitable for obtaining the Feller property for this process  turns out to be particularly delicate (see Remark \ref{prodnofell}). The (near) Feller property plays a key role in analyzing the joint distribution, at multiple time instants, of the component sizes and surplus for bounded-size rules in the critical scaling window. In proving the existence of the standard augmented process, a key role is played by Theorem \ref{theo:aldous-full-gene} which is a generalization of a result of Aldous for the component sizes of an inhomogeneous random graph, to a setting where one considers joint distributions of component sizes and surplus.
We believe that this result is of broader significance and can be used 	to analyze the distribution of surplus in the critical regime for various other  inhomogeneous random graph models, e.g. the rank-1 inhomogeneous random graphs (\cite{bollobas-riordan-janson}). 
	\item In Theorems \ref{thm:suscept-funct} and \ref{thm:suscept-limit} we analyze susceptibility functions (sums of moments of component sizes)  associated with a general bounded-size rule.
	 Spencer and Wormald \cite{spencer2007birth} showed that these susceptibility functions converge to limiting monotonically increasing deterministic functions which are finite only for $t< t_c$ and explode for $t> t_c$.  Theorem \ref{thm:suscept-funct}  uses a dynamic random graph process with immigration and attachment to show that these limiting functions for {\bf all bounded-size rules} have the same critical exponents as the \erdos random graph process.  Theorem \ref{thm:suscept-limit} shows that the susceptibility functions are close to their deterministic analogs in a strong sense even as $t\uparrow t_c$ when the limiting functions explode. 
	\item The analysis of the susceptibility functions gives rise to (rule dependent) constants $\alpha, \beta > 0$ which describe the nature of the explosion of the limiting susceptibility functions as $t\uparrow t_c$. For a given bounded-size rule we consider the rescaled process $\set{{\bar\bfZ}^{\sss (n)}(\lambda): -\infty<  \lambda < \infty} $ where ${\bar \bfZ}^{\sss (n)}(\lambda) = (\bar{\bfC}^{\sss (n)}(\lambda),\bar{\bfY}^{\sss (n)}(\lambda) )$ with $\bar{\bfC}^{\sss (n)}(\lambda)$ denoting the rescaled component sizes and $\bar{\bfY}^{\sss (n)}(\lambda)$ denoting the surplus of these components, namely
	\[\bar{\bfC}^{\sss (n)}(\lambda) := \left(\frac{\beta^{1/3}}{n^{2/3}}\left|\CC_i\left(t_c+\frac{\alpha\beta^{2/3}}{n^{1/3}}\lambda\right)\right| :i\geq 1\right) 
	\mbox{ and }  \bar{\bfY}^{\sss (n)}(\lambda) = \left (\xi_i\left(t_c+\frac{\alpha\beta^{2/3}}{n^{1/3}}\lambda\right) :i\geq 1\right).
	\]  
Using the Feller property proved in Theorem \ref{thm:smc-surplus}, Theorem \ref{theo:aldous-full-gene},  results on the susceptibility functions, and bounds on the maximal component in the barely subcritical regime from \cite{bsr-2012}, we show the  convergence 
of finite dimensional distribution of this process to that of the augmented multiplicative coalescent $\set{{\boldsymbol Z}(\lambda):-\infty< \lambda <\infty}$, namely
for any set of times $-\infty < \lambda_1< \lambda_2< ... < \lambda_m < \infty$,
\[
	\left(\bar{\bfZ}^{\sss (n)}(\lambda_1), \ldots, \bar{\bfZ}^{\sss (n)}(\lambda_m)\right) \convd \left(\bfZ(\lambda_1), \ldots, \bfZ(\lambda_m)\right).
\]
  The result in particular identifies the critical scaling window {\bf for all bounded-size rules} as well as the asymptotic joint distributions of component sizes and surplus for any fixed $\lambda$, implying that such rules belong to the same universality class as the \erdos random graph process. The  convergence for the joint distribution of the surplus and the component sizes for multiple time points $\lambda$ in the critical scaling window is new even in the context of the \erdos random graph process. 
\end{enumeratea}
% Without further ado, we start with a description of the probabilistic models analyzed in this paper and give a statement of our main results.  
% 
% % {\bf Organization:} 
The paper is organized as follows.  In Section \ref{sec:main-defn} we introduce some common notation, give a precise description of bounded-size rules, and give an informal 
description of the augmented multiplicative coalescent. Section \ref{sec:main-results} contains the statements of our main results.  Sections \ref{sec:main-constr-mc}
and \ref{sec:main-st-ag-mc-ex} are devoted to proving the existence and near Feller property of the AMC.  In particular Section \ref{sec:main-st-ag-mc-ex}
contains the proof of Theorem \ref{thm:smc-surplus}.  Section \ref{sec:main-bsr-susceptibility} studies the asymptotics of the susceptibility functions associated
with general bounded-size rules and proves Theorems \ref{thm:suscept-funct} and  \ref{thm:suscept-limit}.
Finally in Section \ref{sec:main-coupling}  we complete the proof of   Theorem \ref{thm:crit-regime}.

\section{Definitions and Notation}
\label{sec:main-defn}

\subsection{Notation}
We collect some common notation and conventions used in this work. 
%\subsubsection{ Graphs  and random graphs}
A graph $\bfG=\{\VV, \EE\}$ consists of a vertex set $\VV$ and an edge set $\EE$, where $\VV$ is a subset of some type space $\XX$. % and $\EE$ is a subset of all possible edges $\{ \{v_1,v_2\}:v_1 \ne v_2 \in \VV\}$.  An example of a type space is $[n]=\{1,2,...,n\}$. % Frequently we will assume $\XX$ to have additional structure, for example to be a measure space $(\XX,\TT,\mu)$.
For a finite set $A$ write $|A|$ for its cardinality. A graph $\bfG$ with no vertices and edges will be called a \textbf{null graph}. 
%$\bfG$ with vertex set $[n]$ is called an \textbf{ empty graph} if  $\EE=\emptyset$, and we write $\bfG = {\bf0}_n$.
% Given two graphs, $\bfG_i=\{\VV_i,\EE_i\}$ for $i=1,2$, $\bfG_1$ is said to be a \textbf{subgraph} of $\bfG_2$ if and only if $\VV_1 \subset \VV_2$ and
% $\EE_1 \subset \EE_2$ and we denote this as 
For graphs $\bfG_1, \bfG_2$, if $\bfG_1$ is a subgraph of $\bfG_2$ we shall write this as $\bfG_1 \subset \bfG_2$. % We write $\bfG_1 =\bfG_2$ if $\bfG_1 \le \bfG_2$ and $\bfG_1 \ge \bfG_2$. 
 % $\CC=\{\VV_0,\EE_0\}$ of a graph $\bfG=\{\VV, \EE\}$ is a subgraph which is connected (i.e. there is a finite path between any two vertices
% in $\CC$). 
 The number of vertices in a connected component $\CC$ of a graph $\bfG$ will be called the size of the component and will be denoted by $|\CC|$.
% and frequently we will denote the size and the component by the
%same symbol.
% Conventionally, we use $\CC=\{\VV_0,\EE_0\}$ denote a connected component of a graph. Given a non-negative function $\phi$ on $\XX$, the \textbf{size} (or \textbf{volume}) of the component is defined as $|\CC|=\sum_{v \in \VV_0} \phi(v)$. The `weight function' $\phi$ will always be clear in the context, and unless specified differently, $\phi \equiv 1$. When $\phi \equiv 1$ we use the terms  \textbf{volume}, \textbf{size}  and  \textbf{number of vertices} interchangeably.  We will use $\CC$ to denote both the components and their sizes.
Let $\GG$ be the set of all  possible graphs $(\VV, \EE)$ on a given type space $\XX$.  When $\VV$ is finite, we will
consider $\GG$ to be endowed with the discrete topology and the corresponding Borel sigma field and refer to a random element of $\GG$ as a random graph. 
% All random graphs in this work are given on a fixed probability space $(\Omega, \mathcal{F}, \mathbb{P})$ which will usually be
% suppressed in our proofs. 
%\subsubsection{Probability and  analysis}

We use $\convp$ and $\convd$ to denote convergence in probability and in distribution respectively.
All the unspecified limits are taken as $n \to \infty$.
Given a sequence of events $\{E_n\}_{n\ge 1}$, we say $E_n$  occurs with high probability (whp) if $\prob\{E_n\} \to 1$.
The notation $O,\Omega,\Theta$ was described in the Introduction.  Furthermore, for a sequence of random variables $\{\xi_n\}_{n\ge 1}$ and a function $f(n)$, we say  $\xi_n = o (f)$ if $\xi_n/f(n) \convp 0$.

 % Given a sequence of random variables $\{\xi_n\}_{n\ge 1}$ and a function $f(n)$, we say $\xi_n = O (f)$ if there is a constant $C$ such that $\xi_n \le C f(n)$ whp. Say that the sequence $\xi_n = \Theta(f)$ if, given any $\eps> 0$, there exist constants $C_1, C_2> 0$ such that 
 % \[\liminf_{n\to\infty} \prob(C_1 f(n)< \xi_n < C_2 f(n)) \geq 1-\eps\]
 %  We say $\xi_n = o (f)$ if $\xi_n/f(n) \convp 0$. 
%\footnote{Replace: We use $o(\cdot), O(\cdot)$ notation in the usual manner. }
% For functions $f, g: \mathbb{N} \to \mathbb{R}$, we write
%  $g=O(f)$ if for some $C \in (0, \infty)$, $\limsup g(n)/f(n) < C$ and $g=\Theta(f)$ if $g=O(f)$ and $f=O(g)$.  Given two sequences of random variables $\{\xi_n\}$ and $\{\zeta_n\}$, we say $\xi_n=O(\zeta_n)$ whp if there is a $C \in (0, \infty) $ such that $\xi_n < C \zeta_n$ whp, and write $\xi_n=\Theta(\zeta_n)$ whp if there exist  $0 < C_1 \le C_2 < \infty$ such that $C_1 \zeta_n<\xi_n<C_2 \zeta_n$ whp. 
% %Occasionally, when clear from the context, we suppress `whp' in the statements. 
% 
%  We also use the following little $o$ notation: For sequences of real numbers $\{g(n)\}$, $\{f(n)\}$, we write $g=o(f)$ if $\limsup|g(n)/f(n)|=0$. For a sequence of random variables $\xi_n$ and a sequence $\{f(n)\}$ of reals, we write ``$\xi_n=o_p(f)$'' if $\xi_n/f(n)$ converges to $0$ in probability.
% For a real measurable function $\psi$ on a measure space $(\XX,\TT,\mu)$, the norms $\|\psi\|_2$ and $\|\psi\|_\infty$ are defined in the usual way.
 
For a Polish space $S$, $\DD([0,T]:S)$ (resp. $\DD([0,\infty):S)$) denote the space of right continuous functions with left limits (RCLL) from $[0,T]$ (resp. $[0,\infty)$) equipped with the usual Skorohod topology.  For a RCLL function $f:[0,\infty) \to \RRR$,
we write $\Delta f(t) = f(t) - f(t-)$, $t >0$.
% We use $=_d$ to denote the equality of random elements in distribution. 
Suppose that $(S, \cals)$ is a measurable space and we are given a partial
ordering on $S$.
% Given two $S$ valued random variables $\xi, \tilde \xi$, we say a pair of $S$ valued random variables $\xi^*, \tilde \xi^*$ given on
% a common probability space define a coupling of $(\xi, \tilde \xi)$ if $\xi =_d\xi^*$ and $\tilde \xi =_d \tilde \xi^*$.
We say 
the $S$ valued random variable $\xi$ \textbf{ stochastically dominates} $\tilde \xi$, and write $\xi \ge_d \tilde \xi$ if there exists a coupling between the two random variables on a common probability space such that $\xi^* \ge \tilde \xi^*$ a.s., where $\xi^*=_d \xi$ and $\tilde \xi^* = \tilde \xi$.  For probability measures $\mu, \tilde \mu$ on $S$, we say $\mu$
stochastically dominates $\tilde \mu$, if $\xi \ge_d \tilde \xi$ where $\xi$ has distribution $\mu$ and $\tilde \xi$ has distribution $\tilde \mu$.
 % For two sequences of $S$ valued random elements $\xi_n$ and $\tilde \xi_n$, we say ``$\xi_n \le_d \tilde \xi_n$ whp.'' if there exist a coupling between $\xi_n$ and $\tilde \xi_n$ for each $n$ (denote as $\xi_n^*$ and $\tilde \xi_n^*$) such that $\xi_n^* \le \tilde \xi_n^*$ whp.
Two examples of $S$ relevant to this work are $\DD([0,T]: \mathbb{R})$ and $\DD([0,T]: \GG)$ with the natural associated partial ordering.
Given a metric space $S$, we denote by $\clb(S)$ the Borel $\sigma$ - field on $S$ and by
 $\mbox{BM}(S), C_b(S), \clp(S)$, the space of bounded (Borel) measurable functions, continuous and bounded function, and probability measures, on $S$,
respectively.  The set of nonnegative integers will be denoted by $\NNN_0$.

%We start with the definition of the dynamic random graph processes which motivated the construction of the limit object in this study. 

\subsection{Bounded-size rules (BSR)}
\label{sec:bsr}
We now define the general class of rules that will be analyzed in this paper. Much of the notation follows \cite{spencer2007birth} which provides a comprehensive analysis of the sub and supercritical regime. 
%Let us first describe the discrete time process which results in a random graph process on a set of $n$ vertices $[n]:=\set{1,2,\ldots, n}$. 
Fix $K \in \NNN$ and let $\Omega_0=\set{\ompar}$ and $\Omega_K =\set{1,2,\ldots, K, \ompar}$ for $K \ge 1$, where $\ompar$ will represent components of size greater than $K$. Given a graph $\bfG$ and a vertex $v\in\bfG$, write $\CC_v(\bfG)$ for the  component that contains $v$. Let 
\begin{equation}
\label{eqn:cv-def}
	c(v) = \left\{\begin{array}{ll}
		|\CC_v(\bfG)| & \mbox{ if }|\CC_v(\bfG)|\leq K\\
		\ompar & \mbox{ if } |\CC_v(\bfG)| > K.
	\end{array}\right.
\end{equation}  
%Call $c(v) $  the $K$-adjusted component size of $v$ in the graph $\bfG$.  
For a quadruple of vertices $v_1, v_2, v_3, v_4$, write $\vec{v} = (v_1,v_2,v_3,v_4)$ and let $c(\vec{v}) = (c(v_1), c(v_2), c(v_3), c(v_4))$. Fix $F\subseteq \Omega_K^4$. We now define the random graph process $\{\BS^{\sss(n)}(k)\}_{k\ge 0}$
on the vertex set $[n]$ evolving through a $F$-bounded-size rule ($F$-BSR)  as follows.
Define $\BS^{\sss(n)}(0) = \mathbf{0}_n$. Having defined $\BS^{\sss(n)}(k)$ for $k \ge 0$, $\BS^{\sss(n)}(k+1)$ is constructed as follows:
 Choose four vertices $(v_1(k), v_2(k), v_3(k), v_4(k))$ uniformly at random amongst all possible $n^4$ vertices uniformly at random and let \[\vec{v}_k = (v_1(k), v_2(k),v_3(k), v_4(k)).\] 
Denote the function $c(\vec{v})$ associated with $\BS^{\sss(n)}(k)$ as $c_k(\vec{v})$.
Define 
 \begin{equation}
	\label{eqn:f-rule-def}
 	\BS^{\sss(n)}(k+1) =\left\{ \begin{array}{ll}
 		\BS^{\sss(n)}(k)\cup (v_1(k), v_2(k)) & \quad \mbox{ if } c_k(\vec{v}_k) \in F\\
\BS^{\sss(n)}(k)\cup (v_3(k), v_4(k)) & \quad \mbox{ otherwise. }
 	\end{array} \right. 
 \end{equation}  
These rules are called bounded-size rules since they treat all components of size greater than $K$ identically. Concrete examples of such rules include
 Erd\H{o}s-R\'{e}nyi random graph (here $K=0$, $F=\Omega_0^4 =\set{\ompar,\ompar,\ompar,\ompar}$) and 
 Bohman-Frieze process (here  $K=1$, $F=\set{(1,1,\alpha, \beta): \alpha, \beta \in \Omega_1}$).
% \\(c) {\bf K-Truncated product rule:} Let $K\geq 1$ and let $F= \set{(i_1, i_2, i_3, i_4): i_l\leq K~ \forall ~l=1,2,3,4,~ i_1 i_2\leq i_3 i_4  }$. In words, add the first edge to the present configuration if the component sizes  of all the chosen vertices are smaller than $K$, else use the second edge. Note that in this case,  ``$K=\infty$'' corresponds to the original product rule. 

\noindent {\bf Continuous time formulation $\set{\BS^{\sss(n)}(t)}_{t\geq 0}$:} It will be more convenient to work in continuous time. % Let $\Ecal = \set{\set{v_1, v_2}: v_i\in [n], v_1\neq v_2}$ denote the space of all possible ${n\choose 2}$ edges on the $n$ vertices. For each pair of edges ${\bf e} =(e_1, e_2) \in \Ecal\times \Ecal$ construct a Poisson process $\PP_{{\bf e}}$ with rate $2/n^{3}$, independent across pairs of edges in $\Ecal\times \Ecal$. 
For every quadruple of vertices $\vec{v} = (v_1, v_2, v_3, v_4)\in [n]^4$, let $\PP_{\vec{v}}$ be a Poisson process with rate $\frac{1}{2n^3}$, independent between quadruples.  
The continuous time random graph process $\{\BS^{\sss(n)}(t)\}_{t\ge 0}$ is constructed recursively as follows.  
We denote the function $c(v)$ [resp. $c(\vec{v})$] associated with $\BS^{\sss(n)}(t-)$ as $c_{t-}(v)$ [resp. $c_{t-}(\vec{v})$].
Given $\BS^{\sss(n)}(t-)$, and that for some  
$\vec{v} \in [n]^4$, $\PP_{\vec{v}}$ has a point at the time instant $t$, we define
\begin{equation}
	\label{eqn:f-rule-def-cts}
 	\BS^{\sss(n)}(t) =\left\{ \begin{array}{ll}
 		\BS^{\sss(n)}(t-)\cup (v_1, v_2) & \quad \mbox{ if } c_{t-}(\vec{v}) \in F\\
\BS^{\sss(n)}(t-)\cup (v_3, v_4) & \quad \mbox{ otherwise. }
 	\end{array} \right. 
 \end{equation}
The rationale behind this scaling for the rate of the Poisson point process is that the total rate of adding edges is
\[\frac{n^4}{2n^3} = \frac{n}{2}.\]
Thus with this scaling, for the $F$-BSR rule corresponding to the Erd\H{o}s-R\'{e}nyi evolution, the giant component emerges at time $t=1$. To simplify notation, when there is no scope for confusion, we will suppress $n$ in the notation. For example, we write $\BS_t:= \BS^{\sss(n)}(t)$. 

Denote $\CC_i^{\sss (n)}(t)$ for the $i$-th largest component of $\BS_t$ at time $t$. 
%For the $i$-th largest component $\CC_i(t) \equiv \CC_i^{\sss (n)}(t)$ of $\BS_t$ at time $t$, let $\xi_i(t) :=\xi_i(\BS_t^{\sss(n)})$ denote the surplus of this component.
%For any graph $\bfG$, let $\CC_i(\bfG)$ denote the $i$-th largest component of $\bfG$. 
The  work of Spencer and Wormald (see \cite{spencer2007birth}) shows that for any given BSR, there exists a (model dependent) {\bf critical time} $t_c>0$ such that
 for $t < t_c$, $|\CC_1^{\sss (n)}(t)|= O(\log n)$ and
 for $t > t_c$, $|\CC_1^{\sss (n)}(t)| \sim f(t)n$ where $f(t) >0$.

%\todo[inline]{(Xuan) The definition of $s_2(t)$ should be re-organized.}
One of the key ingredients in the proof of the above result is an analysis of the susceptibility functions: For any given time $t$ and fixed $k\geq 0$ define the 
$k$-susceptibility function
\begin{equation}\label{eq:eq5.1}\calS_k^{\sss (n)}(t) \equiv \calS_k(t) := \sum_{i\geq 1} \left|\CC_i^{\sss (n)}(t)\right |^k. \end{equation}
%For ease of notation we have suppressed the dependence of $\calS_k$ on $n$. 
Then \cite{spencer2007birth} shows that 
for any bounded-size rule  and for every $k \ge 2$,  there exists a monotonically increasing
  function $s_k : [0,t_c) \to [0, \infty)$ satisfying $s_k(0)=1$ and $\lim_{t \uparrow t_c}s_k(t)=\infty$,  such that
\begin{equation}
\label{eqn:suscept-defn}
	\bar s_k(t) := \frac{\calS_k(t)}{n} \convp s_k(t) \quad \forall t\in [0,t_c).
\end{equation}
%To ease notation, we will typically write $\CC_i^{\sss(n)}(t) = \CC_i^{\sss (n)}(t) $. 

%\noindent{ \bf Complexity of components:}

Along with the size of the components, another key quantity of  interest is the complexity of components. 
Recall the definition of the surplus of a component from \eqref{eqn:surp-def}, and denote $\xi_i^{\sss (n)}(t) :=\spls(\CC_i^{\sss (n)}(t))$ for the surplus of the component $\CC_i^{\sss (n)}(t)$. We will be interested in the joint vector of ordered component sizes and corresponding surplus
\begin{equation*}
%\label{eqn:comp-vec-surp-def}
	((|\CC_i(t)|, \xi_i(t) ): i\geq 1).
\end{equation*} 
%\todo[inline]{Organize the definition of surplus for finite and infinite graphs.}
\subsection{Augmented Multiplicative coalescent}

\subsubsection{The multiplicative coalescent}
\label{sec:mult}
Let $l^2 = \{x = (x_1,x_2,\ldots):  \sum_i x_i^2< \infty\}$.  Then $l^2$ is a separable Hilbert space with the inner product $\langle x, y\rangle = \sum_{i=1}^{\infty}x_iy_i$,
$x=(x_i), y= (y_i) \in l^2$. Let
\begin{equation}
	\ldown = \{(x_1,x_2,\ldots): x_1\geq x_2\geq \cdots \geq 0, \sum_i x_i^2< \infty\}.
	\label{eqn:ldown}
\end{equation}
Then $\ldown$ is a closed subset of $l^2$ which we equip with the metric inherited from $l^2$. In \cite{aldous1997brownian} Aldous introduced a $\ldown$
valued continuous time Markov process, referred to as the {\em standard multiplicative coalescent}, that can be used to describe the asymptotic behavior
of suitably scaled component size vector in Erd\H{o}s-R\'{e}nyi random graph evolution, near criticality. Subsequently, similar results have been shown to hold 
for other random graph models (see \cite{bhamidi-budhiraja-wang2011,aldous2000random} and references therein).  We now give a brief description of this Markov process. \\

Fix $x=(x_i)_{i \in \NNN}$. Let $\{\xi_{i,j}, i,j \in \NNN\}$ be a collection of independent rate one Poisson processes. Given $t \ge 0$, consider the random graph with vertex set $\NNN$ in which there exist $\xi_{i,j}([0,t x_i x_j/2]) + \xi_{j,i}([0,t x_i x_j/2]) $ edges between $(i,j)$,
$1 \le i < j < \infty$, and there are $\xi_{i,i}([0,t x_i^2/2])$ self-loops with the vertex $i \in \NNN$. The volume of a component $\mathcal{C}$ of this graph is defined to be 
$$\vol(\CC) :=\sum_{i\in \mathcal{C}}x_i.$$
Let $X_i(x,t)$ be the volume of the $i$-th largest (by volume) component.  It can be shown that
$X(x,t) = (X_i(x,t), i \ge 1) \in \ldown$, a.s. (see Lemma 20 in \cite{aldous1997brownian}).
  Define
$$T_t: \BM(\ldown) \to \BM(\ldown),$$
as $T_tf(x) = \E(f(X(x,t)))$.  It is easily checked that $(T_t)_{t\ge 0}$ satisfies the semigroup property $T_{t+s} = T_tT_s$, $s,t\ge 0$, and \cite{aldous1997brownian}
shows that $(T_t)$ is Feller, i.e. $T_t(C_b(\ldown)) \subset C_b(\ldown)$ for all $t \ge 0$. The paper \cite{aldous1997brownian} also shows that the semigroup
$(T_t)$ along with an initial distribution $\mu \in \clp(\ldown)$ determines a Markov process with values in $\ldown$ and RCLL sample paths. Denoting by $P^{\mu}$
the probability distribution of this Markov process  on $\DD([0,\infty): \ldown)$, the Feller property says that $\mu \mapsto P^{\mu}$ is a continuous map.
One special choice of initial distribution for this Markov process is particularly relevant for the study of asymptotics
of random graph models.  We now describe this distribution. Let $\{W(t)\}_{t\ge 0}$ be a standard Brownian motion, and for a fixed $\lambda \in \Rbold$, define
\[W_\lambda(t) = W(t)+\lambda t-\frac{t^2}{2},\; t \ge 0.\]
Let $\hat{W}_{\lambda}$ denote the reflected version of $W_{\lambda}$, i.e., 
\begin{equation}
	\hat{W}_\lambda(t) = W_\lambda(t) - \min_{0\leq s\leq t} W_\lambda(s), \; t \ge 0.
	\label{eqn:inh-ref-bm}
\end{equation}
 An excursion of $\hat{W}_\lambda$ is an interval $(l,u) \subset [0,+\infty)$ such that $\hat W_\lambda(l)=\hat W_\lambda(u)=0$ and $\hat W_\lambda(t)>0$ for all $t \in (l,u)$. Define $u-l$ as the size of the excursion. Order the sizes of excursions of $\hat W_\lambda$ as 
\[X^*_1(\lambda)> X^*_2(\lambda)> X^*_3(\lambda)> \cdots\] and write $\bfX^*(\lambda) = (X^*_i(\lambda):i\geq 1).$
Then $\bfX^*(\lambda)$ defines a $\ldown$ valued random variable (see Lemma 25 in \cite{aldous1997brownian})
and let $\mu_{\lambda}$ be its probability distribution.  Using the Feller property and asymptotic connections with
certain non-uniform random graph models, the paper \cite{aldous1997brownian} shows that $\mu_{\lambda}T_t = \mu_{\lambda+t}$, for all $\lambda \in \Rbold$ and $t \ge 0$, where
for $\mu \in \clp(\ldown)$, $\mu T_t \in \clp(\ldown)$ is defined in the usual way: $\mu T_t(A) = \int T_t(1_A)(x) \mu(dx)$, $A \in \clb(\ldown)$.  Using this consistency
property one can determine a unique probability measure $\mu_{\mbox{\tiny{MC}}} \in \clp(\DD((-\infty, \infty): \ldown))$ such that, denoting the canonical coordinate process
on $\DD((-\infty, \infty): \ldown)$ by $\{\pi_t\}_{-\infty < t < \infty}$,
$$
\mu_{\mbox{\tiny{MC}}} \circ (\pi_{t+\cdot})^{-1} = P^{\mu_t}, \; \mbox{ for all } t \in \Rbold ,$$
where $\pi_{t+ \cdot}$ is the process $\{\pi_{t+s}\}_{s\ge 0}$.
The measure $\mu_{\mbox{\tiny{MC}}}$ is known as the {\em standard multiplicative coalescent}.
This measure plays a central role in characterizing asymptotic distribution of component size vectors in the critical window for random graph models \cite{aldous1997brownian,aldous2000random,bhamidi-budhiraja-wang2011}.

\subsubsection{The augmented multiplicative coalescent}
\label{sec:augmented-mc}
We will now augment the above construction and introduce a measure on a larger space that can be used to describe the joint asymptotic behavior of the component size vector and the associated surplus vector, for a broad family of random graph models.

Let $\NNN^{\infty} = \{y=(y_1, \cdots) : y_i \in \NNN, \mbox{ for all } i \ge 1 \}$ and define %  We consider $\NNN^{\infty}$ with the product topology
% which can be metrized so that $\NNN^{\infty}$ is a complete, separable metric space. 
 \[\udown = \{ (x_i, y_i)_{i\geq 1} \in \ldown \times \NNN^\infty : \sum_{i=1}^\infty x_i y_i < \infty \mbox{ and  }
y_m = 0 \mbox{ whenever } x_m =0, m \ge 1\}.\]
We will view $x_i$ as the volume of the $i$-th component and $y_i$ the surplus of the $i$-th component of a graph with vertex set $\NNN$.
Writing $x = (x_i)$ and $y = (y_i)$, we will sometimes denote $(x_i, y_i)$ as $z=(x,y)$.
 We equip  $\udown$ with the metric
\begin{equation}
	\bfd_{\sss \mathbb{U}}((x,y),(x',y'))=\left( \sum_{i=1}^\infty(x_i-x'_i)^2\right)^{1/2}+ \sum_{i=1}^\infty|x_iy_i-x'_iy'_i|. \label{eqn:distance}
\end{equation}
Note that one natural metric on $\udown$, denoted as $\bfd_1$, is the one obtained by replacing the second term in \eqref{eqn:distance}
with
$$\sum_{i=1}^\infty \frac{|y_i-y'_i|}{2^i}\wedge 1.$$
This metric corresponds to the 
topology on $\udown$ inherited from $\ell^2 \times \NNN^{\infty}$ taking the topology generated by the inner product  $\langle \cdot, \cdot\rangle$
 on $\ell^2$ and the product topology on $\NNN^{\infty}$; and then considering the product topology on $\ell^2 \times \NNN^{\infty}$.  Although the metric $\bfd_1$ in some
respects is simpler to work with, it is not a natural metric to consider for the study of the joint distribution of component size and surplus process.  
Another metric (which we denote as $\bfd_2$) that can be considered on $\udown$ corresponds to replacing the second term in \eqref{eqn:distance} with
$\bfd_{vt}(\mu_z, \mu_{z'})$, where $\mu_z = \sum_{i=1}^{\infty} \delta_{z_i}$, $\mu_{z'} = \sum_{i=1}^{\infty} \delta_{z'_i}$ and $\bfd_{vt}$ is the metric corresponding to the
vague topology on the space of $\NNN \cup \{\infty\}$ valued locally finite  measures on $(0,\infty) \times \NNN$.  However this metric as well is not suitable for our
purposes.
These points are further
discussed in Remark \ref{prodnofell}.
%\todo[inline]{Put In.}

Let $\udown^0 = \{(x_i, y_i)_{i\geq 1}\in \udown: \mbox{ if } x_k=x_m, k \le m, \mbox{ then } y_k \ge y_m\}$.
We now introduce the {\em augmented multiplicative coalescent} (AMC). This is a continuous time Markov process with values in $(\udown^0, \bfd_{\sss \mathbb{U}})$, whose dynamics can  heuristically be described as follows: The process jumps at any given time instant from state $(x,y) \in \udown^0$ to:
\begin{itemize}
	\item $(x^{ij}, y^{ij})$ at rate $x_ix_j$, $i \neq j$, where $(x^{ij}, y^{ij})$ is obtained by merging components $i$ and $j$ into a component with volume
	$x_i+x_j$ and surplus $y_i+y_j$ and reordering the coordinates to obtain an element in $\udown^0$.
	\item $(x, y^i)$  at rate $x_i^2/2$, $i \ge 1$, where $(x, y^i)$ is the state obtained by increasing the surplus in the $i$-th component from $y_i$ to $y_{i}+1$ and reordering  the coordinates (if needed) to obtain an element in $\udown^0$. 
\end{itemize}
Whenever $z=(x,y) \in \udown^0$ is such that $\sum_{i=1}^\infty x_i < \infty$, there is a well defined Markov process $\{\bfZ(z,\lambda)\}_{\lambda \ge 0}$ that corresponds to the above transition mechanism, starting at time $\lambda=0$ in the state $z$. In fact in Section \ref{sec:main-constr-mc} (see also Theorem \ref{thm:smc-surplus}) we will see, that there is a well defined Markov process $\{Z(z,\lambda)\}_{\lambda \ge 0}$ corresponding to the above dynamical description for any $z \in \udown^0$.
Define, for $\lambda \ge 0$, $\clt_\lambda: \BM(\udown^0) \to \BM(\udown^0)$ as 
$$(\clt_\lambda f)(z) = \E f(\bfZ(z,\lambda)).$$
As for Aldous' multiplicative coalescent, there is one particular family of distributions that plays a special role. Recall the reflected parabolic Brownian motion $\hat W_\lambda(t)$
from \eqref{eqn:inh-ref-bm}. Let $\PP$ be a  Poisson point process on $[0,\infty) \times [0,\infty)$ with intensity $\lambda_{\infty}^{\otimes2}$ (where $\lambda_{\infty}$ is the Lebesgue measure on $[0, \infty)$) independent of $\hat W_\lambda$. Let $(l_i,r_i)$ be the $i$-th largest excursion of $\hat W_{\lambda}$. Define
$$ X^*_i(\lambda) = r_i-l_i \mbox{ and } Y^*_i(\lambda) = |\PP \cap \{(t,z):0 \le z \le \hat W_\lambda(t), l_i \le t \le r_i  \}|. $$
Then $\bfZ^*(\lambda) = (\bfX^*(\lambda), \bfY^*(\lambda))$ is a.s. a $\udown^0$ valued random variable, where $\bfX^* = (X^*_i)_{i\ge 1}$ and
$\bfY^* = (Y^*_i)_{i\ge 1}$.  Let $\nu_{\lambda}$ be its probability distribution.
In Theorem \ref{thm:smc-surplus}   we will show that there there exists a $\udown^0$ valued stochastic process $(\bfZ(\lambda))_{ -\infty < \lambda < \infty}$
such that $\bfZ(\lambda)$ has probability distribution $\nu_{\lambda}$ for every $\lambda \in (-\infty, \infty)$ and for all $f \in \BM(\udown^0)$, and $ \lambda_1 < \lambda_2 $, we have
$$  \E[ f(\bfZ(\lambda_2))| \{\bfZ(\lambda)\}_{\lambda \le \lambda_1}] =(\clt_{\lambda_2-\lambda_1}f)( \bfZ(\lambda_1) ). $$
% 
% 
% 
% 
% 
% There exist a special version of the process $(\bfZ(\lambda))_{ -\infty < \lambda < \infty}$ such that 
% Thus $(\bfZ(\lambda))_{ -\infty < \lambda < \infty}$ follows the same dynamic as the AMC. We will characterize the marginal distribution $\bfZ(\lambda)$ as follows: 
% We will show in Proposition \ref{pro:existence-smc}, Section \ref{sec:main-constr-mc} that there exist a process $(\bfZ(\lambda))_{ -\infty < \lambda < \infty}$ such that 
% $$ \bfZ(\lambda) = (\bfX^*(\lambda), \bfY^*(\lambda)) = ((X^*_i(\lambda), Y^*_i(\lambda)): i\ge 1).  $$
The process $\bfZ$ will be referred to as the {\em standard augmented multiplicative coalescent}.
We will also show that $\{\clt_{\lambda}\}_{\lambda \ge 0}$ is a  semigroup, namely $\clt_{\lambda_1}\circ\clt_{\lambda_2} = \clt_{\lambda_1+\lambda_2}$,
for $\lambda_1, \lambda_2 \ge 0$ which is {\em nearly Feller}, in the sense made precise in the statement of Theorem \ref{thm:smc-surplus}.
It will be seen that this process plays a similar role in characterizing the asymptotic joint distributions of the component size and surplus vector in the
critical window as Aldous' standard multiplicative coalescent does in the study of asymptotics of the component size vector.

\section{Results}
\label{sec:main-results}
Our first result establishes the existence of the standard augmented coalescent process.
Let $\udown^1 = \{z=(x,y) \in \udown^0: \sum_i x_i = \infty\}$.

% \todo[inline]{Shall we break the following theorem into two, one is about the existence of the nearly Feller process and the other is about the standard version? - Xuan}
\begin{Theorem}
\label{thm:smc-surplus}  
There is a collection of maps  $\{\clt_t\}_{t\ge 0}$, $\clt_t: \BM(\udown^0) \to \BM(\udown^0)$ and a $\udown^0$ valued stochastic process 
 $\{{\boldsymbol Z}(\lambda)\}_{-\infty < \lambda < \infty}= \{(\bfX(\lambda),\bfY(\lambda))\}_{-\infty < \lambda < \infty}$ such that the following hold.
\begin{enumeraten}
	\item $\{\clt_t\}$ is a semigroup: $\clt_t \circ \clt_s = \clt_{t+s}$, $s, t \ge 0$.
	\item  $\{\clt_t\}$ is nearly Feller: For all $t > 0$, $f \in \BM(\udown^0)$ and $\{z_n\} \subset \udown^0$, such that
	$f$ is continuous at all points in $\udown^1$ and 
	$z_n \to z$ for some $z \in \udown^1$, we have $\clt_{t} f(z_n) \to  \clt_tf(z)$.
%	 For all $t > 0$, $f \in C_b(\udown^0)$ and $\{z_n\} \subset \udown^0$, such that	$z_n \to z$ for some $z \in \udown^1$, we have $\clt_{t} f(z_n) \to  \clt_tf(z)$.
	\item  The stochastic process $\{{\boldsymbol Z}(\lambda)\}$ satisfies the Markov property with semigroup $\{\clt_t\}$:
	For all $f \in \BM(\udown^0)$, and $ \lambda_1 < \lambda_2 $, we have
	$$  \E[ f(\bfZ(\lambda_2))| \{\bfZ(\lambda)\}_{\lambda \le \lambda_1}] =(\clt_{\lambda_2-\lambda_1}f)( \bfZ(\lambda_1) ). $$
	\item  Marginal distribution of $\bfZ(\lambda)$ is characterized through the parabolic reflected Brownian motion $\hat W_{\lambda}$:  For each $\lambda \in \RR$,
	 $\bfZ(\lambda)$ has the probability distribution $\nu_{\lambda}$.
	\item  If $f \in \BM(\udown^0)$ is such that $f(x,y)=g(x)$ for some $g \in \BM(\ldown)$, then
	$$  (\clt_t f)(z) = (T_t g)(x), \;\; \forall z=(x,y) \in \udown^0. $$
	Furthermore, $\{\bfX(\lambda) \}_{-\infty < \lambda < \infty}$ is Aldous's standard multiplicative coalescent.
	\end{enumeraten}

\end{Theorem}
A precise definition of $\clt_t$ can be found in Section \ref{sec:main-constr-mc}. We will refer
to $\bfZ$ as the standard augmented multiplicative coalescent.
Theorem  \ref{thm:smc-surplus} will be proved in Section \ref{sec:main-st-ag-mc-ex}.
% \todo[inline]{Cannot probably be on $\Zbold^\infty$. Have to think about the state space more carefully!}

The next few theorems deal with bounded-size rules. Throughout this work we fix $K \in \NNN_0$, $F \in \Omega_K^4$ and consider a $F$ -BSR as introduced in Section \ref{sec:bsr}.

%The first result gives asymptotically tight bounds (upto logarithmic factors), on the size of the largest component in the barely subcritical regime,  which are subsequently used in the asymptotic analysis of  susceptibility functions.  

% Recall  that $\CC_1^{\sss (n)}(t) \equiv \CC_1(t)$ denotes the size of the largest component in  $\BS_t^{\sss(n)}$.  
% \begin{Theorem}[{\bf Barely subcritical regime}]
% \label{thm:subcrit-reg}
% Fix  $\gamma \in (0,1/5)$.  Then there exists a   $B\in (0, \infty)$ such that,
% \[ \prob\left( |\CC_1(t)| \le B\frac{(\log n)^4}{(t_c-t)^2},~ \forall t <t_c-n^{-\gamma} \right) \to 1  \]
% as $n \to \infty$. 
% \end{Theorem}
The first two results consider the asymptotics of the susceptibility functions.  Recall the deterministic functions $s_k$ from \cite{spencer2007birth} introduced above \eqref{eqn:suscept-defn}.

\begin{Theorem}[{\bf Singularity of susceptibility}]
\label{thm:suscept-funct}
There exist  $\alpha, \beta \in  (0, \infty)$ such that 
	\begin{equation} s_2(t) = (1+O(t_c-t))\frac{\alpha}{t_c-t}, \quad s_3(t) = \beta [s_2(t)]^3(1+O(t_c-t)), \label{eq:eq803}\end{equation}
	as $t\uparrow t_c$. 
\end{Theorem}

\begin{Theorem}{\bf (Convergence of susceptibility functions)}
	\label{thm:suscept-limit}
	%Define $\bar s_2^{\sss(n)}(t) = \calS_2^{\sss(n)}(t)$ and $\bar s_3^{\sss(n)}(t) = \calS_3^{\sss(n)}(t)/n$. Then, 
	For every $\gamma \in (1/6, 1/5)$, 
	\begin{align}
		&\sup_{t\in [0, t_n]} \left| \frac{n^{1/3}}{\bar s_2(t)} -  \frac{n^{1/3}}{s_2(t)}\right|  \convp 0\\
		&\sup_{t\in [0, t_n]} \left| \frac{\bar s_3(t)}{(\bar s_2(t))^3} - \frac{ s_3(t)}{(s_2(t))^3} \right|  \convp 0	,	
	\end{align}

	where $t_n = t_c - n^{-\gamma}$.	
\end{Theorem}
% One immediate corollary of the above theorem is as follows, which has also been proved in \cite{spencer2007birth}.
% \begin{Corollary}
% 	For $t \in [0, t_c)$, we have $\calS_2(t)/n \convp s_2(t)$ and $\calS_3(t)/n \convp s_3(t)$.
% \end{Corollary}

%The next theorem makes precise how the susceptibility functions diverge as $t\to t_c$. 

We now state the main result which gives the asymptotic behavior in the critical scaling window as well as merging dynamics for all bounded-size rules. 
\begin{Theorem}[{\bf Bounded-size rules: Convergence at criticality}]
\label{thm:crit-regime}
	Let $\alpha, \beta \in (0, \infty)$ be as in Theorem \ref{thm:suscept-funct}.  For $\lambda \in \Rbold$ define
	\[\bar{\bfC}^{\sss (n)}(\lambda) := \left(\frac{\beta^{1/3}}{n^{2/3}}\left|\CC_i\left(t_c+\frac{\alpha\beta^{2/3}}{n^{1/3}}\lambda\right)\right| :i\geq 1\right) 
	\mbox{ and }  \bar{\bfY}^{\sss (n)}(\lambda) = \left (\xi_i\left(t_c+\frac{\alpha\beta^{2/3}}{n^{1/3}}\lambda\right) :i\geq 1\right).
	\]
	Then $\bar{\bfZ}^{\sss (n)} = (\bar{\bfC}^{\sss (n)}, \bar{\bfY}^{\sss (n)})$ is a stochastic process with sample paths in 
$\DD((-\infty,\infty):\udown)$ and
	% \[\bar \bfZ_n \stackrel{d}{\longrightarrow} \bfZ_{\tmc}\]
	% 	as $n\to\infty$, where $\bfZ_{\tmc}$ is the augmented standard multiplicative coalescent and  $\convd$ denotes finite dimensional convergence in the space $\DD((-\infty, \infty):\udown)$, i.e. 
	for any set of times $-\infty < \lambda_1< \lambda_2< ... < \lambda_m < \infty$
	\begin{equation}
		\label{eq:eq1241}\left(\bar{\bfZ}^{\sss (n)}(\lambda_1), \ldots, \bar{\bfZ}^{\sss (n)}(\lambda_m)\right) \convd \left(\bfZ(\lambda_1), \ldots, \bfZ(\lambda_m)\right)\end{equation}
	as $n\to\infty$, where $\bfZ$ is as in Theorem \ref{thm:smc-surplus}.    
	%This of course gives convergence of $\bar{\bfC}^{\sss (n)}(\lambda)$ for any fixed $\lambda\in \Rbold$. 
\end{Theorem}

\subsection{Background}
We now make some comments on the problem background and future directions.
\begin{enumeratea}
	\item {\bf Critical random graphs:} Starting with the early work of \erdos \cite{er-1,er-2}, there is now a large literature on understanding phase transitions in random graph models, see e.g. \cite{bollobas-rg-book,bollobas-riordan-janson,janson-luczak-bb} and the references therein. Proving and identifying phase transitions in dynamic random graph models such as the bounded-size rule requires a relatively new set of ideas and is much more recent \cite{spencer2007birth}.   The study of  \erdos random graph in the critical regime was carried out
in \cite{janson1994birth,aldous1997brownian}.  In particular   the  paper, \cite{aldous1997brownian} introduced the standard multiplicative coalescent to understand the merging dynamics of the \erdos random graph at criticality. The barely subcritical and supercritical regimes of the Bohman-Frieze process were studied respectively in \cite{bf-spencer-perkins-kang} and \cite{janson2010phase}, with the latter identifying the scaling exponents for the susceptibility functions  for the special case of the Bohman-Frieze (BF) process by using the special form of the differential equations for the BF process. The current work  extends this result to all bounded-size rules (Theorem \ref{thm:suscept-funct}) by viewing such processes as random graph processes with immigration and attachment (see Section \ref{sec:proof-conv-susceptibility}).  
	\item {\bf Unbounded-size rules: } One of the reasons for renewed interest in such models is the recent study of the product rule (\cite{achlioptas2009explosive}), where as before one chooses two edges at random and then uses the edge that minimizes the product of the components at the end points of the chosen edges. This is an example of an unbounded-size rule and simulations in \cite{achlioptas2009explosive} suggest different behavior at criticality as compared to the usual \erdos or BF random graph processes. There has been recent progress in rigorously understanding the continuity at the critical point \cite{riordan2011achlioptas} as well the subcritical regime \cite{riordan2012evolution}. 
Such unbounded rules can be regarded as formal limits of $K$-bounded-size rules analyzed in the current work, as $K \to \infty$.	
	It would be of great interest to identify and understand the critical scaling window of such processes.   
	\item {\bf Related open questions:} In the context of bounded-size rules our results suggest other related questions. In particular, there has been recent progress in understanding structural properties of the component sizes of the \erdos random graph at criticality, in particular see \cite{addario2009continuum,addario2009critical} which use information about the surplus and component sizes in \cite{aldous1997brownian} to prove that the components viewed as metric spaces, converge to random fractals closely related to the continuum random tree \cite{aldous1991-crt} with shortcuts due to surplus edges. Our results strongly suggest the components in any bounded-size rule at criticality belong to the same universality class. Proving this will require substantially new ideas. 
\end{enumeratea}

\subsection{Organization of the paper}
%\todo[inline]{Should we include Theorem \ref{theo:aldous-full-gene} in the Main Results section? }

%The next Section collects notation used in the rest of the paper. Then in Section \ref{sec:main-constr-mc} we dive into the proofs.  
The two main results in this paper are Theorems \ref{thm:smc-surplus} and  \ref{thm:crit-regime}.  
%We give an outline of the proofs of the main results.  
%\begin{itemize}
	 In Section \ref{sec:main-constr-mc}
	we introduce the semigroup
	$\{\clt_t\}_{t\ge 0}$ and, as a first step towards Theorem \ref{thm:smc-surplus},  establish in Theorem \ref{theo:welldef-feller} the existence of a $\udown^0$ valued Markov process associated with this semigroup, starting from an arbitrary
	initial value.
	  Then in Section \ref{sec:main-st-ag-mc-ex} we complete  the proof of Theorem \ref{thm:smc-surplus}.
	We then proceed to the analysis of bounded-size rules in Section \ref{sec:main-bsr-susceptibility} where we study the differential equation systems associated with the BSR process and prove Theorems \ref{thm:suscept-funct} and  \ref{thm:suscept-limit}.
Finally in Section \ref{sec:main-coupling}  we complete the proof of   Theorem \ref{thm:crit-regime}.
% \end{itemize}
 %Without further ado, let us dive into the proofs. 
% \subsection{Standard multiplicative coalescent with surplus}
% \label{sec:mult-coal-surp}

\section{The augmented  multiplicative coalescent}
\label{sec:main-constr-mc}
We begin by making precise the formal dynamics of the augmented multiplicative coalescent process given in Section \ref{sec:augmented-mc}.
 Fix $(x,y)\in \udown^0$. Let $\{\xi_{i,j}\}_{i,j \in \NNN}$ be a collection of i.i.d. rate one Poisson  processes. Let $\bfG(z,t)$, where $z=(x,y)$, be the random graph on vertex set $\mathbb{N}$ given as follows:\\
(I) For $i \in \NNN$, there are $y_i$ self-loops to the vertex $i$. \\
(II) For $i<j \in \NNN$, there are $\xi_{i,j}([0,t x_i x_j/2]) + \xi_{j,i}([0,t x_i x_j/2])$ edges between vertices $i$ and $j$. For $i \in \NNN$, there are $\xi_{i,i}([0,t x_i^2/2])$ self-loops to the vertex $i$. \\
%(II) For $i<j \in \NNN$, there are $\left|\xi_{i,j}  ([0,tx_ix_j/2])\right|$ edges between vertices $i$ and $j$. Note that when $i \neq j$, the edges between the two vertices come from both $\xi_{i,j}$ and $\xi_{j,i}$ and when $i=j$, the edges correspond to self-loops.\\
Let $\FF^x_t = \sigma\{ \xi_{i,j}([0,sx_ix_j/2]): 0 \le s  \le t, \; i,j \in \NNN\}$, $t\ge 0$.

Recall the volume of a component $\CC$ is defined to be $\vol(\CC) = \sum_{i \in \CC} x_i$. The surplus of a finite connected graph was defined in \eqref{eqn:surp-def}. % Thus for a finite graph consisting of $c$ components and contains finite number of vertices ($|\VV|$) and edges ($|\EE|$), then its surplus is $|\EE| - |\VV| + c $. 
For infinite graphs the definition requires some care.
%The mathematical definition of surplus for an infinite connected graph with infinitely many vertices is trickier. 
We define the surplus for a connected graph $\bfG$ with vertex set a subset of $\NNN$ as
$$ \spls(\bfG) := \lim_{k \to \infty} \spls(\bfG^{\sss [k]}),  $$
where $\bfG^{\sss [k]}$ is the {\bf induced subgraph} that has the vertex set $[k]$ (the subgraph with vertex set $[k]$ and all edges between vertices in $[k]$ that
are present in $\bfG$). It is easy to check that this definition of surplus does not depend on the labeling of the vertices. Further note that the surplus of a connected graph might be infinite with this definition.\\
% \todo[inline]{This needs to be properly defined.  What is an induced subgraph?}

 Thus letting $\tilde \CC_i(t)$ be the $i$-th largest component (in volume) in $\bfG(z,t)$,
define $X_i(z,t) : = \vol(\tilde \CC_i(t))$ and $Y_i(z,t):= \spls(\tilde \CC_i(t))$ to  be the {\bf volume} and the {\bf surplus} of the $i$-th largest component at time $t$. In case two components have the same volume, the ordering of $(\tilde \CC_i(t): i \ge 1)$ is taken to be such that $Y_m(z,t) \ge  Y_k(z,t)$ whenever $m \le k$ and $X_m(z,t) = X_k(z,t)$.

Let $\bfX^z(t) := (X_i(z,t): i \ge 1 )$ and $\bfY^z(t) := (Y_i(z,t): i \ge 1)$.
The paper \cite{aldous1997brownian} shows that $\bfX^z(t) \in \ldown$ a.s. for all $t\ge 0$.  The following result shows that
$\bfZ^z(t) = (\bfX^z(t), \bfY^z(t)) \in \udown^0$ a.s., for all $t$.
\begin{Theorem}
	\label{theo:welldef-feller}
	Fix $ z=(x,y)\in \udown^0$ and let $(\bfX^z(t),\bfY^z(t))_{t \ge 0}$ be the stochastic process described  above, then
	for any fixed $t \ge 0$, $(\bfX^z(t),\bfY^z(t)) \in \udown^0$. 
\end{Theorem}	
The above theorem will be proved in Section \ref{sec:well-defined}. 
For $t\ge 0$, define $\clt_t: \BM(\udown^0) \to \BM(\udown^0)$	as
$$\clt_t f(z) = \E f(\bfZ^z(t)),\; z \in \udown^0,\; f \in \BM(\udown^0).$$
The following result shows that $\{\clt_t\}$ is a semigroup that is (nearly) Feller.
\begin{Theorem}	
	\label{theo:welldef-fellerb}
	For $t, s \ge 0$, $\clt_t\circ \clt_s = \clt_{t+s}$.  For all $t > 0$, $f \in \BM(\udown^0)$ and $\{z_n\} \subset \udown^0$, such that
	$f$ is continuous at all points in $\udown^1$ and 
	$z_n \to z$ for some $z \in \udown^1$, we have $\clt_{t} f(z_n) \to  \clt_tf(z)$.
% et $z^{\sss (n)} = (x^{\sss (n)}, y^{\sss (n)}) \to z=(x,y)$ in $\udown^0$. Then, for any $t\ge 0$, $\bfZ^{\sss (n)}(t)=\bfZ(z^{\sss (n)},t)$ converges in distribution, in $\udown$, to  $(\bfZ^z(t))$. 
\end{Theorem}	
 The above theorem will be proved in  Section \ref{sec:feller-property}.   
Throughout we will assume, without loss of generality, that for all $z \in \udown^0$, $\bfZ^z$ is constructed using the same set of Poisson processes $\{\xi_{i,j}\}$.
This coupling of $\bfZ^z$ for different values of $z$ will not be noted explicitly in the statement of various results.

% Firstly, we introduce two tools that will be used repeatedly.\\
% 
% The first tool is a coupling technique so called {\bf $\xi$-coupling}. Recall that given the initial vector $(x,y)\in \udown$, we are using the point processes $\{\xi_i, \xi_{i,j} \}$ to construct the process. For two processes with different initial states $(x,y), (\tilde x, \tilde y) \in \udown$, we can use the same set of point processes to define the augmented multiplicative coalescent processes, and we call such a coupling $\xi$-coupling.\\

We begin with the following elementary lemma.
\begin{Lemma}
	\label{lemma:borel-cantelli}
	Let $\{\FF_m\}_{m \in \NNN_0}$ be a filtration given on some probability space.\\
	(i) Let $\{Z_m\}_{m \ge 0}$ be a $\{\FF_m\}$ adapted sequence of nondecreasing random variables such that $Z_0=0$.  Let 
	 $\lim_{m \to \infty}Z_m=Z_\infty$. Suppose there exists a nonnegative random variable $U$ such that $U<\infty$ a.s. and
	 $ \sum_{m=1}^\infty \E[Z_m-Z_{m-1} | \FF_{m-1} ] \le U.$ Then for any $\epsilon \in (0,1)$,
	$$ \prob\{ Z_\infty > \epsilon \} \le \frac{1+\epsilon}{\epsilon} \E[U \wedge 1]. $$
	(ii) Let $\{A_m\}$ be a sequence of events such that $A_m \in \FF_m$. Suppose there exists a random variable $U<\infty$ a.s. such that $\sum_{m=1}^\infty \E[ {\ind}_{A_m} | \FF_{m-1}] \le U$. Then $\prob\{ A_m \mbox{ i.o.} \}=0$.  Furthermore,
	$$ \prob\{ \cup_{m=1}^\infty A_m \} \le 2 \E[U \wedge 1]. $$
\end{Lemma}
{\bf Proof:} (i) Define $B_m := \sum_{i=1}^m \E[Z_i-Z_{i-1}| \FF_{i-1}]$. Note that $B_m$ is nondecreasing and  $\FF_{m-1}$-measurable.
 Define $\tau=\inf\{ l: B_{l+1} > 1 \}$
where the infimum over an empty set is taken to be $\infty$.
 Since $B_m$ is predictable,  $\tau$ is a stopping time and, for all $m$, $B_{m\wedge \tau} \le 1$.
  Let $B_{\infty} = \lim_{m\to \infty} B_m$.   Since $Z_{m \wedge \tau}-B_{m\wedge\tau}$ is a martingale,
$$\E[Z_\tau] = \lim_{m \to \infty} \E[Z_{m \wedge \tau}] =  \lim_{m \to \infty} \E[B_{m \wedge \tau}] \le \lim_{m \to \infty} \E[B_{m} \wedge 1]=\E[B_\infty \wedge 1]. $$
Thus $$ \prob\{ Z_\infty > \epsilon\} \le \prob\{ \tau < \infty\}+ \frac{1}{\epsilon} \E[B_\infty \wedge 1] = \prob\{ B_\infty > 1\}+ \frac{1}{\epsilon} \E[B_\infty \wedge 1] \le \frac{1+\epsilon}{\epsilon} \E[U \wedge 1].$$
(ii) The first statement is immediate from the  Borel-Cantelli lemma (cf. Theorem 5.3.7 \cite{durrett-book}). For the second statement note that for any $\epsilon \in (0, 1)$, we have $ { \cup_{m=1}^\infty A_m } = \{\sum_{m=1}^\infty {\ind}_{A_m} > \epsilon\} $.  Now applying   part (i) to $Z_m = \sum_{k=1}^m {\ind}_{A_k}$ and taking $\epsilon \to 1$  yields the desired result. \qed\\

Next, we present a result from \cite{aldous1997brownian} that will be used here.  We begin with some notation. For  $x \in \ldown$, we write $x^{\sss [k]}=(x_1,...,x_k, 0, 0, ...)$ for the $k$-truncated version of $x$.
Similarly, for a sequence $x^{\sss (n)}=(x_1^{\sss (n)},x_2^{\sss (n)},...)$ of elements in $\ldown$, $x^{\sss (n)[k]}$ is the $k$-truncation of $x^{\sss (n)}$.
For $z = (x,y), z^{\sss (n)} = (x^{\sss (n)}, y^{\sss (n)}) \in \udown^0$ $z^{\sss [k]}, y^{\sss [k]}, z^{\sss (n)[k]}, y^{\sss (n)[k]}$ are defined similarly.\\
Recall the construction of $\bfG(z,t)$ described in items (I) and (II) at the beginning of the section.
We will distinguish the surplus created in $\tilde \clc_i(t)$ by the action in item (I) and that in item (II). The former  will be referred to as the type I surplus and denoted as 
$\tilde Y_i(z,t)$ while the latter will be referred to as the type II surplus and denoted as $\hat Y_i(z,t) \equiv \hat Y_i(x,t)$.
More precisely,
$$ \tilde Y_i(z,t) = \sum_{j \in \tilde \CC_i(t)} y_j \; \mbox{ and } \;\hat Y_i(z,t)=Y_i(z,t)-\tilde Y_i(z,t). $$
 Also define
$$ \tilde R(z,t) := \sum_{i=1}^\infty X_i(z,t) \tilde Y_i(z,t), \;\; \hat R(x,t) \equiv \hat R(z,t) :=\sum_{i=1}^\infty X_i(z,t) \hat Y_i(z,t) $$
and 
$$R(z,t):=\sum_{i=1}^\infty X_i(z,t)Y_i(z,t), \;\; S(x,t) \equiv S(z,t) := \sum_{i=1}^\infty (X_i(x,t))^2.$$

The following properties of $S$ and $\bfX$ have been established in  \cite{aldous1997brownian}.
\begin{Theorem} 
	\label{theo:aldous-mc-s} [Aldous\cite{aldous1997brownian}]
	(i) For every $x \in \ldown$ and $t\ge 0$, we have $S(x,t)<\infty $ a.s. and $S(x^{\sss [k]}, t) \uparrow S(x,t)$ as $k \to \infty$.\\
	(ii) If $x^{\sss (n)} \to x$ in $\ldown$, then  $\bfX(x^{\sss (n)},t) \convp \bfX(x,t)$ in $\ldown$, as $n \to \infty$. In particular, $\{S(x^{\sss (n)},t)\}_{n \ge 1}$ is tight.
\end{Theorem}
%Now we are ready to start the proof.\\

\subsection{Existence of the augmented MC}
\label{sec:well-defined}
This section proves Theorem \ref{theo:welldef-feller}.
We begin by  considering the type I surplus.
\begin{Proposition}
	\label{lemma:211} For any $t \ge 0$ and $z \in \udown^0$, $ \tilde R(z,t) = \sum_{i=1}^\infty X_i(z,t) \tilde Y_i(z,t) < \infty $ a.s.
\end{Proposition}
Proof of Proposition \ref{lemma:211} is given below Lemma \ref{lemma:765}.  The basic idea is to  bound the truncated version $\tilde R^{\sss [k]}=\tilde R(z^{\sss [k]},t)$ using a martingale argument, and then let $k \to \infty$. The truncation error is  controlled  using Lemma \ref{lemma:248} below and a suitable supermartingale is
constructed  in Lemma \ref{lemma:765}. 

\begin{Lemma} 
	\label{lemma:248}
For every $z \in \udown^0$ and $t \ge 0$, as $k\to \infty$, $\tilde R(z^{\sss [k]},t) \to \tilde R(z,t) \le \infty$ a.s. 
\end{Lemma}
{\bf Proof:} Fix $t\ge 0$. Denote by $A_{ij}$ [resp. $A_{ij}^{\sss [k]}$] the event that there exists a path from $i$ to $j$ in $\bfG(z,t)$ [resp. $\bfG(z^{\sss [k]},t)$], with the
convention that $\prob \set{A_{ii}} = \prob \{ A_{ii}^{\sss [k]} \}= 1$.  Let
$$f_i = \sum_{j=1}^\infty y_j {\ind}_{A_{ij}}, \;\; f_i^{\sss [k]} =\sum_{j=1}^k y_j {\ind}_{A_{ij}^{\sss [k]}}.$$
Then
$$   \tilde R (z,t) = \sum_{i=1}^\infty f_i x_i, \;\; \tilde R(z^{\sss [k]},t)= \sum_{i=1}^\infty f_i^{\sss [k]}x_i.$$
Since $A_{ij}^{\sss [k]} \uparrow A_{ij}$, we have $ f_i^{\sss [k]} \uparrow f_i$.  The result now follows from an application of monotone convergence theorem. \qed\\

\begin{Lemma}\label{lemma:765}
	Suppose that $z = (x,y) = z^{\sss [k]}$ for some $k \ge 1$ and that $\sum_j y_j \neq 0$.  Then
	$$ A_t=A(z,t)= \log \tilde R(z,t) - \int_0^t S(z,u)du$$
 is a supermartingale with respect to the  filtration $\FF_t^x = \sigma\{ \xi_{i,j}([0,sx_ix_j/2]); 0 \le s \le t, \; i,j \in \NNN\}$.
\end{Lemma}
{\bf Proof:} 
From the construction of $\bfZ(z, \cdot)$ we see that $ \tilde R(z,t)$ is a pure jump, nondecreasing process that at any time instant $t$, jumps at rate
$X_i(z,t-)X_j(z,t-)$, $1 \le i < j \le k$, with jump sizes $B_{ij}(t-) = X_i(z,t-)\tilde Y_j(z,t-)+X_j(z,t-) \tilde Y_i(z,t-)$.  Consequently
$\log \tilde R(z,t)$ jumps at the same rate, with corresponding jump size 
$\log ( 1 + \frac{B_{ij}(t-)}{\tilde R(z,t-)})$.  From this and elementary properties of Poisson processes it follows that
$$
\log \tilde R(z,t) = \log \tilde R(z,0) + \sum_{1 \le i < j \le k} \int_0^t \log \left (1 + \frac{B_{ij}(u)}{\tilde R(z,u)}\right ) X_i(z,u)X_j(z,u) du  + M(t),$$
where $M$ is a $\FF_t^x $ martingale.  Consequently, for $0 \le s < t < \infty$
\begin{equation}
	\label{eq:ab619}
\log \tilde R(z,t)  - \log \tilde R(z,s)= \sum_{1 \le i < j \le k} \int_s^t \log \left (1 + \frac{B_{ij}(u)}{\tilde R(z,u)}\right ) X_i(z,u)X_j(z,u) du  + M(t)	- M(s).
	\end{equation}
	Next note that, for $u \ge 0$
	\begin{align*}
	&\sum_{1 \le i < j \le k}  \log \left (1 + \frac{B_{ij}(u)}{\tilde R(z,u)}\right ) X_i(z,u)X_j(z,u)\\
	\le & 	\sum_{1 \le i < j \le k} \frac{B_{ij}(u)}{\tilde R(z,u)} X_i(z,u)X_j(z,u)\\
	= & \sum_{1 \le i < j \le k} \frac{X_i(z,u)\tilde Y_j(z,u)+X_j(z,u) \tilde Y_i(z,u)}{\tilde R(z,u)} X_i(z,u)X_j(z,u)\\
	\le & S(z,u).
	\end{align*}
	Using this observation in \eqref{eq:ab619} we now have
	$$
	\E \left [\log \tilde R(z,t)  - \log \tilde R(z,s) \mid \FF_s^x \right] \le \E \left [ \int_s^t S(z,u) du \mid \FF_s^x \right].$$
	The result follows. 
 \qed\\

{\bf  Proof of Proposition \ref{lemma:211}:} Fix $z = (x,y) \in \udown^0$.  The result is trivially true if $\sum_i y_i = 0$.  Assume now that $\sum_i y_i \neq 0$.
For $k \ge 1$ and $a \in (0, \infty)$, define
$T_{a}^{\sss [k]} =\inf\{ s \ge 0: S(z^{\sss [k]},s) \ge a \}.$ 
Fix $k \ge 1$ and assume without loss of generality that $\sum_{i=1}^k y_i > 0$.
Write $R^{\sss [k]}(t)=R(z^{\sss [k]},t)$, and $A^{\sss [k]}(t)=A(z^{\sss [k]},t)$ where $A$ is as  in Lemma \ref{lemma:765}. 
%Note that we can use a common filtration $\FF^x_t = \sigma\{ \xi_{i,j}([0,tx_ix_j/2]); i,j \in \NNN\}$ for different $k$. 
From the supermartingale property $\E[A^{\sss [k]} (T_{a}^{\sss [k]} \wedge t)] \le \E[A^{\sss [k]}(0)]= \log \tilde R^{\sss [k]}(0)$.  Therefore 
$$ \E\left[ \log \frac {\tilde R^{\sss [k]}(T_{a}^{\sss [k]} \wedge t)}{\tilde R^{\sss [k]}(0)}\right] \le \E\left[ \int_0^{T_{a}^{\sss [k]} \wedge t} S(z^{\sss [k]}, u) du \right] \le ta. $$
Thus
\begin{align*}
	\prob\{ \tilde R^{\sss [k]}(t) > m \} 
	\le \prob\{ \tilde R^{\sss [k]}(t) > m, T_a^{\sss [k]} > t \} + \prob\{  T_a^{\sss [k]} \le t  \}
	\le \frac{ t a} { \log m -\log \tilde R^{\sss [k]}(0)} + \prob \{ T_a^{\sss [k]} \le t\}.
\end{align*}
By Lemma \ref{lemma:248}, $\tilde R^{\sss [k]}(t) \to \tilde R(z,t)$, and by Theorem \ref{theo:aldous-mc-s} (i), $S(z^{\sss [k]},t) \to S(z,t)$ when $k \to \infty$.
Therefore letting $k \to \infty$ on both sides of the above inequality, we have 
\be
 \prob\{ \tilde R(z,t) > m \} \le \frac{ t a} { \log m -\log \tilde R(z,0)} + \prob \{ S(z,t) \ge a\}. \label{eqn:282}
\ee
The result now follows on first letting $m \to \infty$ and then  $a \to \infty$ in the above inequality. \qed \\

The following result is an immediate consequence of the estimate in \eqref{eqn:282} and Theorem \ref{theo:aldous-mc-s}(ii).
\begin{Corollary}\label{corr:ab1157}
	If $z^{\sss (n)} \to z$ in $\udown^0$, then for every $t\ge 0$, $\{\tilde R(z^{\sss (n)},t)\}_{ n \ge 1}$ is tight.
\end{Corollary}

Next we consider the type II surplus.  Let, for $x \in \ldown$
$$\GG_t^x := \sigma\{ \{\xi_{i,j}([0,sx_ix_j/2]) = 0\}:\; 0 \le s \le t,\; i,j \in \NNN \}. $$
The $\sigma$-field $\GG_t^x$ records the information whether or not $i$ and $j$ are in the same component at time $s$, for all $i,j$ and for all $s \le t$.
In particular, components $\{\tilde \clc_i(s), i \ge 1, s \le t\}$ can be determined from the information in $\GG_t^x$ and consequently, 
$\bfX(x,t)$ is $\GG_t^x $ measurable.
\begin{Lemma} 
	\label{lemma:786}
	(i) Fix $x \in \ldown$ and $t \ge 0$.  Then $\hat R(x,t) < \infty$ a.s.\\
	(ii) Let $x^{\sss (n)} \to x$ in $\ldown$. Then the sequence $\{ \hat R(x^{\sss (n)},t) \}_{n \ge 1}$ is tight.
\end{Lemma}
{\bf Proof:} Note that  (i) is an immediate consequence of (ii).  Consider now (ii).
For fixed $x \in \ldown$ and $t \ge 0$, let $\hat \mu_i(x,t)$ denote the probability law of $\hat Y_i(x,t)$, conditioned on $\GG_t^x$.
Then, for a.e. $\omega$, $\hat \mu_i(x,t)$ is a Poisson random variable with parameter
$$
\int_0^t  \sum_{j=1}^\infty (\sum_{k,k'\in \tilde \CC_j(s)} \frac{1}{2}x_k x_{k'} ) {\ind}_{\{\tilde \CC_j(s) \subset \tilde \CC_i(t)\}} ds
= \int_0^t \frac{1}{2} \sum_{j=1}^\infty (X_j(x,s))^2 {\ind}_{\{\tilde \CC_j(s) \subset \tilde \CC_i(t)\}} ds
\le  \frac{t}{2}  (X_i(x,t))^2,
$$
where the last inequality is a consequence of the inequality $\sum_{j:\tilde \CC_j(s) \subset \tilde \CC_i(t) }(X_j(x,s))^2 \le (X_i(x,t))^2$.
Therefore $\hat \mu_i(x,t)\le_d \hat \nu_i(x,t)$, a.s., where $\hat \nu_i(x,t)$ is a random probability measure on $\NNN$ such that for a.e. $\omega$,
$\hat \nu_i(x,t)$ is Poisson with parameter $\frac{t}{2} (X_i(x,t, \omega))^2 $.

A similar argument shows that the conditional distribution of $\sum_{i=1}^{\infty} \hat Y_i(x,t)$, given $\GG_t^x$ is a.s. stochastically dominated by
a random measure on $\NNN$ that, for a.e. $\omega$ has a Poisson distribution with parameter
$\sum_{i=1}^{\infty}\frac{t}{2} (X_i(x,t, \omega))^2  = \frac{t}{2} S(x,t)$.
Also, if $x^{\sss (n)}$ is a sequence converging to $x$ in $\ldown$, we have that for each $n$, the conditional distribution of $\sum_{i=1}^{\infty} \hat Y_i(x^{\sss (n)},t)$, 
given $\GG_t^{x^{\sss (n)}}$ is a.s. stochastically dominated by a Poisson random variable with parameter $\frac{t}{2} S(x^{\sss (n)},t)$.
From Theorem \ref{theo:aldous-mc-s}(ii), $\{S(x^{\sss (n)},t)\}_{n\ge 1}$ is tight.  Combining these facts we have that
$\{\sum_{i=1}^{\infty} \hat Y_i(x^{\sss (n)},t)\}_{n\ge 1}$ is a tight family.  Finally, note that
$\hat R(x^{\sss (n)},t) \le X_1(x^{\sss (n)},t)\left (\sum_{i=1}^{\infty} \hat Y_i(x^{\sss (n)},t) \right)$.  The tightness of $\{ \hat R(x^{\sss (n)},t) \}_{n \ge 1}$
now follows on combining the above established tightness of $\{\sum_{i=1}^{\infty} \hat Y_i(x^{\sss (n)},t)\}_{n\ge 1}$ and the tightness
of $\{X_1(x^{\sss (n)},t)\}_{n\ge 1}$, where the latter is once again a consequence of Theorem \ref{theo:aldous-mc-s}(ii).
\qed\\

We now complete the proof of Theorem \ref{theo:welldef-feller}.\\

{\bf Proof of Theorem \ref{theo:welldef-feller}.} Fix $z=(x,y) \in \udown^0$ and $t \ge 0$. From Lemma \ref{lemma:786} (i)   $\hat R(x,t) < \infty$ a.s. 
Also, from Proposition \ref{lemma:211},  $\tilde R(z,t) < \infty$ a.s.  The result now follows on recalling that
$R(z,t) = \hat R(x,t) + \tilde R(z,t)$. \qed\\

We also record the following consequence of Lemma \ref{lemma:786} and Corollary \ref{corr:ab1157} for future use.
\begin{Corollary} \label{lemma:tight-r}
	If $z^{\sss (n)} \to z$ in $\udown^0$, then
	$ \{R(z^{\sss (n)},t)\}_{n\ge 1}  $ is tight.
\end{Corollary}

\subsection{Feller property of the augmented MC}
\label{sec:feller-property}

In this section, we will prove Theorem \ref{theo:welldef-fellerb}. In fact we will show that 
if $ z^{\sss (n)}= (x^{\sss (n)},y^{\sss (n)})$ converges to $z=(x,y)$ in $\udown^0$, and $z \in \udown^1$, then
\begin{equation}\label{eq:ab1200}
	 (\bfX(z^{\sss(n)},t) ,\bfY(z^{\sss(n)},t)) \convp (\bfX(z,t), \bfY(z,t)). \end{equation}
The following lemma is immediate from the definition of $\bfd_{\sss \mathbb{U}}(\cdot, \cdot)$.
\begin{Lemma}
	\label{lemma:basic-analysis}
	Suppose $(x,y), (x^{\sss (n)}, y^{\sss (n)}) \in \udown$ for $n \ge 1$. Then $\lim_{n \to \infty} \bfd_{\sss \mathbb{U}}((x,y),(x^{\sss (n)},y^{\sss (n)}))=0$ if and only if the following three conditions hold:\\
	(i) $\lim_{n \to \infty}\sum_{i=1}^\infty (x_i^{\sss (n)}-x_i)^2 = 0$.\\
	(ii) $y_i^{\sss (n)}=y_i$ for $n$ sufficiently large, for all $ i \ge 1$.\\
	(iii) $\lim_{n \to \infty}\sum_{i=1}^\infty x_i^{\sss (n)}y_i^{\sss (n)} = \sum_{i=1}^\infty x_iy_i$.
\end{Lemma}
The key ingredient in the proof is the following lemma the proof of which is given after Lemma \ref{strict}.
\begin{Lemma} \label{lemma:feller-conv-prob}
 Let $ z^{\sss (n)}= (x^{\sss (n)},y^{\sss (n)})$ converge to $z=(x,y)$ in $\udown^0$.  Suppose that $z \in \udown^1$. Then\\
	(i) $Y_i(z^{\sss (n)},t) \convp Y_i(z,t)$ for all $i \ge 1$.\\
	(ii) $\sum_{i=1}^\infty X_i(z^{\sss (n)},t)Y_i(z^{\sss (n)},t) \convp \sum_{i=1}^\infty X_i(z,t)Y_i(z,t)$. 
\end{Lemma}
Proof of Theorem \ref{theo:welldef-fellerb} can now be completed as follows.\\

{\bf Proof of Theorem \ref{theo:welldef-fellerb}.}
The first part of the theorem is immediate from the construction given at the beginning of Section \ref{sec:main-constr-mc} and elementary properties of Poisson processes.
For the second part, consider $ z^{\sss (n)}= (x^{\sss (n)},y^{\sss (n)})$, $z=(x,y)$ as in the statement of the theorem.  It suffices to prove \eqref{eq:ab1200}.
From Theorem \ref{theo:aldous-mc-s}(ii), $\bfX(z^{\sss(n)},t) \to \bfX(z,t)$ in probability, in $\ldown$. The result now follows on combining this convergence with
the convergence in Lemma \ref{lemma:feller-conv-prob} and applying Lemma \ref{lemma:basic-analysis}.  \qed \\

Rest of this section is devoted to the proof of Lemma \ref{lemma:feller-conv-prob}.
The key idea of the proof is as follows.  Consider the induced subgraphs on the first $k$ vertices $\bfG^{\sss [k]}=\bfG(z^{\sss [k]},t)$ and $\bfG^{\sss (n)[k]}=\bfG(z^{\sss (n)[k]},t)$. Since there are only finite number of vertices in $\bfG^{\sss [k]}$,  when $n \to \infty$, $\bfG^{\sss (n)[k]}$ will eventually be identical to $\bfG^{[k]}$ almost surely. The main step in the proof is   to control the difference between $\bfG^{\sss (n)[k]}$ and $\bfG^{\sss (n)}$  when $k$ is large, uniformly for all $n$. For this we first analyze the difference between $\bfG^{\sss (n)[k]}$ and $\bfG^{\sss (n)[k+1]}$ in the lemma below.

	Consider the set of vertices $[k+1]=\{1,2,...,k, k+1\}$, and for every  $i \in [k+1]$, let vertex $i$ have label $(x_i, y_i)$ representing its size  and surplus, respectively. Suppose $x_1 \ge x_2 \ge ... \ge x_{k+1}$.  Fix $t > 0$. Define a random graph $\bfG^*$  on the above vertex set as follows. For $i\le k$, the number of edges, $N_i$, between $i$ and $k+1$ follows Poisson$(tx_i x_{k+1})$. In addition, there are $N_0$ = Poisson$(t x^2_{k+1}/2)$ self-loops to the vertex $k+1$. All the Poisson random variables  are taken to be mutually
	independent.\\
	 Denote $X_i$ and $Y_i$ for the component volumes and surplus of the resulting star-like graph if $i$ is the smallest labeled vertex in its component; otherwise let $X_i=Y_i=0$. A precise definition of $(X_i, Y_i)$ is as follows. Write $i \sim k+1$ if there is an edge between $i$ and $k+1$ in $\bfG^*$. By convention $(k+1) \sim (k+1)$. Let
	$\clj_k = \{i \in [k+1]: i \sim k+1\}$, and $i_0 = \min\{i: i \in \clj_k\}$.
	Then
	\[
	(X_i, Y_i) =\left\{
	\begin{array}{cc}
	\left( \sum_{i\in \clj_k}x_i, \sum_{i\in \clj_k}y_i\right)  &   \mbox{ if } i = i_0   \\ \ \\
	(0,0)   &  \mbox{ if } i \in \clj_k\setminus \{i_0\} \\ \ \\
	(x_i,y_i) &  \mbox{ if } i \in [k+1]\setminus\clj_k.
	\end{array}
	\right .
	\]
Define $R_k=\sum_{i=1}^k x_i y_i$, $S_k=\sum_{i=1}^k x_i^2$, $R_{k+1}=\sum_{i=1}^{k+1} X_iY_i$.  Then we have the following result.

\begin{Lemma}
	\label{lemma:basic-dynamic}
	(i) $\prob\{ Y_i \neq y_i \} \le tx_{k+1}y_{k+1}  x_1 + t x_{k+1}^2 \left( 1 + it x_1^2 + tS_k +tR_k x_1 \right)$.\\
	(ii)  $ \E[R_{k+1}-R_k] \le  x_{k+1}y_{k+1}  ( 1 + t S_k) + x_{k+1}^2  (t R_k + t^2 S_k R_k + t^2 S_k x_1) + tx_{k+1}^3  
	(1 + 2tS_k  +t^2S_k^2). $
\end{Lemma}
{\bf Proof:} (i) It is easy to see that, for $i = 1, \cdots k$,
$$
\{Y_i \neq y_i\} \subset \left( \{y_{k+1} > 0\} \cap \{i \in \clj_k\}\right)
\cup \{N_0 \neq 0\} \cup_{j=1}^k \{N_j > 1\} \cup_{j < i} \{N_jN_i \neq 0\} \cup_{j: y_j > 0} \{N_jN_i \neq 0\}.$$
% The event $\{Y_i \neq y_i\}$ can be covered by the union of the following events: (1) $y_{i+1} > 0$ and $k+1$ and $i$ are linked; (2) the Poisson$(t x_{k+1}^2)$ self-loops at vertex $k+1$ is not zero; (3) for some $ j \le k$, $j$ and $i$ are linked by more than one edge; (4) for some $j < i$, $j$ and $i$ are connected through $k+1$; (5) for some $j \le k$ with $y_j >0$, $j$ and $i$ are connected through $k+1$. 
Using the observation that for a Poisson$(\lambda)$ random variable $Z$, $\prob\{Z \ge 1\} < \lambda$ and $\prob\{ Z \ge 2\} < \lambda^2$, we now have
that

\begin{align*}
	\prob\{ Y_i \neq y_i \} 
	\le& tx_i x_{k+1} \cdot y_{k+1} + \frac {t x_{k+1}^2}{2} + \sum_{j=1}^k (tx_jx_{k+1})^2 \\
	+& \sum_{j=1}^{i-1} tx_jx_{k+1} \cdot tx_i x_{k+1} + \sum_{j=1}^k tx_jx_{k+1} \cdot tx_i x_{k+1} \cdot y_j.
\end{align*}
 Proof is now completed on collecting all the terms and using the fact that $x_i \ge x_1$  for every $i$.\\
(ii) Note that  $$X_0 =x_{k+1} + \sum_{j=1}^k x_j {\ind}_{\{ N_j \ge 1 \}}, \; Y_0 =y_{k+1} + \sum_{j=1}^k y_{j} {\ind}_{\{ N_j \ge 1 \}} 
+ N_{0} + \sum_{j=1}^k (N_j - 1)^+ .$$ Then
\begin{align*}
	R_{k+1}-R_k =& X_0Y_0 - \sum_{j \in \clj_k} x_jy_j \\
	=& x_{k+1}y_{k+1} + \sum_{j=1}^k (x_jy_{k+1}+x_{k+1}y_j){\ind}_{\{ N_j \ge 1 \}} + \sum_{1\le j < l \le k} (x_j y_l +x_l y_j) {\ind}_{\{ N_j \ge 1 \}}{\ind}_{\{ N_l \ge 1 \}}\\
	+& N_{0} X_0 + x_{k+1} \sum_{j=1}^k (N_j-1)^+ + \sum_{j=1}^k x_j (N_j-1)^+ \\
	+& \sum_{1 \le j < l \le k}(x_j {\ind}_{\{ N_j \ge 1 \}} (N_l-1)^+ + x_l {\ind}_{\{ N_l \ge 1 \}} (N_j-1)^+).
\end{align*}
The result now follows on taking expectations in the above equation  and using the fact that $\E[(N_j-1)^+] < (tx_jx_{k+1})^2 $. \qed \\

Recall that, by construction, $X_i(z,t) \ge X_{i+1}(z,t)$ for all $z \in \udown$, $t\ge 0$ and $i \in \NNN$.  The following lemma which is a key ingredient in the proof of Lemma \ref{lemma:feller-conv-prob} says
that if $z \in \udown^1$, ties do not occur, a.s.

% \todo[inline]{The proof below still needs work.}
\begin{Lemma}
	\label{strict}
	Let $z \in \udown^1$.  Then for every $t > 0$ and $i \in \NNN$,
	$X_i(z,t) > X_{i+1}(z,t)$ a.s.
\end{Lemma}
\noindent \textbf{Proof:} Fix $t > 0$.  Consider the graph $\bfG(z,t)$ and write $\CC_{x_i} \equiv \CC_{x_i}(t)$ for the component of vertex $(x_i,y_i)$ at time $t$. It suffices to show for all $i\neq j$
\begin{equation}
	\label{eqn:1257}
	\prob\set{ |\CC_{x_i}| = |\CC_{x_j}|, \CC_{x_i} \ne  \CC_{x_j} } = 0.
\end{equation}
The key property we shall use is that for $z=(x,y) \in \udown^1$, $\sum_{i=1}^\infty x_i =\infty$. Now fix $i\geq 1$. It is enough to show that $|\CC_{x_i}|$ has no atom i.e for all $(x,y)\in \udown^1$ 
\begin{equation}
\label{eqn:no-atom-one}
	\prob(|\CC_{x_i}| =a) = 0, \qquad \mbox{ for any } a\geq 0. 
\end{equation} 
To see this, first note that since $|\CC_{x_i}|< \infty$ a.s.,  conditional on $\CC_{x_i}$ the vector $z^* = ((x_k,y_k): x_k\notin \CC_{x_i}) \in \udown^1$ almost surely. Thus on the event $x_j\notin \CC_{x_i}$, conditional on $\CC_{x_i}$, using \eqref{eqn:no-atom-one} with $a= |\CC_{x_i}|$ implies that $\prob(|\CC_{x_j}| = |\CC_{x_i}|\mid \CC_{x_i}) = 0$ and this completes the proof. Thus it is enough to prove \eqref{eqn:no-atom-one}. For the rest of the argument, to ease notation let $i=1$.  Let us first show the simpler assertion that the volume of direct neighbors of $x_1$ has a continuous distribution. More precisely, let  $N_{i,j}(t):=\xi_{i,j}([0,tx_ix_j/2])+\xi_{j,i}([0,tx_ix_j/2])$, $1 \le i < j$, denote the number of edges between any two vertices  $x_i$ and $x_j$ by time $t$. Then the volume of \emph{direct} neighbors of the vertex $x_1$ is $L := \sum_{i =  2 }^\infty x_i {\ind}_{\set{N_{1,i}(t) \ge 1}} $ and we will first show that  $L$ has no atom, namely
\begin{equation}
	\label{eqn:1261}
\prob(L = a) = 0, \qquad \mbox{ for all } a\geq 0. 	
\end{equation}
 For any random variable $X$ define the maximum atom size of $X$ by 
$${\bf atom}(X) : = \sup _{ a \in \RRR} \prob\set{X = a}.$$
For two independent random variables $X_1$ and $X_2$ we have ${\bf atom}(X_1 + X_2) \le \min \set{ {\bf atom}(X_1), {\bf atom}(X_2) }$. For $m \ge 2$, define $L_m = \sum_{i=m}^\infty x_i {\ind}_{\set{N_{1,i}(t) \ge 1}}$. Since $L_m$ and $L-L_m$ are independent, we have
$ {\bf atom}(L) \le {\bf atom}(L_m).$
Define the event
$$ E_m:=\set{ N_{1,i}(t) \le 1 \mbox{ for all } i \ge m },  $$  
and write 
$$ L^*_m(t) := \sum_{i=m}^\infty x_i N_{1,i}(t).$$
Then $L^*_m(t)$ is a pure jump Levy process with Levy measure $\nu(du) = \sum_{i=m}^\infty x_1 x_i \delta_{x_i}(du)$. By \cite{hartman1942infinitesimal}, such a Levy process has continuous marginal distribution since the Levy measure is infinite ($ \nu(0,\infty) = (\sum_{i=m}^\infty x_i) x_1 = \infty$) . Thus $L^*_m(t)$ has no atom. Next, for any $a \in \RRR$,
\begin{align*}
	\prob \set{ L_m = a} 
	\le& \prob \set{ E_m^c} + \prob \set{ E_m, L_m =a }
	= \prob \set{ E_m^c} + \prob \set{ E_m, L^*_m(t) =a }\\
	\le& \sum_{i=m}^\infty \frac{(tx_1x_i)^2}{2} + 0
	= \frac{t^2 x_1^2}{2} \sum_{i=m}^\infty x_i^2.
\end{align*}
Thus ${\bf atom}(L) \le {\bf atom}(L_m) \le \frac{t^2 x_1^2}{2} \sum_{i=m}^\infty x_i^2$. Since $m$ is arbitrary, we have ${\bf atom}(L)=0$. Thus $L$ is a continuous variable, and \eqref{eqn:1261} is proved. 

Let us now strengthen this to prove \eqref{eqn:no-atom-one}. % Then denote $\tilde \CC_{x_1}$ for the volume of component in $\bfG(z,t)$ that contains $x_1$.
 Let $\tilde \bfG$ be the subgraph of $\bfG(z,t)$ obtained by deleting the vertex $x_1$ and all related edges. Let $\tilde X_i$ be the volume of the $i$-th largest component of $\tilde \bfG$. Note that $\sum_{i=1}^\infty \tilde X_i = \sum_{i=2}^\infty x_i = \infty$ a.s. Conditional on $(\tilde X_i)_{i\geq 1}$, let $\tilde N_{1,i}$ have Poisson distribution with parameter $t x_1 \tilde X_i$. Then
\[ \CC_{x_1} \stackrel{d}{=} x_1 + \sum_{i=1}^\infty {\tilde X}_i\ind_{\set{\tilde N_{1,i}\geq 1}},\]
where the second term has the same form as the random variable $L$. Using \eqref{eqn:1261} completes the proof. \qed
 %  by the above argument we have 
 % $$ \prob \{ \tilde \CC_{x_1} = a | \tilde X_i : i \ge 1 \} = 0 \mbox{ for all } a \in \RRR, $$
 % and this implies $\tilde \CC_{x_1}$ has no atom.
 % 
 % Now we are ready to show \eqref{eqn:1257}. Given $\tilde \CC_{x_1}$ and the fact that $\tilde \CC_{x_2} \ne \tilde \CC_{x_1}$, the rest of the graph has the same law as $\bfG(\tilde z,t)$, where $\tilde z$ is gained by deleting those vertices that are in $\tilde \CC_{x_1}$. Thus we have a.s. 
 % $$ \prob \{  \tilde \CC_{x_2} = a |  \tilde \CC_{x_1} \} {\ind}_{ \{ x_2 \notin \tilde \CC_{x_1} \}} =0 \mbox{ for all } a\in \RRR. $$
 % This implies a.s.
 % $$ \prob \{  \tilde \CC_{x_2} =  \tilde \CC_{x_1}|  \tilde \CC_{x_1} \} {\ind}_{ \{ x_2 \notin \tilde \CC_{x_1} \}} =0. $$
 % Taking expectation on both sides of the above expression will prove \eqref{eqn:1257}. \qed\\

% \todo[inline]{Now it is super rigorous..not quite!..now it is!}

We now proceed to the proof of Lemma \ref{lemma:feller-conv-prob}.

{\bf Proof of Lemma \ref{lemma:feller-conv-prob}.} 
Fix $t > 0$ and $z^{\sss (n)}, z$ as in the statement of the lemma.  Denote
 $ Y^{\sss [k]} =Y(z^{\sss [k]},t), \;\; Y^{\sss (n)[k]}= Y(z^{\sss (n)[k]},t)$.
Similarly, denote $\CC^{\sss [k]}_i$ and $\CC_i^{\sss (n)[k]}$ for the corresponding $i$-th largest component; and  $X_i^{\sss [k]}$ and $X_i^{\sss (n)[k]}$ for their respective sizes. Also, write $X^{\sss (n)} = X(x^{\sss (n)},t)$ and define $Y^{\sss (n)}, R^{\sss (n)}, S^{\sss (n)}$ similarly.

For  $i \in \NNN$, define the event $E_i^{\sss (n)[k]}$ as,
$$ E_i^{\sss (n)[k]} := \{\omega: X_j^{\sss (n)[k]}(\omega) > X_{j+1}^{\sss (n)}(\omega), \mbox{ for } j=1,2,...,i \}, $$
and define $E_i^{\sss [k]}$ similarly. 
Then
\begin{align}
	\prob\{ Y_i^{\sss (n)} \neq Y_i(t) \}  
	\le& \prob\{ Y_i^{\sss (n)} \neq Y_i^{\sss (n)[k]} \} + \prob\{Y_i^{\sss (n)[k]} \neq Y_i^{\sss [k]}\} +\prob\{ Y_i^{\sss [k]} \neq Y_i(t)\} \nonumber\\
	\le& \prob\{ Y_i^{\sss (n)} \neq Y_i^{\sss (n)[k]}, E_i^{\sss (n)[k]} \} + \prob\{ (E_i^{\sss (n)[k]})^c \} + \prob\{Y_i^{\sss (n)[k]} \neq Y_i^{\sss [k]}\} +\prob\{ Y_i^{\sss [k]} \neq Y_i(t)\}. \label{eqn:314}
\end{align}
Note that
$$
E_i^{\sss (n)[k]} \subset \{\omega: \CC_j^{\sss (n)[k]}(\omega) \subset \CC_j^{\sss (n)[m]}(\omega) \subset \CC_j^{\sss (n)}(\omega), \mbox{ for all } j=1,2,...,i \mbox{ and } m \ge k\}.$$
Thus the probability of the event $\{ Y_i^{\sss (n)[m+1]} \neq Y_i^{\sss (n)[m]}, E_i^{\sss (n)[k]} \}$, for $m \ge k$, can be estimated using Lemma \ref{lemma:basic-dynamic} (i). More precisely, let  $\FF^{\sss [m]}=\sigma \{ \xi_{i,j}; {i , j \le m} \}$ for $ m \ge 1$. Then by Lemma \ref{lemma:basic-dynamic} (i),
\begin{align*}
	\prob\{ Y_i^{\sss (n)[m+1]} \neq Y_i^{\sss (n)[m]}, E_i^{\sss (n)[k]} | \FF^{\sss [m]} \} \le& t x^{\sss (n)}_{m+1}y^{\sss (n)}_{m+1} X^{\sss (n)[m]}_1 \\
	+&
	t (x^{\sss (n)}_{m+1})^2 \left( 1 + it (X^{\sss (n)[m]}_1)^2 + tS^{\sss (n)[m]} +t R^{\sss (n)[m]} X^{\sss (n)[m]}_1 \right), %\label{eqn:324}
\end{align*}
where $S^{\sss (n)[m]} = \sum_i (X_i^{\sss (n)[m]})^2$ and $R^{\sss (n)[m]} = \sum_i (X_i^{\sss (n)[m]}Y_i^{\sss (n)[m]})$.

Note that $X_1^{\sss (n)[k]} \le X_1^{\sss (n)}$, $R^{\sss (n)[k]} \le R^{\sss (n)}$ and $S^{\sss (n)[k]} \le S^{\sss (n)}$. Thus we have
\begin{align*}
	\sum_{m=k}^\infty \prob\{ Y_i^{\sss (n)[m+1]} \neq Y_i^{\sss (n)[m]}, E_i^{\sss (n)[k]} | \FF^{\sss [m]} \} 
	 \le& t 
	\left(\sum_{m=k+1}^\infty x^{\sss (n)}_{m}y^{\sss (n)}_{m}\right) X_1^{\sss (n)} \\
	+& t \left(\sum_{m=k+1}^\infty (x^{\sss (n)}_{m})^2\right) \left( 1 + it (X_1^{\sss (n)})^2 + 
	tS^{\sss (n)} +tR^{\sss (n)}X_1^{\sss (n)} \right).
\end{align*}
Denote the right hand side of the above inequality as  $U^{\sss (n)[k]}$. Then by Lemma \ref{lemma:borel-cantelli}(ii), we have
\be
 \prob \{ Y_i^{\sss (n)} \neq Y_i^{\sss (n)[k]}, E_i^{\sss (n)[k]} \} = \prob \left(
\cup_{m=k}^\infty \{ Y_i^{\sss (n)[m+1]} \neq Y_i^{\sss (n)[m]}, E_i^{\sss (n)[k]} \}\right) \le 2 \E[U^{\sss (n)[k]} \wedge 1] \label{eqn:879}
\ee
and therefore
\begin{align}
	\prob\{ Y_i^{\sss (n)} \neq Y_i(t) \}  
	\le 2\E[ U^{\sss (n)[k]} \wedge 1 ] + \prob\{ (E_i^{\sss (n)[k]})^c \} + \prob\{Y_i^{\sss (n)[k]} \neq Y_i^{\sss [k]}\} +\prob\{ Y_i^{\sss [k]} \neq Y_i(t)\}.
	\label{eqn:879.5}
\end{align}
% We prove the tightness property of $U^{\sss (n)[k]}$ in the next lemma:
% \begin{Lemma}
% 	\label{lemma:tight-uniformly-u}
% 	
% \end{Lemma}

Next note that $X^{\sss (n)}_1$, $S^{\sss (n)}$ and $R^{\sss (n)}$ are all tight sequences by Corollary \ref{lemma:tight-r} and Theorem \ref{theo:aldous-mc-s}(ii).
Thus $(1+ it  (X_1^{\sss (n)})^2 + tS^{\sss (n)} + tR^{\sss (n)} X_1^{\sss (n)})$ is also tight. Also, since $z^{\sss (n)} \to z$,
$$ \limsup_{k \to \infty} \limsup_{ n \to \infty} \sum_{i=k+1}^\infty x^{\sss (n)}_i y^{\sss (n)}_i =0 \mbox{ and }  \limsup_{k \to \infty} \limsup_{ n \to \infty}  \sum_{i=k+1}^\infty(x^{\sss (n)}_i )^2 =0.$$
Combining the above observations we have that
 $\limsup_{k \to \infty} \limsup_{ n \to \infty} \prob\{U^{\sss (n)[k]}>\epsilon \}=0$ for all $\epsilon > 0$. From the inequality
$$\E [ U^{\sss (n)[k]} \wedge 1]  \le  \prob\{ U^{\sss (n)[k]} > \epsilon \} + \epsilon$$ we now see that
\begin{equation} \label{lemma:tight-uniformly-u} \limsup_{k \to \infty} \limsup_{n \to \infty} \E [ U^{\sss (n)[k]} \wedge 1]=0. \end{equation}
	Next, from a straightforward extension of Proposition 5 of Aldous \cite{aldous1997brownian} we have that
	 $(\bfX^{\sss (n)}, X_1^{\sss (n)[k]}, ..., X_i^{\sss (n)[k]} ) \convd (\bfX(t), X_1^{\sss [k]}, ..., X_i^{\sss [k]} )$ in $\ldown \times \mathbb{R}^i$ when $n \to \infty$, for each fixed $i$ and $k$.  Combining this with Lemma \ref{strict} we now see that for fixed $i$
$$ \limsup_{k \to \infty}\limsup_{n \to \infty} \prob\{ (E_i^{\sss (n)[k]})^c \} = 0.$$
Also, for each fixed $k$
%Since $\lim_{n \to \infty} Y_i^{\sss (n)[k]}=Y_i^{\sss [k]}$ for a finite graph, then we have that for fixed $k$
$$
\limsup_{n \to \infty} \prob\{ Y_i^{\sss (n)[k]} \neq Y_i^{[k]}\} = 0.$$
Observing that $\lim_{k \to \infty} Y_i^{\sss [k]}=Y_i(t)$ and the last term in \eqref{eqn:879.5} does not depend on $n$,  we have that
$$
\limsup_{k \to \infty} \limsup_{n \to \infty} \prob\{ Y_i^{\sss [k]} \neq Y_i(t)\} = 0.$$
Part (i) of the lemma now follows on combining the above observations and taking limit as $n \to \infty$ and then $k \to \infty$ in \eqref{eqn:879.5}. 
% 
% 
% 
% 
% 
% Taking the limit $\limsup_{k \to \infty} \limsup_{ n \to \infty}$ on both side and then letting $\epsilon \to 0$ will complete the proof. \qed\\
% 
% {\bf Completing the proof of Lemma \ref{lemma:feller-conv-prob}(i): } By \eqref{eqn:314} \eqref{eqn:879}, we have
% 
% Take $\limsup_{n \to \infty}$ on both side,
% 
% 
% 
% 
% 
% 
% \begin{align*}
% \limsup_{k \to \infty} 	\limsup_{n \to \infty} \prob\{ Y_i^{\sss (n)}(t) \neq Y_i(t) \} 0 + \prob\{ (E_i^{\sss [k]})^c \} + 0 +\prob\{ Y_i^{\sss [k]} \neq Y_i\}.
% \end{align*}
%  The third term goes to zero since the related graph is finite. Then we let $k \to \infty$ for the right hand side. The first term goes to zero by Lemma \ref{lemma:tight-uniformly-u} andthe rest terms all go to zero since $ \lim_{k \to \infty} X_i^{\sss [k]}=X_i$ and $\lim_{k \to \infty} Y_i^{\sss [k]}=Y_i$ a.s. The proof is thusly completed. \qed\\

% 
% 
% \subsubsection{ Proof of Lemma \ref{lemma:feller-conv-prob}(ii) }
% 
% {\bf Proof of  Lemma \ref{lemma:feller-conv-prob}(ii):} 

We now prove part (ii) of the lemma.
% 
% 
% Actually we are going to prove the following lemma:
% \begin{Lemma}
% 	\label{lemma:tilde-s-converge}
% 	Under the $\xi$-coupling, recall $R^{\sss (n)}=\sum_{i=1}^\infty X^{\sss (n)}_i  Y^{\sss (n)}_i$ and $R=\sum_{i=1}^\infty X_i Y_i$, then\\
% 	(i) $ \liminf_{n \to \infty} R^{\sss (n)} \ge R$ a.s.\\
% 	(ii) For any $\epsilon > 0$,  $ \lim_{n \to \infty} \prob\{ R^{\sss (n)} > R + \epsilon \} =0$.
% \end{Lemma}
% {\bf Proof:} (i) This is immediate once one notices that
Note that
$$ \liminf_{n \to \infty} R^{\sss (n)} \ge \lim_{n \to \infty} R^{\sss (n)[k]} = R^{[k]}. $$
% \todo[inline]{Not quite true: Also, since $Y_i^{\sss [k]}$ increases to $Y_i(t)$ and $X_i^{\sss [k]}$ increases to $X_i(t)$, as $k\to \infty$, we have by monotone convergence that $R^{[k]} \to R(z,t)$ as $k \to \infty$. }
With a similar argument as in Lemma \ref{lemma:248}, we have $R^{[k]} \to R(z,t)$ as $k \to \infty$. Thus sending $k \to \infty$ in the above display we have
\begin{equation}\label{eq:ab1746}\liminf_{n \to \infty} R^{\sss (n)} \ge R(z,t).\end{equation}
To complete the proof, it suffices to show that
\begin{equation}
	\mbox{ For any } \epsilon > 0,\,   \lim_{n \to \infty} \prob\{ R^{\sss (n)} > R(z,t) + \epsilon \} =0.\label{eq:ab1928}
\end{equation}
Note that
\begin{align}
	\prob\{ R^{\sss (n)} - R(z,t) > \epsilon \} 
	\le& \prob\{ R^{\sss (n)} - R^{\sss (n)[k]} > \epsilon/2 \} + \prob\{ R^{\sss (n)[k]} - R(z,t) > \epsilon/2 \}\nonumber\\
	\le& \prob\{ R^{\sss (n)} - R^{\sss (n)[k]} > \epsilon/2 \} + \prob\{ R^{\sss (n)[k]} - R^{\sss [k]} > \epsilon/2 \}.\label{eq:ab1956}
\end{align} 
The second term on the right side above goes to zero for each fixed $k$, as $n \to \infty$. For the first term, note that
% \be
% \limsup_{k \to \infty} \limsup_{n \to \infty} \prob\{ R^{\sss (n)} - R^{\sss (n)[k]} > \epsilon \} =0. \label{eqn:471}
% \ee
% 
% recall the filtration $\FF^{\sss [k]}=\sigma \{ \xi_{i,j}; i,j \le k \}$, then 
by Lemma \ref{lemma:basic-dynamic}(ii), for all $m \ge k$
\begin{align*}
	\E[ R^{\sss (n)[m+1]}- R^{\sss (n)[m]} | \FF^{\sss [m]} ] 
	\le  x^{\sss (n)}_m y^{\sss (n)}_m  U_1^{\sss (n)} +  (x^{\sss (n)}_m)^2 U_2^{\sss (n)} + (x^{\sss (n)}_{m+1})^3 U_3^{\sss (n)},
\end{align*}
where $U_1^{\sss (n)}= 1 + t S^{\sss (n)}$, $ U_2^{\sss (n)}= t R^{\sss (n)} + t^2 S^{\sss (n)}R^{\sss (n)}+t^2 S^{\sss (n)} X_1^{\sss (n)}$ and $U_3^{\sss (n)}= t(1 + 2t S^{\sss (n)} + t^2 (S^{\sss (n)})^2)$. Thus by Lemma \ref{lemma:borel-cantelli} (i),
$$ \prob\{ R^{\sss (n)} - R^{\sss (n)[k]} > \epsilon \} \le (1 + 1/\epsilon) \E[ U^{\sss (n)[k]} \wedge 1], $$
where $U^{\sss (n)[k]}=( \sum_{m=k+1}^\infty x^{\sss (n)}_m y^{\sss (n)}_m ) U_1^{\sss (n)} +  (\sum_{m=k+1}^\infty(x^{\sss (n)}_m)^2) U_2^{\sss (n)} + (\sum_{m=k+1}^\infty(x^{\sss (n)}_{m+1})^3) U_3^{\sss (n)} $. Note that $U_1^{\sss (n)}$, $U_2^{\sss (n)}$ and $U_3^{\sss (n)}$ are all tight sequences
and $z^{\sss (n)} \to z$.  An argument similar to the one used to prove \eqref{lemma:tight-uniformly-u} now shows that, for all $\epsilon > 0$,
$$ \limsup_{k \to \infty} \limsup_{n \to \infty} \prob\{ R^{\sss (n)} - R^{\sss (n)[k]} > \epsilon \}  \le \left(1 + \frac{1}{\epsilon}\right) \limsup_{k \to \infty} \limsup_{n \to \infty}  \E[ U^{\sss (n)[k]} \wedge 1] =0.$$
The statement in \eqref{eq:ab1928} now follows on using the above convergence in \eqref{eq:ab1956} and combining it with the observation below \eqref{eq:ab1956}.
This completes the proof of part (ii). \qed\\
\begin{Remark}
	\label{prodnofell}
	Lemma \ref{lemma:feller-conv-prob} is at the heart of the (near) Feller property in Theorem \ref{theo:welldef-fellerb}
	which is crucial for the proof of the joint convergence in \eqref{eq:eq1241}. 
	The proof of the lemma reveals the reason for considering the metric $\bfd_{\sss \mathbb{U}} $ on $\udown$ rather than $\bfd_1$ or $\bfd_2$
%	\footnote{Replaced: $\bfd_{vt}$} (see Section \ref{sec:augmented-mc}).  
	The proof hinges upon the convergence of
	$\sum_{m=1}^\infty x^{\sss (n)}_m y^{\sss (n)}_m $ to $\sum_{m=1}^\infty x_m y_m $, as $n\to \infty$, even for the proof of convergence of
	$Y_i(z^{\sss (n)},t) \convp Y_i(z,t)$.  This suggests that the convergence in $\bfd_1$ or
	$\bfd_2$ is ``too weak'' to yield the desired Feller property.
\end{Remark}

% \subsubsection{Right continuity of the process}
% \label{sec:right-cont}
% \todo[inline]{Xuan can you see if you can fill in the details?}

\section{The standard augmented multiplicative coalescent.}
\label{sec:main-st-ag-mc-ex}
In this section we prove Theorem \ref{thm:smc-surplus}.  Proposition 4 of \cite{aldous1997brownian} proves a very useful result on convergence of component size 
vectors of a general family of non-uniform random graph models to the ordered excursion lengths of $\hat W _{\lambda}$. We begin in this section by extending this
result to the joint convergence of component size and component surplus vectors in $\udown$, under a slight strengthening of the conditions assumed in
\cite{aldous1997brownian}.  

% We want to construct a special version of the augmented multiplicative coalescent $(\bfX(\lambda), \bfY(\lambda))_{\lambda \in \RRR}$, so called the {\bf standard augmented multiplicative coalescent}, which is the limit object in the critical window as stated in Theorem \ref{thm:crit-regime}. The first coordinate process $\bfX(\cdot)$ is just the standard multiplicative coalescent as introduced in \cite{aldous1997brownian}. 
% \todo[inline]{check where those objects are defined.}
% 
% \subsubsection{The marginal distribution for a fixed $\lambda$}
Recall the excursion lengths and  mark count process $\bfZ^*(\lambda) = (\bfX^*(\lambda), \bfY^*(\lambda))$ defined in Section \ref{sec:augmented-mc}. 
%\footnote{Add: In this section we will use $q$ in $Z(z,q)$ to denote the time parameter}In this section we will use $q$ in $Z(z,q)$ to denote the time parameter. 
Our first
result below shows that, for fixed $\lambda \in \RRR$, $\bfZ^*(\lambda)$ arises as a limit of $\bfZ(z^{\sss (n)}, q^{\sss (n)})$ in $\udown$ for all sequences $\{z^{\sss (n)}\} \subset \udown$ and
$q^{\sss (n)}=q^{\sss (n)}_\lambda \subset (0, \infty)$ that satisfy certain regularity conditions.

% We will identify this as the marginal distribution $(\bfX(\lambda), \bfY(\lambda))$ for the standard augmented MC. The first step is to show $(\bfX^*(\lambda), \bfY^*(\lambda))$ is an universal limit for under very general setting. More specifically, we will show that 
% $$(\bfX(x^{\sss (n)}, y^{\sss (n)}, t^{\sss (n)}),\bfY(x^{\sss (n)}, y^{\sss (n)}, t^{\sss (n)})) \convd (\bfX^*(\lambda), \bfY^*(\lambda)),$$
% given some general regularity conditions on $x^{\sss (n)}, y^{\sss (n)}$ and $t^{\sss (n)}$. Note that we have studied the limit when $(x^{\sss (n)},y^{\sss (n)}) \to (x,y)$ and $t$ is fixed in section \ref{sec:feller-property}, but here $t$ also changes with $n$. We uses similar ideas as in \cite{aldous1997brownian} but requires more refined analysis due to the special metric on $\udown$. To keep our notation consistent with \cite{aldous1997brownian} and save $t$ for the time parameter in the random walk we are going to define, we use $q^{\sss (n)}$ instead of $t^{\sss (n)}$ in this section, and we often drop $n$ in the notation for simplicity. \\

For $n \ge 1$, let $z^{\sss (n)} = (x^{\sss (n)}, y^{\sss (n)}) \in \udown^0$.
 %(z_i^{\sss (n)}: i \ge 1 )\in \udown^0$ be such that $z^{\sss (n)}_i = (0,0)$ for all $i >n$.
Writing $z_i^{\sss (n)} = (x_i^{\sss (n)}, y_i^{\sss (n)})$, $i \ge 1$, define 
$$x^{*\sss (n)}= \sup_{i \ge 1} x^{\sss (n)}_i, \; s^{\sss (n)}_r = \sum_{i=1}^\infty (x_i^{\sss (n) })^r, \; r \ge 1.$$
 Let $\{q^{\sss (n)}\}$ be a nonnegative sequence. We will suppress $(n)$ from the notation unless needed.
%The main theorem of this section is
\begin{Theorem} 
	\label{theo:aldous-full-gene}
	Let $z^{\sss (n)} =  (z_1^{\sss (n)}, \cdots )\in \udown^0$ be such that $z^{\sss (n)}_i = (0,0)$ for all $i >n$.
	Suppose that, as $n\to \infty$,
	\begin{equation}
		\label{eqn:qsigma2}		
	\frac{s_3}{(s_2)^3} \to 1, \quad
	q- \frac{1}{s_2} \to \lambda , \quad
	\frac{x^*}{s_2} \to 0, 	
	\end{equation}
	and, for some $\varsigma \in (0, \infty)$,
	\begin{equation}
		\label{eqn:additional-condition}
	s_1 \cdot \left(\frac{x^*}{s_2}\right)^\varsigma \to 0.
	\end{equation}
	
	Further suppose that $y_i^{\sss(n)} = 0$ for all $i$. Then $\bfZ^{\sss (n)} = \bfZ(z^{\sss (n)},q^{\sss (n)})$ converges in distribution in $\udown$ to $\bfZ^*(\lambda)$.
	% $$ (\bfX^{\sss (n)}, \bfY^{\sss (n)}) := (\bfX(x^{\sss (n)}, y^{\sss (n)}, q^{\sss (n)}),\bfY(x^{\sss (n)}, y^{\sss (n)}, q ^{\sss (n)})) \convd (\bfX^*(\lambda), \bfY^*(\lambda)).$$
\end{Theorem}
{\bf Remark:} The convergence assumption in \eqref{eqn:qsigma2}	is  the same as that in Proposition 4 of \cite{aldous1997brownian}.  The additional assumption in
\eqref{eqn:additional-condition}  is not very stringent as will be seen in Section \ref{sec:main-coupling} when this
result is applied to a general family of bounded-size rules.

Given Theorem \ref{theo:aldous-full-gene}, the proof of Theorem \ref{thm:smc-surplus} can now be completed as follows.

{\bf Proof of Theorem \ref{thm:smc-surplus}.}
The first two parts of the theorem were shown in Theorem \ref{theo:welldef-fellerb}.  
Also, part (v) of the theorem is immediate from the definition of $\{T_t\}$ in Section \ref{sec:mult}. 
Recall the definition of $\nu_{\lambda}$ from
Section \ref{sec:augmented-mc}.  In order to prove parts (iii)-(iv) it suffices to  show that 
\begin{equation}
	\label{eq:ab2121}
	\mbox{ for any } \lambda_1, \lambda_2 \in \RRR, \, \lambda_1 \le \lambda_2,\; \nu_{\lambda_1}\clt_{\lambda_{2}-\lambda_1} = \nu_{\lambda_2}.
	\end{equation}
Indeed, using the semigroup property of $(\clt_{\lambda})$ and the above relation, it is  straightforward to define a consistent
family of finite dimensional distributions $\mu_{\lambda_1, \cdots \lambda_k}$ on $(\udown)^{\otimes k}$, $-\infty < \lambda_1 < \lambda_2, \cdots
\lambda_k<\infty$, $k \ge 1$, such that $\mu_{\lambda} = \nu_{\lambda}$ for every $\lambda \in \RRR$.  The desired result then follows from Kolmogorov's consistency 
theorem.  

We now prove \eqref{eq:ab2121}.  
% Consider the \erdos random graph model, where  we have $n$ isolated vertices at time $0$, and for each pair of vertices, edges are added at the rate $1/n$. In addition, for each vertex, self-loops are added at the rate $1/2n$. 
% Denote by $\CC_i^{\sss (n)}(t)$  the size of the $i$-th largest component at time $t$ and define
% $$ \bar\bfX^{\sss (n)}(\lambda)= \left( \frac{1}{n^{2/3}} \CC_i^{\sss (n)}\left(1 + \frac{\lambda}{n^{1/3}}\right): i \ge 1\right). $$
% Also, let $\bar Y_i^{\sss (n)}(\lambda)$ be the surplus of the $i$-th largest component at time $1 + \lambda/n^{1/3}$ and let 
% $\bar\bfY^{\sss (n)}(\lambda) = (\bar Y_i^{\sss (n)}(\lambda))_{i\ge 1}$, $\bar\bfZ^{\sss (n)}(\lambda) = (\bar\bfX^{\sss (n)}(\lambda), \bar\bfY^{\sss (n)}(\lambda))$.
% It is easy to verify that for any $\lambda_1 \le \lambda_2$, the joint distribution of $(\bar \bfZ^{\sss (n)}(\lambda_1), \bar \bfZ^{\sss (n)}(\lambda_1))$
% is same as that of 
Let
$$
z^{\sss (n)} = (x^{\sss (n)}, y^{\sss (n)}), \; x_i^{\sss (n)} = n^{-2/3}, \; y_i^{\sss (n)} = 0, \; i = 1, \cdots n, \; q^{\sss (n)}_{\lambda_j} = \lambda_j + n^{1/3}, \; j = 1,2.$$
We set $z^{\sss (n)}_i =0$ for $i > n$.
% $(\bfZ(z^{\sss (n)}, q^{\sss (n)}_{\lambda_1}), \bfZ(z^{\sss (n)}, q^{\sss (n)}_{\lambda_2}))$, where
% $$
% z^{\sss (n)} = (x^{\sss (n)}, y^{\sss (n)}), \; x_i^{\sss (n)} = n^{-2/3}, \; y_i^{\sss (n)} = 0, \; i = 1, \cdots n, \; q^{\sss (n)}_{\lambda_j} = \lambda_j + n^{1/3}, \; j = 1,2.$$
% 
% 
% In Theorem \ref{theo:aldous-full-gene}, this corresponds to starting with $n$ vertices with weight $x_i = n^{-2/3}$ and $q=n^{1/3}+\lambda$. To verify assumptions \eqref{eqn:qsigma2} and \eqref{eqn:additional-condition} 
%\todo[inline]{Instead of using $\nu_\lambda$ everywhere, shall we just write $\mathcal{L}(Z^*(\lambda))$ to denote the law of $Z^*(\lambda)$}
Note that with this choice of $x^{\sss (n)}$, $s_1=n^{1/3}, s_2 = n^{-1/3}, s_3 = n^{-1}$ and so clearly \eqref{eqn:qsigma2} and \eqref{eqn:additional-condition} (with
any $\varsigma > 1$) are satisfied with $q = q_{\lambda_j}$, $\lambda = \lambda_j$, $j=1,2$.  Thus, denoting
the distribution of $\bfZ(z^{\sss (n)}, q^{\sss (n)}_{\lambda_j})$ by $\nu_{\lambda_j}^{\sss (n)}$, we have by Theorem \ref{theo:aldous-full-gene} that
\begin{equation}
	\label{eq:eqab2153}
	\nu_{\lambda_j}^{\sss (n)} \to \nu_{\lambda_j}, \mbox{ as } n \to \infty. 
\end{equation}
 Also, from the construction of $\bfZ(z,t)$ in Section \ref{sec:main-constr-mc}, it is clear that
$\nu_{\lambda_1}^{\sss (n)}\clt_{\lambda_{2}-\lambda_1} = \nu_{\lambda_2}^{\sss (n)}$.  The result now follows on combining the convergence in
\eqref{eq:eqab2153} with Theorem \ref{theo:welldef-fellerb} and observing that $\bfZ^*(\lambda) \in \udown^1$ a.s. for every $\lambda \in \RRR$. 
%\todo[inline]{add the proof for the fifth item in the theorem.} 
\qed\\

% so that the assumptions are verified.  Thus
% $$ (X^{\sss (n)}(\lambda), Y^{\sss (n)}(\lambda)) \convd (\bfX^*(\lambda), \bfY^*(\lambda)) ~ \mbox{in}~ \udown. $$
% Let $T_t$ be the generator of the coalescing paired process as defined in the very beginning, since the \erdos process construct above is meant to be this way, then for any fixed $n$, $\lambda_1 < \lambda_2$ we have
% $$T_{\lambda_2-\lambda_1} (\bfX^{\sss (n)}(\lambda_1), \bfY^{\sss (n)}(\lambda_1)) =_d (\bfX^{\sss (n)}(\lambda_2), \bfY^{\sss (n)}(\lambda_2)).$$
% Thus by the Feller property, immediately we have
% $$T_{\lambda_2-\lambda_1} (\bfX^*(\lambda_1), \bfY^*(\lambda_1)) =_d (\bfX^*(\lambda_2), \bfY^*(\lambda_2)).$$
% Then by Kolmogorov consistency theorem, we can construct the standard coalescing paired process on $\lambda \in (-\infty, +\infty)$.
%\todo[inline]{Xuan: typically surplus edges in components means new edges between vertices. In your previous definition when construction the non-uniform random graph this was fine since we thought of each of the vertex as a component (with more than one vertex) and self loops meant a new edge being created within a component.  In this construction one can get self loop at a single vertex. Does this change things?}

% We include \eqref{eqn:additional-condition} in the assumptions in order to prove the weak convergence in $\udown$. Note that this condition is very weak and it is easy to be satisfied.\\

Rest of this section is devoted to the proof of Theorem \ref{theo:aldous-full-gene} and is organized as follows.  Recall the random graph process
$\bfG(z,q)$, for $z \in \udown$, $q\ge 0$, defined at the beginning of Section \ref{sec:main-constr-mc}.  In Section \ref{sec:conv-in-prod} we will
give an equivalent in law construction of $\bfG(z,q)$, from \cite{aldous1997brownian}, that defines the random graph
simultaneously with a certain breadth-first-exploration random walk.  The excursions of the reflected version of this walk encode
the component sizes of the random graph while the area under the excursions gives the parameter of the Poisson distribution describing
the (conditional) law of the surplus associated with the corresponding component. Using this construction, in Theorem \ref{theo:generalized-aldous},
we will first prove a weaker result than Theorem \ref{theo:aldous-full-gene} which proves the convergence in distribution of $\bfZ^{\sss (n)}$ to $\bfZ^*(\lambda)$
in $\ldown \times \NNN^{\infty}$, where we consider the product topology on $\NNN^{\infty}$.  This result is proved in Section \ref{sec:sec6.2}.
In Section \ref{sec:conv-d-metric} we will give the proof of Theorem \ref{theo:aldous-full-gene} using Theorem \ref{theo:generalized-aldous} and an auxiliary
tightness lemma (Lemma \ref{lemma:tightness-sup}).  Finally, proof of Lemma \ref{lemma:tightness-sup} is given in Section \ref{sec:sec6.3new}.

% {\bf Proof:} We will prove this theorem in two steps. By Proposition 4 of \cite{aldous1997brownian}, we already have $\bfX^{\sss (n)} \convd \bfX^*(\lambda)$ in $\ldown$. So the first step is generalizing this to convergence of $(\bfX^{\sss (n)}, \bfY^{\sss (n)}) \convd (\bfX^*,\bfY^*)$ in $\ldown \times \NNN^\infty$, namely, the joint convergence of $(\bfX^{\sss (n)}, Y^{\sss (n)}_i) \convd (\bfX^{\sss (n)}, Y^*_i)$ for any fixed $i$. Then second step is establishing the convergence in $\udown$.\\
% The first step will be completed in Section \ref{sec:conv-in-prod} and the second step will be done in Section \ref{sec:conv-d-metric}. \qed 

\subsection{Breadth First Exploration Walk. }
\label{sec:conv-in-prod}

In this section, following \cite{aldous1997brownian}, we will give an equivalent in law construction of $\bfG(z,q)$ that defines the random graph
simultaneously with a certain breadth-first-exploration random walk.  Given $q \in (0, \infty)$ and $z \in \udown^0$ such that $x_i=0$ for all $i > n$ and $y_i=0$ for all $i$,
we will construct a random graph $\bar \bfG(z,q)$ that is equivalent in law to  $\bfG(z,q)$, in two stages, as follows.
We begin with a graph on $[n]$ with no  edges.
Let $\{\eta_{i,j}\}_{i,j \in \NNN}$ be independent Poisson point processes on $[0,\infty)$ such that $\eta_{ij}$ for $i\neq j$ has intensity
$qx_j$; and for $i=j$ has intensity $qx_i/2$.\\
{\bf Stage I: The breadth-first-search forest and associated random walk:} 
Choose a vertex $v(1) \in [n]$ with $\prob(v(1) = i)\propto x_i$.  Let
$$\III_1 = \{ j \in [n]: j \neq v(1) \mbox{ and } \eta_{v(1),j} \cap [0, x_{v(1)}] \neq \emptyset \}.$$
Form an edge between $v(1)$ and each $j \in \III_1$.  Let $c(1) = \# (\III_1)$.
Let $m_{v(1),j}$ be the first point in $\eta_{v(1),j}$ for each $j \in \III_1$.  Order the vertices in $\III_1$ according to increasing 
values of $m_{v(1),j}$ and label these as $v(2), \cdots v(c(1) + 1)$.  Let 
$$\VV_1 =\{v(1)\}, \; \NN_1 = \{v(2), \cdots , v(c(1) +1)\}, \;  l_1 = x_{v(1)} \mbox{ and } d_1 = c(1).$$
Having defined $\VV_{i'}$, $\NN_{i'}$, $l_{i'}$, $d_{i'}$  and the edges up to step $i'$, with 
$\VV_{i'} = \{v(1), \cdots v(i')\}$, $\NN_{i'} = \{ v(i'+1), v(i'+2), \cdots v(d_{i'} +1)\}$  for $1 \le i' \le i-1$, define,
if $\NN_{i-1} \neq \emptyset$
$$
\III_i = \{ j \in [n]: j \not \in \NN_{i-1}\cup\VV_{i-1} \mbox{ and } \eta_{v(i),j} \cap [0, x_{v(i)}] \neq \emptyset \}$$
and form an edge between $v(i)$ and each $j \in \III_i$.  Let $c(i) = |\III_i|$ and
let $m_{v(i),j}$ be the first point in $\eta_{v(i),j}$ for each $j \in \III_i$.  Order the vertices in $\III_i$ according to increasing 
values of $m_{v(i),j}$ and label these as $v(d_{i-1}+2), \cdots v(d_i + 1)$, where $d_i = d_{i-1} + c(i)$.  Set
$$
 l_i = l_{i-1} + x_{v(i)}, \;
\VV_{i} = \{v(1), \cdots v(i)\}, \; \NN_{i} = \{ v(i+1), v(i+2), \cdots v(d_{i} +1)\}.$$
In case $\NN_{i-1} = \emptyset$, we choose $v(i)\in [n]\setminus \VV_{i-1}$ with probability proportional to $x_j$, $j \in [n]\setminus \VV_{i-1}$
and define $\III_i, c(i), d_i, l_i, \VV_i, \NN_i$ and the edges at step $i$ exactly as above.

This procedure terminates after exactly $n$ steps at which point we obtain a forest-like graph with no surplus edges.  We will include
surplus to this graph in stage II below.

Associate with the above construction an (interpolated) random walk process $Z^{\sss (n)}(\cdot)$ defined as follows.  $Z^{\sss (n)}(0) = 0$ and  
\begin{equation}
	\label{eqn:def_zn}
	Z^{\sss (n)}(l_{i-1}+u)=Z^{\sss (n)}(l_{i-1})-u+\sum_{j \notin \VV_i \cup \NN_{i-1}} x_{j}{\ind}_{\{m_{v(i),j}<u\}} \quad \mbox{for}~ 0< u< x_{v(i)},\;
	i = 1, \cdots n,
\end{equation}
 where by convention $l_0 =0$ and $\NN_0 = \emptyset$.  This defines $Z^{\sss (n)}(t)$ for all $t \in [0, l_n)$.  Define
$Z^{\sss (n)}(t)=Z^{\sss (n)}(l_n-)$ for all $t \ge l_n$.\\

% Finally  update $\NN_i:= \NN_{i-1} \cup \{v \in \VV : v \notin  \VV_i \cup \NN_{i-1} , \eta^*_{(v(i),v)}< x_i \} \setminus \{v(i)\}$, and let $\VV_{i}=\VV_{i-1} \cup \{v(i)\}$. One technique issue is that we can only define $Z^{\sss (n)}(t) $ up to time $l_n=\sum_{i=1}^n x_i$ by the above construction, so for $t > l_n$, define $Z^{\sss (n)}(t)=Z^{\sss (n)}(l_n)$. \\ 

{\bf Stage II: Construction of surplus edges:}
%************ begin old version ************\\
%We construct surplus edges on the graph obtained in Stage I and a point process $\clp_x$ on $[0, l_n]$, simultaneously, as follows.
%For $i = 1, \cdots n$ and each $\tau $ in $\eta_{v(i),j}\cap [0, x_{v(i)}]$, $\tau \neq m_{v(i),j}$, $j \in \III_i$
%[resp. in $\eta_{v(i), v(i)} \cap [0, x_{v(i)}]$; in $\eta_{v(i), v(j)} \cap [0, x_{v(i)}]$, $v(j) \in \NN_{i-1} \setminus \{v(i)\}$ ] 
%construct an edge between $v(i)$ and $j$ [resp. between $v(i)$ and itself; between $v(i)$ and $v(j)$] and construct points for the point process
%$\clp_x$ at time instant $l_{i-1} + \tau \in [0, l_n]$.\\
%************** end old version *************** \footnote{I suggest the following new version of construction of the surplus, which have more explanation. }\\
For each $i=1, \cdots, n$, we construct surplus edges on the graph obtained in Stage I and a point process $\clp_x$ on $[0, l_n]$, simultaneously, as follows. \\
(i) For each $v \in \III_i$ and $\tau  \in \eta_{v(i),v}\cap [0, x_{v(i)}] \setminus \{ m_{v(i),v} \}$, construct an edge between $v(i)$ and $v$. This corresponds to multi-edges between the two vertices $v(i)$ and $v$.\\
(ii) For each $\tau \in \eta_{v(i), v(i)} \cap [0, x_{v(i)}]$, construct an edge between $v(i)$ and itself. This corresponds to self-loops at the vertex $v(i)$.\\
(iii) For each $v(j) \in \NN_{i-1} \setminus \{v(i)\}$ and $\tau \in \eta_{v(i), v(j)} \cap [0, x_{v(i)}]$, construct an edge between $v(i)$ and $v(j)$. This corresponds to additional edges between two vertices, $v(i)$ and $v(j)$, that were indirectly connected in stage I.\\
For each of the above cases, we also construct points for the point process $\PP_x$ at time $l_{i-1}+\tau \in [0,l_n]$.\\

% Note that the above construction only gives the \emph{new} vertices not yet found by the exploration process that are connected to $v(i)$. To complete constructing the entire graph, we have to add the following three sources of surplus:\\
% (i) the self-loops at $v(i)$ for each point in $\eta_{v(i), v(i)} \cap [0, x_{v(i)}]$;\\
% (ii) the edges between $v(i)$ and $v(j) \in \NN_{i-1} \setminus \{v(i)\}$ for each point in $\eta_{v(i), v(j)} \cap [0, x_{v(i)}]$;\\
% (iii) the additional edges between $v(i)$ and $v(j) \in \NN_i \setminus \NN_{i-1}$ for each point in $\eta_{v(i), v(j)} \cap [0, x_{v(i)}] \setminus \{ \eta^*_{v(i),v(j)} \}$.\\
% Those surplus edges are also added sequentially according to the time $t$. For example, in the situation (ii), if $\eta^* \in \eta_{v(i), v(j)} \cap [0, x_{v(i)}]$, then add one edge between $v(i)$ and $v(j)$ at the time $l_{i-1}+\eta^*$. \\

This completes  the construction of the graph $\bar \bfG(z,q)$ and the random walk $Z^{\sss (n)}(\cdot)$. This graph has the same law as $\bfG(z,q)$, so the associated component sizes and surplus vector denoted as $(\bar \bfX(z,q), \bar \bfY(z,q))$ has the same law as that of $( \bfX(z,q),  \bfY(z,q))$.
Furthermore, conditioned on $Z^{\sss (n)}$, $\clp_x$ is Poisson point process on $[0, l_n]$ whose intensity we denote by $r_x(t)$.

Using the above construction we will show in next section, as a first step,  a  weaker result than Theorem \ref{theo:aldous-full-gene}.

\subsection{Convergence in $\ldown \times \NNN^{\infty}$.}
\label{sec:sec6.2}
The following is the main result of this section.
\begin{Theorem}
	\label{theo:generalized-aldous}
	Let $z^{\sss (n)} \in \udown^0$ and $q^{\sss (n)} \in (0,\infty)$ be sequences that satisfy the conditions in Theorem \ref{theo:aldous-full-gene}.  Then
	\begin{equation}
	\label{eqn:fd-con-surpl}
		(\bfX^{\sss (n)},\bfY^{\sss (n)}) \convd (\bfX^*(\lambda), \bfY^*(\lambda)) 
	\end{equation}
	in the space $\ldown \times \NNN^\infty$ as $n\to \infty$, where we consider the product topology on $\NNN^\infty$.
\end{Theorem}	
The key ingredient in the proof is the following result.  With $z^{\sss (n)}$ and $q^{\sss (n)}$ as in the above theorem, define
$\bar \bfX^{\sss (n)} = \bar \bfX(z^{\sss (n)}, q^{\sss (n)})$, $\bar \bfY^{\sss (n)} = \bar \bfY(z^{\sss (n)}, q^{\sss (n)})$ and $r^{\sss (n)}(t) = r_{x^{\sss (n)}}(t)1_{[0, l_n]}(t)$, $t \ge 0$.
Denote   the random walk process from Section \ref{sec:conv-in-prod} constructed using $(x^{\sss (n)}, q^{\sss (n)})$ (rather than $(x,q)$), once more, by $Z^{\sss (n)}(\cdot)$.

Define the rescaled process $\bar{Z}^{\sss (n)}(\cdot)$ and its reflected version $\hat Z^{\sss (n)} (\cdot)$ as follows
 \begin{equation}
 \label{eqn:reflec-def}
 	\bar{Z}^{\sss (n)}(t) := \sqrt{\frac{s_2}{s_3}} Z^{\sss (n)}(t), \quad \hat Z^{\sss (n)} (t):= \bar{Z}^{\sss (n)}(t) - \min_{0\leq u\leq t} \bar{Z}^{\sss (n)}(u).
 \end{equation}
\begin{Lemma}\label{lem:lem1325}
(i) As $n \to \infty$, the process $\bar{Z}^{\sss (n)} \convd W_\lambda$ in $\DD ([0,\infty): \RRR)$.

(ii) For $n \ge 1$,
\begin{equation}
	\label{eqn:rt-q-small}
	\sup_{t \ge 0} \left|r^{\sss (n)}(t)-q \sqrt{\frac{s_3}{s_2}} \hat Z^{\sss (n)}(t) \right|
	\le \frac{3}{2} q x_*.
\end{equation}
\end{Lemma}
Given Lemma \ref{lem:lem1325}, the proof of Theorem \ref{theo:generalized-aldous} can be completed as follows.\\

{\bf Proof of Theorem \ref{theo:generalized-aldous}:}
%\todo[inline]{Move proof to appendix?}
The paper \cite{aldous1997brownian} shows that the vector $\bar \bfX^{\sss (n)}$ can be represented as the ordered sequence of excursion lengths of the 
process $\hat{Z}^{\sss (n)} $.  Also, weak convergence of  $\bar{Z}^{\sss (n)}$ to $W_\lambda$ in Lemma \ref{lem:lem1325} (i) implies the convergence
of $\hat{Z}^{\sss (n)} $ to $\hat W_{\lambda}$.  Using these facts, Proposition 4 of \cite{aldous1997brownian} shows that
$\bar \bfX^{\sss (n)}$ converges in distribution to the ordered excursion length sequence of $\hat W_{\lambda}$, namely $\bfX^*(\lambda)$, in $\ldown$.
Also, conditional on $\hat{Z}^{\sss (n)}$, $\clp_x$ is a Poisson point process on $[0, \infty)$ with rate $r^{\sss (n)}(t)$ and for $i \ge 1$,
$\bar Y_i^{\sss (n)}$ has a Poisson distribution with parameter
$\int_{[a_i^{\sss (n)}, b_i^{\sss (n)}]} r^{\sss (n)}(s) ds$, where $a_i^{\sss (n)}, b_i^{\sss (n)}$ are the left and right endpoints of the $i$-th ordered excursion
of $\hat Z^n$.  From conditions in \eqref{eqn:qsigma2} it follows that $ q  x^* \to 0$ and  $ q \sqrt{s_3/s_2} \to 1$.   Lemma \ref{lem:lem1325} (ii) 
then shows that $\int_{[a_i^{\sss (n)}, b_i^{\sss (n)}]} r^{\sss (n)}(s) ds$ converges in distribution to $\int_{[a_i, b_i]} \hat W_{\lambda}(s) ds$,
where $a_i, b_i$ are the left and right endpoints of the $i$-th ordered excursion
of $\hat W_{\lambda}$.
In fact we have the joint convergence of $\left(\hat{Z}^{\sss (n)}, \left(\int_{[a_i^{\sss (n)}, b_i^{\sss (n)}]} r^{\sss (n)}(s) ds\right)_{i\ge 1}\right)$
to $\left(\hat W_{\lambda}, \left(\int_{[a_i, b_i]} \hat W_{\lambda}(s) ds\right)_{i\ge 1}\right)$.  This proves 
the convergence of 	$(\bar\bfX^{\sss (n)},\bar\bfY^{\sss (n)})$ to $(\bfX^*(\lambda), \bfY^*(\lambda))$ in  $\ldown \times \NNN^\infty$.
The result follows since  $(\bar \bfX^{\sss (n)}, \bar \bfY^{\sss (n)})$ has the same law as  $( \bfX^{\sss (n)},  \bfY^{\sss (n)})$. \qed\\

{\bf Proof of Lemma \ref{lem:lem1325}} Part (i) was proved in Proposition 4 of \cite{aldous1997brownian}. Consider now (ii).\\
% The conditions in \eqref{eqn:qsigma2} imply that $ q  x^* \to 0$ and $ q \sqrt{\sigma_3/\sigma_2} \to 1$. Thus using results of part (i) and \eqref{eqn:rt-q-small} and classic weak convergence theory (see e.g.\cite{kallenberg}) now gives the joint convergence of the components and surplus. We will now prove \eqref{eqn:rt-q-small}.\\
% We start by making two observations about the associated walk:
% \\(a) The random walk above is a jump process with negative unit drift. 
% \\(b) For each $i \ge 1$, there is a jump of size $x_{v(i)}$, unless $v(i)$ is the first vertex in its component. Define $\delta_{v(j)}=1$ if $v(j)$ is the first vertex in its component and $\delta_{v(j)}=0$ otherwise.\\
% With the above observations, no matter $\NN_{i-1} = \emptyset$ or $\NN_{i-1}=\{v(i), v(i+1), ... , v(i+l)\}$, $Z^{\sss (n)}(l_{i-1})$ can be expressed as:
% $$ Z^{\sss (n)}(l_{i-1})= - \sum_{j=1}^{i-1} \delta_{v(j)} x_{v(j)} + \sum_{v \in \NN_{i-1}} x_v.$$
% Thus $\inf_{j \le i-1} Z^{\sss (n)}(l_{j-1}) = - \sum_{j=1}^{i-1} \delta_{v(j)} x_{v(j)}$. In addition, for $t \in (l_{i-1},l_i]$, by the recursive definition in \eqref{eqn:def_zn}, we have
It is easy to verify that $Z^{\sss (n)}$ satisfies
$$ Z^{\sss (n)}(l_{i})= - \sum_{j=1}^{i} \delta_{v(j)} x_{v(j)} + \sum_{v \in \NN_{i}} x_v,\; i = 1, \cdots n.$$
The above equation implies that for all $k \le i$,  $Z^{\sss (n)}(l_{k}) \ge - \sum_{j=1}^{i} \delta_{v(j)} x_{v(j)}$. In addition, taking $k_0 = \sup \set{j \le i: \delta_{v(j)}=1} $ we have $Z^{\sss (n)}(l_{k_0}) = - \sum_{j=1}^{i} \delta_{v(j)} x_{v(j)}$. In particular, this implies that $\inf_{j \le i} Z^{\sss (n)}(l_{j}) = - \sum_{j=1}^{i} \delta_{v(j)} x_{v(j)}$.
   Also, from \eqref{eqn:def_zn} we have that for $t \in (l_{i-1},l_i]$,
$ Z^{\sss (n)}(t) \ge Z^{\sss (n)}(l_{i-1}) - x^*$. Consequently
\begin{equation}
	\label{eqn:rw_minimum}
	\left|\inf_{0 \le u \le t}Z^{\sss (n)}(u) +\sum_{j=1}^{i-1} \delta_{v(j)} x_{v(j)}\right| = \left|\inf_{0 \le u \le t}Z^{\sss (n)}(u) - \inf_{\set{j : l_{j} \le t}} Z^{\sss (n)}(l_{j})\right|   \le x^*.
\end{equation}
Let $\NN_{i-1}=\{v(i), v(i+1), ... , v(i+l)\}$.  From the above expression for $Z^{\sss (n)}(l_{i})$, we have that for $t \in (l_{i-1},l_i]$
\begin{equation}
	\label{eqn:rw_zn}
	Z^{\sss (n)}(t) =  \left(- \sum_{j=1}^{i-1} \delta_{v(j)} x_{v(j)} + \sum_{j=i}^{i+l} x_{v(j) } \right) - \left(t -l_{i-1} \right) + \sum_{j \notin \VV_i \cup \NN_{i-1}} x_j {\ind}_{\{m_{v(i),j}<t-l_{i-1}\}},
\end{equation}
Also, accounting for the three sources of surplus described in Stage II of the construction, one has
 the following formula for $r^{\sss (n)}(t)$ at time $t \in (l_{i-1},l_i]$:
$$ r^{\sss (n)}(t) = q \cdot \left( \frac{x_{v(i)}}{2}+ \sum_{j=i+1}^{i+l}x_{v(j)} + \sum_{j \notin \VV_i \cup \NN_{i-1}} x_j {\ind}_{\{m_{v(i),j}<t-l_{i-1}\}} \right). $$

The three terms in the above expression correspond to self-loops;  edges between vertices that in stage I were only connected indirectly; and  additional edges between two vertices that were directly connected in stage I. Combining the above expression with \eqref{eqn:rw_zn} and \eqref{eqn:rw_minimum}, we have
\begin{align}
	\label{eqn:rate-surplus}
	\left|r^{\sss (n)}(t)-q \cdot \left(Z^{\sss (n)}(t)-\min_{0 \le s \le t}Z^{\sss (n)}(s)\right) \right|
	\le q \cdot \left( \left| \inf_{0 \le s \le t}Z^{\sss (n)}(s) +\sum_{j=1}^{i-1} \delta_{v(j)} x_{v(j)} \right| + \frac {x_{v(i)}}{2} \right)
	\le \frac{3}{2} q x_*.
\end{align}
The result follows. \qed

%\todo[inline]{continue from here}

\subsection{Proof of Theorem \ref{theo:aldous-full-gene}. }
\label{sec:conv-d-metric}
In this section we complete the proof of Theorem \ref{theo:aldous-full-gene}.  The key step in the proof is the following lemma whose
proof is given in Section \ref{sec:sec6.3new}.
\begin{Lemma}
	\label{lemma:tightness-sup}
	Let $z^{\sss (n)} \in \udown^0$ and $q^{\sss (n)} \in (0,\infty)$ be as in Theorem \ref{theo:aldous-full-gene}.  Let $\hat Z^{\sss (n)}$ be as 
	introduced in \eqref{eqn:reflec-def}.  Then
	$ \{\sup_{t \ge 0} \hat Z^{\sss (n)}(t) \}_{n \ge 1}$ is a tight family of $\RRR_+$ valued random variables.
\end{Lemma}
\begin{Remark}
	In fact one can establish a stronger statement, namely $ \sup_{u \ge t} \sup_{n \ge 1} \hat Z_u^{\sss (n)} \to 0$ in probability
	as $t \to \infty$.  Also, although not used in this work, using very similar techniques as in the proof of Lemma \ref{lemma:tightness-sup}, it can be shown
	that $ \sup_{u \ge t} \hat W_\lambda(u) $ converges a.s. to $0$, as $t \to \infty$.
\end{Remark}

{\bf Proof of Theorem \ref{theo:aldous-full-gene}.}  Since $(\bfX^{\sss (n)}, \bfY^{\sss (n)})$ has the same distributions as
$(\bar\bfX^{\sss (n)}, \bar\bfY^{\sss (n)})$, we can equivalently consider the convergence of the latter sequence.
From Theorem \ref{theo:generalized-aldous} we have that $(\bar\bfX^{\sss (n)}, \bar\bfY^{\sss (n)})$ converges to
$(\bfX^*(\lambda), \bfY^*(\lambda))$, in distribution, in $\ldown \times \NNN^{\infty}$ (with product topology on $\NNN^{\infty}$).
By appealing to Skorohod representation theorem, we can assume without loss of generality that the convergence is almost sure.  
In view of Lemma \ref{lemma:basic-analysis}, it now suffices to argue that 
	$$ \sum_{i=1}^\infty \left|\bar X^{\sss (n)}_i \bar Y^{\sss (n)}_i- X^*_i(\lambda) Y^*_i(\lambda)\right| \convp 0. $$
	Fix $\epsilon > 0$.  Then, for any $k \in \NNN$,
\begin{align}
\prob\set{ \sum_{i=1}^{\infty} \left|\bar X^{\sss (n)}_i \bar Y^{\sss (n)}_i- X^*_i(\lambda) Y^*_i(\lambda)\right| > \epsilon }\le &
\prob\set{ \sum_{i=1}^k \left|\bar X^{\sss (n)}_i \bar Y^{\sss (n)}_i- X^*_i(\lambda) Y^*_i(\lambda)\right| > \frac{\epsilon}{3} } \nonumber \\
+& \prob\set{ \sum_{i=k+1}^\infty \bar X^{\sss (n)}_i \bar Y^{\sss (n)}_i > \frac{\epsilon}{3}} +\prob\set{ \sum_{i=k+1}^\infty X^{*}_i(\lambda) Y^{*}_i(\lambda) > \frac{\epsilon}{3}}. \label{eqn:1108}	
\end{align}
From the convergence of $(\bar\bfX^{\sss (n)}, \bar\bfY^{\sss (n)})$  to
$(\bfX^*(\lambda), \bfY^*(\lambda))$ in $\ldown \times \NNN^{\infty}$ we have that
$$
\lim_{n\to \infty} \prob\set{ \sum_{i=1}^k \left|\bar X^{\sss (n)}_i \bar Y^{\sss (n)}_i- X^*_i(\lambda) Y^*_i(\lambda)\right| > \frac{\epsilon}{3} }  = 0.$$
Consider now the second term in \eqref{eqn:1108}.  Let 	$E^{\sss (n)}_L = \{ \sup_{t\ge 0} r_t^{\sss (n)} \le L \}$.  Then
$$
\prob\set{ \sum_{i=k+1}^\infty \bar X^{\sss (n)}_i \bar Y^{\sss (n)}_i > \frac{\epsilon}{3}}
\le \prob\set{( E^{\sss (n)}_L)^c} +  \frac{3}{\epsilon} \E \left ({\ind}_{E^{\sss (n)}_L}\left[
\sum_{i=k+1}^\infty \bar X^{\sss (n)}_i \bar Y^{\sss (n)}_i \wedge 1\right] \right).
$$
Let $\GG = \sigma\{\hat Z^{\sss (n)}(t): t \ge 0\}$.  Since $r_t^{\sss (n)}$ is $\GG$ measurable  for all $t \ge 0$,  $E^{\sss (n)}_L \in \GG$. 
%\footnote{Add: Since $r_t^{\sss (n)}$ is measurable in $\GG$ for all $t \ge 0$, thus $E^{\sss (n)}_L \in \GG$.} 
 Then
\begin{align*}
	\E \left ({\ind}_{E^{\sss (n)}_L}\left[
	\sum_{i=k+1}^\infty \bar X^{\sss (n)}_i \bar Y^{\sss (n)}_i \wedge 1 \right]\right)
	=& \E \left ({\ind}_{E^{\sss (n)}_L}\E\left[
	\sum_{i=k+1}^\infty \bar X^{\sss (n)}_i \bar Y^{\sss (n)}_i \wedge 1 \mid \GG \right] \right)\\
	\le& \E \left ({\ind}_{E^{\sss (n)}_L}
	\left(\sum_{i=k+1}^\infty \E\left[\bar X^{\sss (n)}_i \bar Y^{\sss (n)}_i \mid \GG \right]\wedge 1 \right) \right)\\
	\le& L \E \left[\sum_{i=k+1}^\infty (\bar X^{\sss (n)}_i)^2 \wedge 1  \right],
\end{align*}
where the last inequality follows on observing that, conditionally on $\GG$, $\bar Y^{\sss (n)}_i$ has a Poisson distribution with
rate that is dominated by $\bar X^{\sss (n)}_i \cdot (\sup_{t\ge 0} r_t^{\sss (n)})$.  Using the convergence of $\bar \bfX^{\sss (n)}$ to $\bfX^*$, we now have
$$
\limsup_{n\to \infty} 
\E \left ({\ind}_{E^{\sss (n)}_L}\left[
\sum_{i=k+1}^\infty \bar X^{\sss (n)}_i \bar Y^{\sss (n)}_i\right] \wedge 1 \right) \le  L \E \left[\sum_{i=k+1}^\infty ( X^*_i(\lambda))^2 \wedge 1  \right].$$
Let $\delta > 0$ be arbitrary.  Using Lemma \ref{lemma:tightness-sup} and Lemma \ref{lem:lem1325} (ii) we can choose $L \in (0, \infty)$ such that
$\prob\set{( E^{\sss (n)}_L)^c} \le \delta$.
Finally, taking limit as $n \to \infty$ in \eqref{eqn:1108} we have that 
\begin{align}
	\limsup_{n\to \infty}\prob\set{ \sum_{i=1}^{\infty} \left|\bar X^{\sss (n)}_i \bar Y^{\sss (n)}_i- X^*_i(\lambda) Y^*_i(\lambda)\right| > \epsilon }
	\le & \delta + L \E \left[\sum_{i=k+1}^\infty ( X^*_i(\lambda))^2 \wedge 1  \right]\nonumber\\
	 +&  \prob\set{ \sum_{i=k+1}^\infty X^{*}_i(\lambda) Y^{*}_i(\lambda) > \frac{\epsilon}{3}}.
\end{align}
The result now follows on sending $k \to \infty$ in the above display and recalling that
$\sum_{i=1}^\infty ( X^*_i(\lambda))^2 < \infty$ and $\sum_{i=1}^\infty X^{*}_i(\lambda) Y^{*}_i(\lambda) < \infty$ a.s. and $\delta > 0$ is arbitrary. \qed

\subsection{Proof of Lemma \ref{lemma:tightness-sup}.}
\label{sec:sec6.3new}
In this section we prove Lemma \ref{lemma:tightness-sup}.  We will only treat the case  $\lambda = 0$.  The general case can be treated similarly. The key step in the proof is the following proposition whose proof is given at the end of the section.
 
 Note that
$\sup_{t \ge 0}|\bar Z^{\sss (n)}(t) - \bar Z^{\sss (n)}(t-)| \le x^* \sqrt{s_2/s_3} \to 0$ as $n \to \infty$.  
Also, as $n \to \infty$, $qs_2 \to 1$.
Thus, without loss of generality, we will assume that
\begin{equation}\label{eq:eq1547}
	\sup_{n \ge 1}\sup_{t \ge 0}|\bar Z^{\sss (n)}(t) - \bar Z^{\sss (n)}(t-)| \le 1, \; \sup_{n \ge 1} q^{\sss (n)}s_2^{\sss (n)} \le 2 .
\end{equation}
Fix $\vartheta \in (0, 1/2)$ and define $t^{* \sss (n)} = \left(\frac{s_2}{x^*}\right)^{\vartheta}$. 
Denote by $\{\FF_t^{\sss (n)}\}$  the  filtration generated by $\{\bar Z^{\sss (n)}(t)\}_{t \ge 0}$. 
For ease of notation, we write $\sup_{t \in [a,b]}=\sup_{[a,b]}$.  We will suppress $(n)$ in the notation, unless needed.\\
% Throughout the proof we shall deal with $\lambda=0$, the general case follows in the exact same manner. We often hide $n$ in the notation. Let $\{ \bar Z(t)\}_{t \ge 0}$ be the scaled random work and $\hat Z(t) = \bar Z(t) - \inf_{0 \le u \le t } \bar Z_u$ be the reflected version.  and denote $c$ for the largest jump size.
\begin{Proposition}
	\label{prop:tails-bound}
	There exist $\Theta \in (0, \infty)$, events $G^{\sss (n)}$, increasing $\FF_t^{\sss (n)}$-stopping times $1= \sigma^{\sss (n)}_0 < \sigma^{\sss (n)}_1 < ... $, and 
	a real positive sequence $\set{\kappa_i}$ with $\sum_{i=1}^\infty \kappa_i < \infty$, such that the following hold:\\
	(i) For every $i \ge 1$, $\set{\sigma_i^{\sss (n)}}_{n \ge 1}$ is tight.\\
	(ii) For every $i \ge 1$,
	$$ \prob\left(\set{ \sup_{[{\sigma^{\sss (n)}_{i-1}},{\sigma^{\sss (n)}_i}]} \hat Z^{\sss (n)}(t) > 2\Theta + 1 } \cap \set{ \sigma^{\sss (n)}_{i-1} < t^{* \sss (n)} } \cap G^{\sss (n)}\right) \le \kappa_i. $$
	(iii) As $n \to \infty$, $ \prob\set{ \sup_{[\sigma^{*\sss (n)}, \infty)} \hat Z^{\sss (n)}(t) > \Theta; G^{\sss (n)}} \to 0$, where $\sigma^{*\sss(n)}= \inf\set{ \sigma_i^{\sss (n)}: \sigma_i^{\sss (n)} \ge t^{*\sss (n)}}  $.	\\
	(iv) As $n \to \infty$, $\prob(G^{\sss (n)}) \to 1$.
\end{Proposition}
Given Proposition \ref{prop:tails-bound}, the proof of Lemma \ref{lemma:tightness-sup} can be completed as follows.\\ \ \\

{\bf Proof of Lemma \ref{lemma:tightness-sup}:} 
% First, we will simplify the situation. Note that $\prob \set{ G_\epsilon^{\sss (n)}} \to 0$ as defined in Lemma \ref{lemma:mn-an-prop} (i), thus the tightness of $\{\sup_{t\ge 0 }\hat Z^{\sss (n)}(t) \}_{n \ge 0}$ is equivalent to the tightness of $\{{\ind}_{G_\epsilon^{\sss (n)}} \sup_{t\ge 0 }\hat Z^{\sss (n)}(t) \}_{n \ge 0}$. Moreover, instead of taking $\lambda=0$, we can take $\lambda = -\epsilon$, and this is will add an $-\epsilon$ drift to the process $\bar Z(t)$. Thus we can assume the drift of $\bar Z(t)$ satisfies
% $$ A(t) \le -t \mbox{ for all } t \le t^*. $$
% Another simplification is to take $c=1$, since we only need $c$ to be the upper bound of the jump sizes. 

 Fix $\epsilon \in (0,1)$.  Let   $\Theta \in (0, \infty)$, $G^{\sss (n)}$, $\sigma_i^{\sss (n)}$, $\kappa_i$ be as in Proposition \ref{prop:tails-bound}.  
Choose $i_0 > 1$ such that $\sum_{i\ge i_0} \kappa_i \le \epsilon$. Since $\set{\sigma_{i_0-1}^{\sss (n)}}$ is tight, 
 there exists $T \in (0,\infty)$ such that $\limsup_{n \to \infty} \prob\set{ \sigma_{i_0-1}^{\sss (n)} > T } \le \epsilon$. 
Thus for any $M' > 2\Theta + 1$, we have
\begin{align*} \prob\set{\sup_{[1,\infty)}\hat Z^{\sss (n)}(t) > M'} \le & \prob\set{\sup_{[1,T]}\hat Z^{\sss (n)}(t) > M'} + \prob\set{\sigma_{i_0-1}^{\sss (n)} > T } + \prob\set{(G^{\sss (n)})^c} \\
	+& \prob\set{\sup_{[\sigma^{\sss (n)}_{i_0-1},\sigma^{* \sss (n)}]}\hat Z^{\sss (n)}(t) > 2\Theta+1; G^{\sss (n)}} + 
	\prob\set{ \sup_{[\sigma^{*\sss (n)}, \infty)} \hat Z^{\sss (n)}(t) > \Theta; G^{\sss (n)}}. \end{align*}
Taking $\limsup_{n \to \infty}$ on both sides 
$$ \limsup_{n \to \infty}\prob\set{\sup_{[1,\infty)}\hat Z^{\sss (n)}(t) > M'} \le \limsup_{n \to \infty}\prob\set{\sup_{[1,T]}\hat Z^{\sss (n)}(t) > M'} + \epsilon + 0 + \epsilon + 0. $$
Since $\set{\sup_{[1,T]}\hat Z^{\sss (n)}(t)}_{n \ge 1}$ is tight, we have,
$$ \limsup_{M' \to \infty} \limsup_{n \to \infty}\prob\set{\sup_{[1,\infty)}\hat Z^{\sss (n)}(t) > M'}  \le 2\epsilon . $$
Since  $\epsilon > 0$ is arbitrary, the result follows. \qed \\

We now proceed to the proof of Proposition \ref{prop:tails-bound}.  The following lemma is key.\\

\begin{Lemma}
	\label{lemma:mn-an-prop} 
	There are $\{\FF_t^{\sss (n)}\}$ adapted processes $\{A^{\sss (n)}(t)\}$, $\{B^{\sss (n)}(t)\}$ and $\FF_t^{\sss (n)}$-martingale $\{M^{\sss (n)}(t)\}$
	such that\\
	(i) $A^{\sss (n)}(\cdot)$ is a non-increasing function of $t$, a.s.  For all $t \ge 0$, $\bar Z^{\sss (n)}(t) = \int_0^t A^{\sss (n)}(u) du + M^{\sss (n)}(t)$.\\
	%The process can be decomposed by $\bar Z(t) = \int_0^t A(u) du + M_t$, where $M_t$ is a martingale. Then\\
	(ii) For $t \ge 0$, $\langle M^{\sss (n)},M^{\sss (n)} \rangle_t = \int_0^t B^{\sss (n)}(u) du$.\\
	(iii) $\sup_{n \ge 1} \sup_{u \ge 0} B^{\sss (n)}(u) \le 2$.\\
	(iv) With   $  G^{\sss (n)} = \set{ A(t) < -t/2  \mbox{ for all } t \in [1, t^{*(n)}]  } $, $\prob(G^{\sss (n)}) \to 1$ as $n\to\infty$.\\
	% (ii) Write $M^2_t = \int_0^t B(u) du+ M_t^*$, where $M^*_t$ is a martingale and $\int_0^t B(u) du = \langle M,M \rangle_t$ is the predictable quadratic variation. Then there exists $n_0 \ge 0$ such that for all $n \ge n_0$ and $t \ge 0$, we have $\langle M,M \rangle_t \le 2t$.\\
	(v) For any $\alpha \in (0, \infty)$ and  $t > 0$, 
	\begin{equation}
\prob \set{ \sup_{ u \in [0,t]} |M^{\sss (n)}(u)| > \alpha } \le 2\exp\left\{ \alpha \right\} \cdot \exp\left\{ - \alpha \log \left( 1+ \frac{\alpha}{2t} \right) \right\}. \label{eqn:1117}	
	\end{equation}
\end{Lemma}

{\bf Proof:} Recall the notation from Section \ref{sec:conv-in-prod}.  Parts (i) and (ii) are proved in \cite{aldous1997brownian}.  Furthermore, from
Lemma 11 of \cite{aldous1997brownian} it follows that, for  $ t \in [l_{i-1}, l_i)$, writing $Q_2(t) = \sum_{j=1}^i (x_{v(j)})^2 $, we have
\begin{equation*}
%\label{eqn:an-def}
	A(t)  \le \sqrt{\frac{s_2}{s_3}}(-1+q s_2 - q Q_2(t)), \; 
	B(t)  \le q s_2 .
\end{equation*}
Part (iii) now follows on recalling  from \eqref{eq:eq1547} that $qs_2 \le 2$.  
To prove (iv) it suffices to show that
\begin{equation}
\label{eqn:sn-s-to-zero}
	\sup_{t\leq t^*} \left|\frac{s_2}{s_3} Q_2(t) - t\right|\convp 0.
\end{equation}
To prove this we will use the estimate on  Page 832, Lemma 13 of \cite{aldous1997brownian}, which says that for any fixed $\epsilon \in (0,1)$, and $ L \in (0,\infty)$
$$ \prob \set{ \sup_{ t \in [0,L]}\left|\frac{s_2}{s_3} Q_2(t) - t \right| > \epsilon }=O\left(  \frac{L^2 x^*}{s_2}+ \sqrt{\frac{L (x^*)^2 s_2}{s_3}} + \frac{L^2 s_3}{s_2^2}+ \sqrt{\frac{Ls_3}{s_2}} + \frac{s_2^2}{(1-2 Ls_2)^+}
\right). $$
%\todo[inline]{where does the fifth term comes from?}
Note that the first term on the right hand side determine its order when $L \to \infty$.
%\footnote{Note that the first term on the right hand side determine its order when $L \to \infty$.}
Taking $L = t^*$ in the above estimate we see that, since $\vartheta \in (0, 1/2)$, the expression on the right side above goes to $0$ as $n \to \infty$.  This proves
\eqref{eqn:sn-s-to-zero} and thus completes the proof of (iv).
Finally, proof of (v) uses standard concentration inequalities for martingales.  Indeed, recalling that the maximal jump size of
$\bar Z$, and consequently that of $M$, is bounded by $1$ and $\langle M,M \rangle_t \le 2t$, we have from Section 4.13, Theorem 5 of \cite{liptser-mart-book} that, for any fixed $\alpha>0$ and $t>0$, 

$$ \prob\{ \sup_{u \in [0,t]} |M_u| > \alpha \} \le 2 \exp \set{-\sup_{\lambda > 0}\left[\alpha \lambda  - 2t \phi(\lambda)\right] },$$
where $\phi(\lambda) = (e^{\lambda }-1-\lambda ).$ A straightforward calculation shows
$$ \sup_{\lambda > 0 } [\alpha \lambda - 2t \phi(\lambda)] = \alpha \log\left(1 + \frac{\alpha}{2t}\right)  - \left( \alpha - 2t \log\left(1 + \frac{\alpha } {2t}\right) \right) \ge \alpha \log\left(1 + \frac{\alpha}{2t}\right)  - \alpha.$$
The result follows. \qed\\

The bound \eqref{eqn:1117} continues to hold  if we replace $M(u)$ with $M(\tau+u)-M(\tau)$ for any finite stopping time $\tau$. 
From this observation we immediately have the following corollary.
\begin{Corollary} \label{lemma:mart-concentration} 
	Let $M$ be as in Lemma \ref{lemma:mn-an-prop}.  Then,  for any finite stopping time $\tau$:\\
	(i)   $\prob \set{ \sup_{ u \in [0,t]} |M(\tau+u)-M(\tau)| > \alpha } \le 2 e^{-\alpha}$, whenever  $\alpha > 2(e^2-1)t$. \\
	(ii)  $\prob \set{ \sup_{ u \in [0,t]} |M(\tau+u)-M(\tau)| > \alpha } \le 2 (2e/\alpha)^{\alpha} t^{\alpha}$, for all $t >0$ and $\alpha > 0$. 
\end{Corollary}
Part (i) of the corollary  is useful when  $\alpha$ is large and
part (ii)  is useful when  $t$ is small.
Finally we now give the proof of Proposition \ref{prop:tails-bound}.\\

{\bf Proof of Proposition \ref{prop:tails-bound}:} From Lemma \ref{lem:lem1325} (i)  we have that $\hat Z^{\sss (n)}$ converges in distribution to
$\hat W_0$ (recall we assume that $\lambda = 0$) as $n \to \infty$. Let $\{\epsilon_i\}_{i \ge 1}$ be a  positive real sequence bounded by $1$ and fix $\Theta \in (2, \infty)$.  Choice of $\Theta$ and $\epsilon_i$  will be specified later in the proof. 
Let $\sigma^{\sss (n)}_0 < \tau^{\sss (n)}_1 \le \sigma^{\sss (n)}_1 < \tau^{\sss (n)}_2 \le \sigma^{\sss (n)}_2 < ... $ be a sequence of stopping times such that $\sigma^{\sss (n)}_0=1$, and for $i \ge 1$,
\begin{equation}
	\tau^{\sss (n)}_i = \inf \{t \ge \sigma^{\sss (n)}_{i-1} + \epsilon_i : \hat Z^{\sss (n)}(t)\ge \Theta\} \wedge (\sigma^{\sss (n)}_{i-1}+1), \;\; \sigma^{\sss (n)}_i = \inf \{t \ge \tau^{\sss (n)}_{i}  : \hat Z^{\sss (n)}(t) \le 1\}. \label{eqn:def-stopping-times}
\end{equation}
Similarly define stopping times 
$1 = \bar \sigma_0 < \bar \tau_1 \le \bar \sigma_1 < \bar \tau_2 \le \bar \sigma_2 < ... $ 
by replacing $\hat Z^{\sss (n)}$ in \eqref{eqn:def-stopping-times} with $\hat W_0$.  Due to the negative quadratic drift in the definition of
$W_0$ it follows that $\bar \sigma_i < \infty$ for every $i$ and from the
weak convergence of $\hat Z^{\sss (n)}$ to $\hat W_0$ it follows that $\sigma^{\sss (n)}_i \to  \bar \sigma_i $ and 
$\tau^{\sss (n)}_i \to  \bar \tau_i $, in distribution, as $n \to \infty$.  Here we have used the fact that if $\zeta$ denotes the first time $W_0$
hits the level $\alpha \in (0, \infty)$ then, a.s., for any $\delta > 0$, there are infinitely many crossings of the level $\alpha$ in $(\zeta, \zeta + \delta)$.
In particular we have that $\{\sigma^{\sss (n)}_i\}_{n \ge 1}$ is a tight sequence, and this proves part (i) of Proposition \ref{prop:tails-bound}.\\

For the rest of the proof we suppress $(n)$ from the notation.  Since the jump size of $\hat Z$ is bounded by $1$, we have that  $  \sup_{[{\sigma_{i-1}},{\sigma_{i-1}+ \epsilon_i}]} \hat Z(t) \le \Theta  $ implies $ \sup_{[\sigma_{i-1},\tau_i]} \hat Z(t) \le \Theta+1 $ and  thus, in this case, when $t \in [\tau_i,\sigma_i]$, we have $\hat Z(t) = \hat Z(\tau_i) + ( \bar Z(t) - \bar Z(\tau_i) )\le \Theta+1 +( \bar Z(t) - \bar Z(\tau_i))$.
  Let $G \equiv G^{\sss (n)}$ be as in Lemma \ref{lemma:mn-an-prop} (iv) and let
$H_i = G \cap \{\sigma_{i-1} < t^*\}$, then writing $\prob (\cdot \cap H_i)$ as $\prob_i(\cdot)$,
\begin{align}
	\prob_i \set{ \sup_{[{\sigma_{i-1}},{\sigma_i}]} \hat Z(t) > 2\Theta + 1 } 
	\le& \prob_i \set{ \sup_{[{\sigma_{i-1}},{\sigma_{i-1}+ \epsilon_i}]} \hat Z(t) > \Theta }\\
	 +& \prob_i \set{ \sup_{[\tau_i, \sigma_i]} \left[\Theta + 1 + (\bar Z(t) - \bar Z(\tau_i))\right] > 2\Theta+1 }.
\end{align} 
Denote the two terms on the right side by $\TTT_1$ and $\TTT_2$ respectively. Recalling that $\hat Z (\sigma_{i-1})\le 2$, we have 
from the decomposition in Lemma \ref{lemma:mn-an-prop} (i) and Corollary \ref{lemma:mart-concentration}(ii) that
 % 
 % we have the decomposition $\bar Z(t) = \int_0^t A(u)du + M_t$. Then $\sup_{[{\sigma_{i-1}},{\sigma_{i-1}+ \epsilon_i}]} \hat Z(t) > \Theta$ implies that there exist two random times $\sigma_{i_1} \le \alpha \le \beta \le \sigma_{i-1}+ \epsilon_i$ such that $ \bar Z_{\beta} - \bar Z_\alpha > \Theta$, thus $ M_\beta - M_\alpha > \Theta$ and this implies
\begin{align}
	\TTT_1 \le \prob\set{ \sup_{[\sigma_{i-1},\sigma_{i-1}+\epsilon_i]} |M(t)-M(\sigma_{i-1})| > \frac{\Theta-2}{2} } \le  C_{\frac{\Theta-2}{2}} {\epsilon_i}^{(\Theta-2)/2}, \label{eqn:1167}
\end{align}
Here, for $\alpha > 0$, $C_{\alpha} = 2 (2e/\alpha)^{\alpha}$ and we have used the fact that on $H_i$,
$A(t) \le -t/2\le 0$ for all $t \in [\sigma_{i-1}, \sigma_{i-1}+ \epsilon_i]$.
%where the last bound and the constant $C_{\Theta/2}$ is due to Lemma \ref{lemma:mart-concentration} (ii).\\

Next, let $\{\delta_i\}_{i\ge 1}$ be a sequence of positive reals bounded by $1$.  Setting $d_i = \sum_{j=1}^{i-1} \epsilon_i$, we have
% In order to bound $P_2$, let $\set{d_i}_{i\ge 1}$ be a deterministic sequence such that $A(t) \le -d_i$ for all $t \in [\tau_i,\sigma_i]$. Since $\tau_i \ge \sum_{j=1}^{i-1}\epsilon_i $ and $A(t) \le -t$, we will take $d_i = \sum_{j=1}^{i-1}\epsilon_i$ later in the proof of part (ii).  Let $\{\delta_i\}_{i\ge 1}$ be another deterministic sequence such that $\delta_i < 1$, whose values will also be specified later. Thus
\begin{align}
	\TTT_2 
	\le& \prob_i \set{ \sup_{[\tau_i,\tau_i+\delta_i]}(\bar Z(t) - \bar Z(\tau_i)) > \Theta } + \prob_i \set{\sup_{[\tau_i + \delta_i, \tau_i + 1]}(\bar Z(t) - \bar Z(\tau_i)) > \Theta} + \prob\set{ \sigma_i > \tau_i+1} \nonumber \\
	\le& \prob \set{ \sup_{[\tau_i,\tau_i+\delta_i]}(M(t) - M(\tau_i)) > \Theta } + 
	\prob \set{\sup_{[\tau_i + \delta_i, \tau_i + 1]}(M(t) - M(\tau_i)) > \Theta + \frac{\delta_i d_i}{2}} \nonumber\\
	+& \prob\set{ M(\tau_i + 1)-M(\tau_i) > -\Theta + \frac{d_i}{2}} \nonumber\\
	\le& C_{\Theta}\delta_i^{\Theta} + 2 e^{-\delta_i d_i/2} + 2 e^{\Theta-d_i/2}, \label{eqn:1174}
\end{align} 
whenever 
\begin{equation}\label{eq:ins431} \min\{\delta_i d_i/2 , d_i/2 - \Theta\} > 2(e^2-1).\end{equation}
 Fix $\Theta > 14$.  Then  $\max\set{C_{\Theta}, C_{(\Theta-2)/2}}\le 2$. We will impose additional conditions on $\Theta$ later in the proof.   Combining \eqref{eqn:1167} and \eqref{eqn:1174}, we have
\begin{align} 
	\prob_i \set{ \sup_{[{\sigma_{i-1}},{\sigma_i}]} \hat Z(t) > 2\Theta + 1 } \le 2 ( {\epsilon_i}^{(\Theta-2)/2} + \delta_i^{\Theta} + e^{-\delta_i d_i/2} + e^{\Theta-d_i/2}) \equiv \kappa_i. \label{eqn:basic-bound}
\end{align}
%Define the right-hand-side of \eqref{eqn:basic-bound} to be $\kappa_i$. Now we specify our choice of those values $\epsilon_i, d_i$ and $\delta_i$ as follows:
Let
$$ \epsilon_i= i^{-1/2}, \;\; d_i = \sum_{j=1}^{i-1} \epsilon_i \sim i^{1/2} , \;\; \delta_i = 1/\sqrt{d_i} \sim i^{-1/4} .$$
Then, \eqref{eq:ins431} holds for $i$ large enough, and 

$$\kappa_i \sim 2 ( i^{-(\Theta-2)/4}+ i^{-\Theta/4} + e^{-i^{1/4}/2} + e^\Theta  e^{-i^{1/2}/2}),$$
which, since  $\Theta >14$, is summable. 
%Note that we require $\inf_{[\sigma_{i-1},\sigma_i]} A(t) \le -d_i$, and this is true when $\sigma_{i-1} < t^*$, thus condition (ii) of the proposition is proved.\\
This proves part (ii) of the Proposition.\\

Now we consider part (iii).  
We will construct another sequence of stopping times with values in $[t^*, \infty)$, as follows.
 Define
 $\sigma_0^* :=\inf\set{ \sigma_i : \sigma_i \ge t^*} = \inf\set{t \ge t^*: \hat Z(t) \le 1 }$, then define $\tau^*_i, \sigma^*_i$ for $i \ge 1$ similarly as in \eqref{eqn:def-stopping-times}. Similar arguments as before give a bound as \eqref{eqn:basic-bound} with $d_i$ replaced by $t^*$, $\delta_i$ replaced by
$1/\sqrt{t^*}$, $\epsilon_i$ replaced with $1/t^*$ and $\Theta$ replaced by any $\Theta_0 > 14$. Namely,
\begin{align} 
	\prob \set{ \sup_{[{\sigma_{i-1}^*},{\sigma_i^*}]} \hat Z(t) > 2\Theta_0 + 1;\; G^{\sss(n)} } \le 2 ( (1/t^*)^{(\Theta_0-2)/2} + (1/\sqrt{t^*})^{\Theta_0} + e^{-\sqrt{t^*}/2 } + e^{\Theta_0-t^*/2}) . \label{eqn:basic-boundnew}
\end{align}
Here we have used the fact that
 since $A(t)$ is non-increasing, on $G^{\sss(n)}$, $A(t) \le -t^*/2$ for all $t \ge t^*$.\\

%and \eqref{eqn:basic-bound} is indeed true when replacing $\sigma_i$ with $\sigma^*_i$.\\
Recall that, by construction, 
 $\hat Z(t) = 0$ when $t \ge s_1$. So there exist $i_0$ such that $\tau^*_{i_0}= \infty$, in fact since $\sigma^*_i \ge \sigma^*_{i-1}+ \epsilon$, we have that $i_0 \le s_1/\epsilon$.  Thus, we have from the above display that
$$ \prob\set{ \sup_{[\sigma^{*}_0, \infty)} \hat Z(t) > 2\Theta_0+1} \le \frac{2s_1}{\epsilon}  ( (1/t^*)^{(\Theta_0-2)/2} + (1/\sqrt{t^*})^{\Theta_0} + e^{-\sqrt{t^*}/2 } + e^{\Theta_0-t^*/2}).$$ 
 Taking $\Theta > 29$, we have on setting $\Theta_0 = \frac{\Theta -1}{2}$ in the above display

$$ \prob\set{ \sup_{[\sigma^{*}_0, \infty)} \hat Z(t) > \Theta} \le 2s_1 \left( {\left(\frac{1}{t^*}\right)}^{(\Theta-1)/4-2} + {\left(\frac{1}{t^*}\right)}^{(\Theta-1)/4-1} + \frac{1}{t^*} e^{-\sqrt{t^*}/2} + \frac{1}{t^*}e^{(\Theta-1)/2-t^*/2}\right).$$

From \eqref{eqn:additional-condition} we have that  $s_1 \cdot (\frac{1}{t^*})^{\varsigma/\vartheta} \to 0$.
So if $\Theta \ge 4 (\frac{\varsigma}{\vartheta} +2) + 1$, the above expression approaches $0$ as $n \to \infty$.
The result now follows on taking $\Theta = \max\{29, 4 (\frac{\varsigma}{\vartheta} +2) + 1\}$. \qed \\

\section{Bounded-size rules at time $t_c- n^{-\gamma}$}
\label{sec:main-bsr-susceptibility}

 Throughout Sections \ref{sec:main-bsr-susceptibility} and \ref{sec:main-coupling} we take $T = 2t_c$ which is a convenient upper bound for the time parameters of interest.
In this section we prove Theorems   \ref{thm:suscept-funct} and \ref{thm:suscept-limit}.

% After defining the standard augmented multiplicative coalescent process properly, in this section and the next, we will study the BSR process and eventually complete the proof of Theorem \ref{thm:crit-regime}.
% 
% The goal of this section is to analyze the component sizes and susceptibility of $\BS^{\sss(n)}(t)$ at time $t_n =t_c-n^{\gamma}$ for $\gamma \in (1/6,1/5)$, thus verifying the conditions in Theorem \ref{theo:aldous-full-gene} for the proper rescaled component sizes of $\BS^{\sss(n)}(t_c-n^{\gamma})$. We will also prove Theorem \ref{thm:suscept-limit} and Theorem \ref{thm:suscept-funct} along the way. Then in Section \ref{sec:main-coupling} we will construct a coupling between $\BS^{\sss(n)}(t)$ and an \erdos process through the critical window, thus applying Theorem \ref{theo:aldous-full-gene} and completing the proof of Theorem \ref{thm:crit-regime}.\\

We begin with some notation associated with BSR processes, which closely follows  \cite{spencer2007birth}.
Recall from Section \ref{sec:bsr} the set $\Omega_K$ and the random graph process $\BS^{\sss (n)}(t)$ associated  with a given $K$-BSR $F \subset \Omega_K^4$. Frequently we will suppress $n$ in the notation. 
Also recall the definition of $c_t(v)$ from Section \ref{sec:bsr}.

%that for $v \in \VV$, $c(v) = \CC(v)$ if $\CC(v) \le K$ and $c(v)=\ompar $ otherwise. 
For $i \in \Omega_K$, define 
\be
X_i(t)=|\{v \in \BS_t^{\sss (n)} : c_t(v)=i  \}| \mbox{ and } \barx_i(t)=X_i(t)/n. \label{eqn:def-x}
\ee
% Denote $\CC^{\sss(n)}_i(t)$ for the size of the $i$-th largest component of $\BS^{\sss(n)}(t)$, then define the susceptibilities $\calS_k(t)$ and $\bar s_k(t)$, $k\ge 1$, as follows
% $$ \calS_k(t) = \sum_{i=1}^\infty (\CC_i^{\sss(n)}(t) )^k \mbox{ and } \bar s_k(t) = \calS_k(t)/n.$$
Denote by $\BS^*(t)$ the subgraph of $\BS(t)$ consisting of all components of size greater than $K$, and define, for $k=1,2,3$
$$ \calS_{k,\ompar}(t):= \sum_{\{\CC \subset \BS^*(t)\}} |\CC|^k \mbox{ and } \bar s_{k,\ompar}(t) = \calS_{k,\ompar}(t)/n,$$
where $\{\CC \subset \BS^*(t)\}$ denotes the collection of all components in $\BS^*(t)$. For notational convenience in long formulae, we sometimes write $\BS(t) = \BS_t$ and similarly $\BS^*(t) = \BS_t^*$.  Similar notation will be used throughout the paper.

Clearly
\be
\calS_k(t) = \calS_{k,\ompar} + \sum_{i=1}^K i^{k-1} X_i(t), \; \bar s_k(t) = \bar s_{k,\ompar} + \sum_{i=1}^K i^{k-1} \bar x_i(t). \label{eqn:s2omega-to-s2}
\ee
Also note  that $\calS_1(t) = n$ and $\calS_{1,\ompar}(t) = X_\ompar(t)$.\\

Recall the Poisson processes $\PP_{\vec{v}}$ introduced in Section \ref{sec:bsr}.  Let
$\FF_t = \sigma \{\PP_{\vec{v}}(s): s \le t, \vec{v} \in [n]^4\}$.
For $T_0 \in [0,T]$ and a $\{\FF_t\}_{0\le t < T_0}$ semimartingale $\{J(t)\}_{0 \le t < T_0}$ of the form 
 \begin{equation} \label{eq:semimart} dJ(t)  = \alpha(t) dt + dM(t), \langle M, M \rangle_t = \int_0^t \gamma(s) ds,  \end{equation}
where $M$ is a $\{\clf_t\}$ local martingale and $\gamma$ is a progressively measurable process,  we write $\alpha = \bfd(J)$, $M = \bfm(J)$
and $\gamma = \bfv(J)$.

% For any two pure jump processes $J_1(t)$ and $J_2(t)$, denote  $\E[ \Delta J_1(t) \Delta J_2(t)| \FF_t]$ for the $\FF_t$-adapted process $g(t)$ such that
% $$ \langle J_1, J_2 \rangle_t = \int_0^t g(s)ds. $$
% So saying $\E[ \Delta J(t) | \FF_t] = A(t)$ and $\E[ (\Delta J(t))^2 | \FF_t] = B(t)$ is equivalent to saying that there exists decomposition
% $$ J(t) = \int_0^t A(s) ds + M(t), $$
% such that $A(t)$ is $\FF_t$-adapted and $M(t)$ is a $\FF_t$-martingale with 
% $$\langle M \rangle(t) = \int_0^t B(s)ds.$$
% 
% \todo[inline]{Shall we apply this notation to all the rest of the analysis? In that case, the analysis of $Y(t)$ and $Z(t)$ in Section 7.4 and 7.5 can also be written in this form and looks more compact.---Xuan}

{\bf Organization:} Rest of this section is  organized as follows.
%\begin{itemize}
 In Section \ref{sec:prelim-bsr}, we state a recent result on BSR models and certain deterministic maps associated with the evolution of $\BS_t^*$ from \cite{bsr-2012} that will be used in this work.
	 In Section \ref{sec:bsr-diff-s2s3A}, we will study the asymptotics of $\bar s_{2,\ompar}$ and $\bar s_{3,\ompar}$.
	% analyze the susceptibilities $\bar s_{2,\ompar}(t)$ and $\bar s_{3,\ompar}(t)$, in preparation of proving Theorem \ref{thm:suscept-funct}.
In Section \ref{sec:proof-thm-alphabeta}, we will complete the proof of Theorem \ref{thm:suscept-funct}.
 In Section \ref{sec:decompose-s2-s3}, we will obtain some useful semimartingale decompositions for certain functionals of $\bar s_2$ and $\bar s_3$.
	%further give more information about $\bar s_2(t)$, $\bar s_3(t)$ and the related ratios shown up in Theorem \ref{thm:suscept-limit}.
 In Section \ref{sec:proof-conv-susceptibility}, we will complete the proof of Theorem \ref{thm:suscept-limit}.
%\end{itemize}

\subsection{Evolution of  $\BS^*_t$.}
\label{sec:prelim-bsr}
We begin with the following lemma from   \cite{bsr-2012} (see also 
\cite{spencer2007birth}). 
\begin{Lemma}
	\label{lemma:approx-xi}\ 
(a) For each $i \in \Omega_K$, there exists a continuously differentiable function $x_i:[0,T] \to [0,1]$ such that  for any  $\delta \in (0,1/2)$, there exist $C_1, C_2 \in (0, \infty)$ such thar
 for all $n$, 
		\[\prob\left( \sup_{i \in \Omega_K} \sup_{s \in [0,T]} |\barx_i(t)-x_i(t) | > n^{-\delta}\right) <C_1 \exp\left( -C_2 n^{1-2\delta}\right).\]
(b) There exist polynomials $\{F^x_i(\bfx)\}_{i\in \Omega_K}$,  $\bfx=(x_i)_{i \in \Omega_K} \in \RRR^{K+1}$, such that  $\bfx(t)= (x_i(t))_{i \in \Omega_K}$  is the unique solution to the differential equations:
\begin{equation}
x_i^\prime(t)= F^x_i(\bfx(t)), \;\; i \in \Omega_K, \;\; t \in [0,T]	 \mbox{ with initial values } \bfx(0)=(1,0,..., 0).	\label{eqn:sys-dif-eqns}
\end{equation}
Furthermore, $\bar x_i$ is a $\{\clf_t\}_{0\le t < T}$ semimartingale of the form \eqref{eq:semimart} and
$$\sup_{0\le t < T} |\bfd(\bar x_i)(t) - F_i^x(\bar \bfx(t)) | \le \frac{K^2}{n}.$$
Also,
for all $i \in \Omega_K$ and $t \in  (0,T]$, we have $x_i(t)>0$ and $ \sum_{i\in \Omega_K} x_i(t) = 1$. \\
\end{Lemma}

Recall that $\BS^*(t)$ is the subgraph of $\BS(t)$ consisting of all components of size greater than $K$.   The evolution of this graph is governed by three type of events:\\

{\bf Type 1 (Immigrating vertices): } This corresponds to the merger of two components of size bounded by $K$ into a component of size larger than $K$.  Such an event leads to the
appearance of a new component in $\BS^*(t)$ which we view as the immigration of a `vertex' into $\BS^*(t)$.  Denote by  $n a^*_i(t)$ the  rate at which a component of size $K+i$  immigrates into $\BS_t^*$ at time $t$.
In \cite{bsr-2012} it is shown that there are polynomials $F_i^a(\bfx)$ for $1\le i\le K$ such that, with
$\bar \bfx(t) = (\bar x_i(t))_{i \in \Omega_K}$
\be
\sup_{t \in [0,\infty) }| a^*_i(t)-F^a_i(\bar \bfx(t))| \le \frac{K}{n}.	\label{eqn:error-a}
\ee
We define, with $\bfx(t)$ as in Lemma \ref{lemma:approx-xi},
\be
a_i(t) := F^a_i(\bfx(t)), \; i = 1, \cdots K. \label{eqn:def-fa-a}
\ee

{\bf Type 2 (Attachments):} This event corresponds to a component of size at most $K$ getting linked with some component of size larger than $K$.  For $1\le i \le K$, denote by $|\clc| c^*_i(t)$ the rate at which a component of size $i$ attaches to  a  component $\clc$ in $\BS^*_{t-}$.  
  Then (see \cite{bsr-2012}) there exist polynomials  $F^c_i(\bfx)$ for $1\le i\le K$,  such that 
$c_i^*(t) = F^c_i(\bar \bfx(t))$.  Define
\be
 c_i(t) := F^c_i(\bfx(t)), i = 1, \cdots K. \label{eqn:def-fc-c}
\ee

{\bf Type 3 (Edge formation):} This event corresponds to the addition of an edge between  components in $\BS^*_t$.  The occurrence of this event adds one edge between two vertices in $\BS^*_{t-}$,  the vertex set stays unchanged, whereas the edge set
has one additional element.  From \cite{bsr-2012}, there is a polynomial $F^b(\bfx)$ such that, defining $b^*(t) = F^b(\bar \bfx(t))$, the rate at which each pair of components $\CC_1 \neq \CC_2 \in \BS_t^*$ merge at time $t$, equals  $|\CC_1|   |\CC_2| b^*(t)/n$.
Furthermore 
\be
b(t) :=  F^b(\bfx(t)) \label{eqn:def-fb-b}
\ee
 satisfies $b(t_c) \in (0, \infty)$.

% One last comment is that, as part of the proof of Theorem \ref{lemma:approx-xi}, also see \cite{spencer2007birth}, we have
% $$ |\E[\Delta \bar x_i(t)| \FF_t] - F^x_i(\bar \bfx(t))| \le K^2/n. $$
% This $O(n^{-1})$ error along with the one in \eqref{eqn:error-a} have trivial contribution in the analysis through this section, so we will just ignore the $O(n^{-1})$ error and assume $\E[\Delta \bar x_i(t)| \FF_t] = F^x_i(\bar \bfx(t))$ and $a^*_i(t) = F^a_i(\bar \bfx(t))$.

\subsection{Analysis of $\bar s_{2,\ompar}(t)$ and $\bar s_{3,\ompar}(t)$}
\label{sec:bsr-diff-s2s3A}

We begin by recalling a result from \cite{spencer2007birth}.
Define  functions $F^s_{2,\ompar} : [0,1]^{K+1}\times \RRR \to \RRR$ 
and $F^s_{3,\ompar} : [0,1]^{K+1}\times \RRR^2 \to \RRR$
as 
\be
F^s_{2,\ompar}(\bfx, s_2) 
:= \sum_{j=1}^K (K+j)^2 F^a_{j}(\bfx) + y_2\sum_{j=1}^K 2j F^c_{j}(\bfx)  + x_\ompar \sum_{j=1}^K j^2 F^c_{j}(\bfx) + (s_2)^2  F^b (\bfx) , 
\label{eqn:fs2omega}
\ee
for $(\bfx, s_2) \in [0,1]^{K+1}\times \RRR $ and, for $(\bfx, s_2, s_3) \in [0,1]^{K+1}\times \RRR^2 $
\begin{align}
F^s_{3,\ompar}(\bfx, s_2, s_3) 
:=& \sum_{j=1}^K (K+j)^3 F^a_{j}(\bfx) + s_3 \sum_{j=1}^K 3jF^c_{j}(\bfx) + 3s_2 \sum_{j=1}^K j^2 F^c_{j}(\bfx)\nonumber\\
 +& x_\ompar \sum_{j=1}^K  j^3 F^c_j(\bfx) +  3s_2 s_3  F^b(\bfx). \label{eqn:fs3omega}
\end{align}
\begin{Lemma}\cite{spencer2007birth}
	\label{lemma:s2omega-s3omegaA}
	As $n \to \infty$, $\bar s_{j,\ompar}(t) \to s_{j,\ompar}(t)$ in probability, $j=2,3$, $t \in [0, t_c)$.  Furthermore, $s_{j,\ompar}$ are continuously
	differentiable on $[0, t_c)$ and can be characterized as the unique solutions of the equations
	\begin{align}
		s'_{2,\ompar}(t) =& F^s_{2,\ompar}(\bfx(t), s_{2,\ompar}(t)),\; s_{2,\ompar}(0)=0, \label{eqn:diff-s2}\\
		s'_{3,\ompar}(t) =& F^s_{3,\ompar}(\bfx(t), s_{2,\ompar}(t), s_{3,\ompar}(t)),\; s_{3,\ompar}(0)=0. \label{eqn:diff-s3}
	\end{align}
 Furthermore,
	$\lim_{t \to t_c} s_{2,\ompar}(t) = \lim_{t \to t_c} s_{3,\ompar}(t) = \infty.$
\end{Lemma}

% In this section, we will analysis $\bar s_{2,\ompar}(t)$ and $\bar s_{3,\ompar}(t)$ closely propose the limit objects for $\lim_{n \to \infty}\bar s_{2,\ompar}(t)$ and $\lim_{n \to \infty}\bar s_{3,\ompar}(t)$, denote by $s_{2,\ompar}(t)$ and $s_{3,\ompar}(t)$ respectively, which are solutions of certain differential equations. These differential equations has been studied in \cite{spencer2007birth}. From $s_{2,\ompar}(t)$ and $s_{3,\ompar}(t)$, one can easily get the differential equations satisfied by $s_2(t)$ and $s_3(t)$.\\
% 
% The next lemma analyze the trend part of $\bar s_{2,\ompar}(t)$ and $\bar s_{3,\ompar}(t)$.
The following lemma gives additional information on the convergence of $\bar s_{j,\ompar}$ to $s_{j,\ompar}$.  For $T_0 \in [0, T]$, a stochastic process
$\{\xi(t)\}_{0 \le t < T_0}$, and a nonnegative sequence $\alpha(n)$, the quantity $O_{T_0}(\xi(t)\alpha(n))$ will represent a stochastic process
$\{\eta(t)\}_{0\le t < T_0}$ such that for some $d_1 \in (0, \infty)$, $\eta(t) \le d_1 \xi(t)\alpha(n)$,  for all $0 \le t < T_0$ and $n \ge 1$.

\begin{Lemma}
	\label{lemma:s2omega-s3omega}
	The processes $\bar s_{j,\ompar}$, $j=2,3$, are $\{\clf_t\}_{0\le t < t_c}$ semimartingales of the form \eqref{eq:semimart} and
	\begin{align*}
	&|\bfd(\bar s_{2,\ompar})(t) - F^s_{2,\ompar}(\bar \bfx(t), \bar s_{2,\ompar}(t))| =  O_{t_c}(\calS_4(t)/n^2)\\
	&|\bfd(\bar s_{3,\ompar})(t) - F^s_{3,\ompar}(\bar \bfx(t), \bar s_{2,\ompar}(t), \bar s_{3,\ompar}(t))|=  O_{t_c}(\calS_5(t)/n^2).
	\end{align*}	
\end{Lemma}

{\bf Proof:} Note that $\calS_{2, \ompar}$ and $\calS_{3, \ompar}$ have jumps at time instant $t$  with rates and values $ \Delta \calS_{2, \ompar}(t)$, $\Delta \calS_{3, \ompar}(t)$, respectively, given as follows.
% 
% 
% For (a), the existence of solution and explosion at $t_c$ has been proved in \cite{spencer2007birth}, so we only prove (b) and derive the expression of the two functions $F^s_{2,\ompar}$ and $F^s_{3,\ompar}$.\\
% For $k=2,3$ denote $\Delta \calS_{k,\ompar}$ for the jump size of $\calS_{k,\ompar}(t)$ when one edge is added to $\BS_t^*$, then according to the three cases in Section \ref{sec:prelim-bsr} we have

\begin{itemize}
	\item  for each $1 \le i \le K$, with rate $n  a^*_i(t)$, 
	$$ \Delta \calS_{2, \ompar}(t) = (K+i)^2, \;\; \Delta \calS_{3, \ompar}(t) = (K+i)^3.$$
	\item  for each $1 \le i \le K$ and $\CC \subset \BS_{t-}^*$, at rate $|\CC|  c^*_i(t)$,
	$$ \Delta \calS_{2, \ompar}(t) = 2|\CC| i+ i^2, \;\; \Delta \calS_{3, \ompar}(t) = 3|\CC|^2 i + 3 |\CC| i^2 +i^3. $$
	\item  for all unordered pair  $\CC,\tilde \CC \subset \BS_{t-}^*$, such that   $\CC \neq \tilde \CC$, at rate $ |\CC|  |\tilde\CC|   b^*(t)/n$, 
	$$ \Delta \calS_{2, \ompar}(t) = 2|\CC| |\tilde \CC| , \;\; \Delta \calS_{3, \ompar}(t) = 3|\CC|^2|\tilde \CC| + 3 |\CC| |\tilde \CC|^2. $$
\end{itemize}
Thus 
\begin{align}
	\bfd(\calS_{2,\ompar})(t)
	=& \sum_{j=1}^K (K+j)^2 n a^*_j(t) + \sum_{j=1}^K \sum_{\CC \subset \BS_t^*} (2j|\CC| + j^2) |\CC| c^*_j(t) + \sum_{\CC \neq \tilde \CC \subset \BS_t^*} 2|\CC| |\tilde \CC|  \frac{b^*(t) |\CC| |\tilde \CC|}{n} \nonumber  \\
	=& \sum_{j=1}^K (K+j)^2 n a^*_j(t) + \sum_{j=1}^K 2j c^*_j(t) \calS_{2,\ompar}(t)
	 +  \sum_{j=1}^K j^2 c^*_j(t) X_\ompar(t)\nonumber\\ +& \frac{b^*(t)}{n} (\calS_{2,\ompar}^2(t) - \calS_{4,\ompar}(t))\nonumber\\
	=& n\left(F_{2,\ompar}^s(\bar \bfx, \bar s_{2,\ompar}) + O(1/n) + O_{t_c}(S_{4,\ompar}(t)/n^2)\right) \label{eqn:2072}
\end{align}
and
\begin{align*}
	&\bfd(\calS_{3,\ompar})(t) \nonumber\\
	=& \sum_{j=1}^K (K+j)^3 n a^*_j(t) + \sum_{j=1}^K \sum_{\CC \subset \BS_t^*} (3j|\CC|^2 + 3 j^2 |\CC| + j^3) |\CC| c^*_j(t) \\
	&+ \sum_{\CC \neq \tilde \CC \subset \BS_t^*} (3|\CC|^2 |\tilde \CC| + 3|\CC| |\tilde \CC|^2 ) \frac{b^*(t) |\CC| |\tilde \CC|}{n} \\
	=& \sum_{j=1}^K (K+j)^3 n a^*_j(t) + \sum_{j=1}^K 3j c^*_j(t) \calS_{3,\ompar}(t)  
	+ \sum_{j=1}^K 3j^2 c^*_j(t) \calS_{2,\ompar}(t) \\ +& \sum_{j=1}^K j^3 c^*_j(t) X_\ompar(t) + \frac{3b^*(t)}{n} (\calS_{3,\ompar}(t)\calS_{2,\ompar}(t) - \calS_{5,\ompar}(t))\\
	=& n\left(F_{3,\ompar}^s(\bar \bfx(t), \bar s_{2,\ompar}(t), \bar s_{3,\ompar}(t)) + O(1/n) + O_{t_c}(S_{4,\ompar}(t)/n^2)\right).
\end{align*}
% Note that by the nature of $\calS_{k,\ompar}(t)$ and $\bar s_{k,\ompar}(t)$, we have $\E[\Delta \bar s_{k,\ompar}(t)|\FF_t] = \E[\Delta \calS_{k,\ompar}(t)|\FF_t]/n$.
% 
% Then by \eqref{eqn:2072} and \eqref{eqn:error-a} and the fact $\calS_4(t)/n \ge 1$ we have
% \begin{align}
% 	|\E[\Delta \bar s_{2,\ompar}(t)|\FF_t] - F^s_{2,\ompar}(\bar \bfx(t), \bar s_{2,\ompar}(t))| 
% 	\le& \sum_{j=1}^K (K+j)^2|a^*_j(t)-a_j(t)| + \frac{b^*(t)\calS_{4,\ompar}(t)}{n^2} \nonumber\\	
% 	\le& \sum_{j=1}^K (K+j)^2\cdot \frac{K}{n} + \frac{b^*(t)\calS_{4}(t)}{n^2}
% 	= O\left(\frac{\calS_{4}(t)}{n^2}\right).
% \end{align}
% The other bound for $\E[\Delta \bar s_{3,\ompar}(t)|\FF_t]$ can be established similarly.
The result follows.  \qed \\

% The key results in the above theorem can be summarized as
% \begin{align*}
% 	\E[\Delta \bar s_{2,\ompar}(t)|\FF_t] 
% 	&\approx F^s_{2,\ompar}(\bar \bfx(t), \bar s_{2,\ompar}(t)) \approx F^s_{2,\ompar}(\bfx(t), s_{2,\ompar}(t)) =s'_{2,\ompar}(t)\\
% 	\E[\Delta \bar s_{3,\ompar}(t)|\FF_t] 
% 	&\approx F^s_{3,\ompar}(\bar \bfx(t), \bar s_{2,\ompar}(t), \bar s_{3,\ompar}(t)) \approx F^s_{3,\ompar}(\bfx(t), s_{2,\ompar}(t), s_{3,\ompar}(t)) =s'_{3,\ompar}(t).	
% \end{align*}
% \ \\
As an immediate consequence  of \eqref{eqn:s2omega-to-s2} and the convergence of $(\bar s_k, \bar s_{k,\ompar}, \bar \bfx)$
to $( s_k, s_{k,\ompar}, \bfx)$
we have the following formula.
\be
s_k(t) := s_{k,\ompar}(t) + \sum_{i=1}^K i^{k-1} x_i(t), \mbox{ for } k=2,3. \label{eqn:2134}
\ee
This, along with Lemma \ref{lemma:s2omega-s3omegaA} and Lemma \ref{lemma:approx-xi}(b), yields the following differential equations for $s_2$ and $s_3$.
\begin{Lemma}
	\label{lemma:diff-s2s3}
	The functions $s_2,s_3$ are continuously differentiable on $[0, t_c)$ and can be characterized as the unique solutions of the following differential equations
	\begin{align*}
		s'_2(t) &= F^s_2(\bfx(t), s_2(t)),  \;\;\;\;  s_2(0)=1,\\
		s'_3(t) &= F^s_3(\bfx(t), s_2(t) ,s_3(t)), \;\; \;\;  s_3(0)=1.
	\end{align*}
	where the function $F^s_2(\cdot)$ and $F^s_3(\cdot)$ are defined as
	\begin{align*}
		 F^s_2(\bfx,s_2) &:= F^s_{2,\ompar}\left(\bfx, s_2 - \sum_{i=1}^K i x_i\right) + \sum_{i=1}^K i F^x_{i}(\bfx),\\
		F^s_3(\bfx,s_2,s_3) &:= F^s_{3,\ompar}\left(\bfx,\; s_2- \sum_{i=1}^K i x_i,\; s_3-\sum_{i=1}^K i^2 x_i\right) + \sum_{i=1}^K i^2 F^x_{i}(\bfx).
	\end{align*}	
\end{Lemma}

\subsection{Proof of Theorem \ref{thm:suscept-funct}}
\label{sec:proof-thm-alphabeta}
In this section we prove  Theorem \ref{thm:suscept-funct}. We begin with the following lemma which defines the two parameters $\alpha$ and $\beta $ that appear in Theorems \ref{thm:suscept-funct} and  \ref{thm:crit-regime}.
Recall from Section \ref{sec:prelim-bsr} that $b(t_c) \in (0, \infty)$.
\begin{Lemma} 
	\label{lemma:alpha-beta}
	There following two limits exist, 
	$$ \alpha := \lim_{t \to t_c-} (t_c-t)s_2(t), \;\; \beta := \lim_{t \to t_c-} \frac{s_3(t)}{(s_2(t))^3}.$$
	Furthermore, $\alpha, \beta \in (0, \infty)$ and $\alpha = 1/{b(t_c)} $.
\end{Lemma}
{\bf Proof:} By \eqref{eqn:2134}, for $k = 2,3$, $|s_k(t)-s_{k,\ompar}(t)| \le K^{k}$. Since $s_k(t) \to \infty$, we thus have that $\lim_{t\to t_c-} s_k(t)/s_{k,\ompar}(t) = 1$.
% Therefore it suffices to prove
% $$ \alpha = \lim_{t \to t_c-} (t_c-t)s_{2,\ompar}(t) = \frac{1}{b(t_c)} \mbox{ and }\; \beta = \lim_{t \to t_c-} \frac{s_{3,\ompar}(t)}{(s_{2,\ompar}(t))^3} \in (0, \infty).$$
Write $y_\ompar(t)=1/s_{2,\ompar}(t)$ and $z_\ompar(t)=y^3_\ompar(t) s_{3,\ompar}(t)$, it suffices to show that:
\begin{equation}\label{eq:eq742}
	\lim_{t \to t_c-} \frac{t_c-t}{y_\ompar(t)} = \lim_{t \to t_c-} -\frac{1}{y'_\ompar(t)} = \frac{1}{b(t_c)}, \mbox{ and }  \lim_{t \to t_c-} z_\ompar(t) \in (0, \infty). \end{equation}
Define $A_l(t)=\sum_{i=1}^K(K+i)^l a_i(t)$ and $C_l(t)=\sum_{i=1}^{K}i^l c_i(t)$ for $l=1,2,3$. Then by Lemma \ref{lemma:s2omega-s3omegaA},  \eqref{eqn:fs2omega} and \eqref{eqn:fs3omega}, the derivative of $y_\ompar(t)$ and $z_\ompar(t)$ can be written as follows (we omit $t$ from the notation):
\begin{align}
	y_\ompar' =& -(A_2+C_2 x_\ompar) y_\ompar^2 - 2 C_1 y_\ompar - b, \label{eqn:diff-y}\\
	z_\ompar' =& y_\ompar^3 \left[ A_3+3C_1s_{3,\ompar} +3C_2s_{2,\ompar}+C_3x_\ompar + 3bs_{2,\ompar}s_{3,\ompar}  \right]\nonumber \\
	-& 3y_\ompar^2s_{3,\ompar} \left[(A_2+C_2 x_\ompar) y_\ompar^2 + 2 C_1 y_\ompar + b\right]\nonumber\\
	=&  -(3y_\ompar A_2+3y_\ompar C_2x_\ompar +3C_1)  z_\ompar +  (y_\ompar^3A_3+ 3y_\ompar^2C_2 +y_\ompar^3C_3 x_\ompar)  \nonumber\\
	=& - B_1 z_\ompar  + B_2,\label{eqn:diff-z}
\end{align}

where $B_1(t) = (3y_\ompar(t) A_2(t)+3y_\ompar(t) C_2(t)x_\ompar(t) +3C_1(t))$ and $B_2(t) = (y_\ompar^3(t)A_3(t)+ 3y_\ompar^2(t)C_2(t) +y_\ompar^3(t)C_3(t) x_\ompar(t)) $.
Since $\lim_{t \to t_c-} y_\ompar(t)=0$,  we have $\lim_{t \to t_c-}y'_\ompar(t) = - b(t_c)$ which proves the first statement in \eqref{eq:eq742}. 

Choose $t_1 \in (0,t_c)$ such that $y_{\ompar}(t), z_{\ompar}(t) \in (0,\infty)$ for all $t \in (t_1, t_c)$.
Then from \eqref{eqn:diff-z}, for all such $t$
$$ z_\ompar(t) =  \int_{t_1}^{t} e^{ -\int_s^{t} B_1(u)du } B_2(s)ds + z_\ompar(t_1)e^{ -\int_{t_1}^{t} B_1(u)du } .$$
Since $B_1, B_2$ are nonnegative and $\sup_{t \in [t_1, t_c]} \{ B_1(t) + B_2(t)\} < \infty$, we have $\lim_{t \to t_c-} z_\ompar(t) \in (0, \infty)$.  This completes the proof of \eqref{eq:eq742}.  
The result follows.
\qed\\

We now complete the proof of Theorem \ref{thm:suscept-funct}.\\

{\bf Proof of Theorem \ref{thm:suscept-funct}:} Let $\alpha, \beta$ be as introduced in Lemma \ref{lemma:alpha-beta}. 
From Lemma \ref{lemma:diff-s2s3} it follows that $y(t) = 1/ s_2(t)$ and $z(t) = y^3(t) s_3(t)$, for $0 \le t < t_c$, solve
the differential equations
\begin{align}
	\label{eq:eq1746}
	y'(t) = F^y(\bfx(t),y(t)),\;  z'(t) = F^z(\bfx(t), y(t), z(t)), \; y(0) =z(0) = 1,
\end{align}
where
 $F^y: [0,1]^{K+2} \to \RRR$  and $F^z: [0,1]^{K+2} \times \RRR \to \RRR$ are defined as
\begin{equation}
F^y(\bfx, y) := -y^2 F^s_2(\bfx , 1/y), \; F^z(\bfx, y, z) := 3 z F^y(\bfx, y)/y + y^3 F^s_3(\bfx, 1/y, z/y^3), \label{eqn:def-fy}
\end{equation}
$(\bfx, y, z) \in [0,1]^{K+2} \times \RRR \to \RRR$.
It is easy to check that $F^y$ and $F^z$ are polynomials in $(x_1, ...,x_K, x_\ompar, y)$ and $(x_1, ...,x_K, x_\ompar, y, z)$ respectively,
furthermore for each fixed $(\bfx, y) \in [0,1]^{K+2}$ 
the map $z \mapsto F^x(\bfx, y, z)$ is linear. 
%In particular  the function $g(t,z):=F^z(\bfx(t),y(t),z)$ is  Lipschitz  in $z$ a fact that will be used in the proof of Theorem \ref{theo:conv-pure-jump}.
Thus \eqref{eq:eq1746} has a unique solution.  Also, defining
 $y(t_c) = \lim_{t\to t_c-}y(t) = \lim_{t\to t_c-}y_{\ompar}(t)$ and $z(t_c) = \lim_{t\to t_c-}z(t) = \lim_{t\to t_c-}z_{\ompar}(t)$, we see that 
$y, z$ are twice continuously differentiable (from the left) at $t_c$.  Furthermore, $y'(t_c-) = -\alpha^{-1}$ and
$z(t_c-) = \beta$.  Thus we have 
$$ y(t) = \frac{1}{\alpha} (t_c-t)(1 + O(t_c-t)), \;\; z(t) = \beta (1 + O(t_c-t)), \mbox{ as } t \uparrow t_c. $$
The result follows. \qed\\

%\todo[inline]{(Xuan)  Differential equation theory needed. If the above assertion is not true, we can also just state the weaker result as in Lemma \ref{lemma:alpha-beta}.}

\subsection{Asymptotic analysis of $\bar s_2(t)$ and $\bar s_3(t)$}
\label{sec:decompose-s2-s3}
In preparation for the proof of Theorem \ref{thm:suscept-limit}, in this section we will obtain some useful semimartingale decompositions for
$ Y(t):=\frac{1}{\bar s_2(t)}$ and $Z(t):=\frac{\bar s_3(t)}{(\bar s_2(t))^3}$.  Throughout this section and next we will denote $ |\CC_1^{\sss (n)}(t)|$ as $I(t)$.
% 
% 
% $\calS_k(t) = \sum_{\CC \subset \BF_t} \CC^k$, $\bar s_k(t) = \calS_k(t)/n$ for $k=1,2,...$ and denote $\bfx=(x_i)_{i \in \Omega_K}$, $\bar \bfx(t) = (\bar x_i(t))_{i \in \Omega_K}$ and $\bfx(t) = (x_i(t))_{i \in \Omega_K}$.\\
% 
% 
% In this section and the next, we denote $I(t) := \CC_1^{\sss (n)}(t)$ and 
% $$ Y(t):=\frac{1}{\bar s_2(t)}, \;\;\; Z(t):=\frac{\bar s_3(t)}{(\bar s_2(t))^3},\;\;\; y(t):= \frac{1}{s_2(t)}, \;\;\; z(t) := \frac{ s_3(t)}{( s_2(t))^3},$$
% where $s_2(t)$ and $s_3(t)$ are as in Lemma \ref{lemma:diff-s2s3}. 
Recall the functions $F^s_2$, $F^s_3$ introduced in Lemma \ref{lemma:diff-s2s3}. 
%The goal of this section is to decompose $Y(t)$ and $Z(t)$ in the way that the conditions in Theorem \ref{thm:suscept-limit} are satisfied. The main decomposition results are stated in Lemma \ref{lemma:decomp-s2-y} and Lemma \ref{lemma:decomp-s3-z}.\\

%First we summarized the analysis of $\bar s_2(t)$ and $\bar s_3(t)$ in the following lemma.
\begin{Lemma} 
	\label{lemma:decompose-s2-jumps}
The processes $\bar s_2$ and $\bar s_3$ are $\{\clf_t\}_{0\le t < t_c}$ semimartingales of the form \eqref{eq:semimart} and the  following equations hold.\\
\begin{enumerate}[(a)]
\item	$\bfd(\bar s_2)(t) = F^s_2( \bar \bfx(t), \bar s_2(t)) + O_{t_c}\left(I^2(t)\bar s_2(t)/{n}\right).$  \\
\item	$\bfd(\bar s_3)(t) = F^s_3(\bar \bfx(t), \bar s_2(t) ,\bar s_3(t)) + O_{t_c}\left( {I^3(t) \bar s_2(t)}/{n}\right).$\\
\item	$\bfv(\bar s_2)(t)  = O_{t_c}(I^2(t)\bar s_2^2(t)/n).$\\
% \item    	$\langle \bfm(\bar s_2), \bfm(\bar s_3)\rangle_s = O_{t_c}\left({I^3(t)\bar s_2^2(t)}/{n}\right).$\\
% \item	$\langle \bfm(\bar s_3), \bfm(\bar s_3)\rangle_t   = O_{t_c}(I^4(t)\bar s_2^2(t)/n).$
\end{enumerate}
\end{Lemma}
\textbf{Proof: } Parts (a) and (b) are immediate from \eqref{eqn:s2omega-to-s2}, Lemma \ref{lemma:approx-xi}(b) and Lemma \ref{lemma:s2omega-s3omega}.
For part (c), recall the three types of events described in Section \ref{sec:prelim-bsr}.
For type 1, $\Delta \bar s_2(t)$ is bounded by $2K^2/n$ and the total rate of such events is bounded by $n/2$.  For type 2, the attachment
of a size $j$ component, $1\le j \le K$, to a component $\CC$ in $\BS^*_{t-}$ occurs at rate $|\CC| c^*_j(t)$ and produces a jump
$\Delta \bar s_2(t) = 2j|\CC|/n$.  For type 3, components $\CC$ and $\tilde \CC$ merge at rate $|\CC| |\tilde \CC| b^*(t)/n$ and produce a jump  $ \Delta \bar s_2(t) = 2 |\CC| |\tilde \CC|/n$.  %In addition, a jump in $\bar s_2$ can be produced by a merger of two components of size at most $K$ into a new component of sizebounded by $K$.  For such events $\Delta \bar s_2(t)$ is bounded by $2K^2/n$ and the total rate of such events is bounded by $n/2$.
Thus  for $t \in [0, t_c)$, $\bfv(\bar s_2)(t)$ can be estimated as
\begin{align*}
	\bfv(\bar s_2)(t) \le &
	   \frac{n}{2} \left( \frac{2K^2 }{n}\right)^2  
	 + \sum_{j=1}^K \sum_{\CC \subset \BS^*_{t}} \left(\frac{2j|\CC|}{n}\right)^2 |\CC| c^*_j(t)
	 + \sum_{\CC \neq \tilde \CC \subset \BS^*_t} \left(\frac{2|\CC| |\tilde \CC|}{n}\right)^2    \frac{b^*(t) |\CC| |\tilde \CC|}{n} \\ 
	\le& \frac{2K^4}{n} + \frac{4K^2 \calS_3 }{n^2}+ \frac{4(\calS_3)^2}{n^3} = O_{t_c}\left(\frac{I^2(t)\bar s_2^2(t)}{n}\right).
\end{align*}
This proves (c). 
%Parts (d) and (e) are immediate from (c) on observing that $\Delta \bar s_3(t) \le 3 I(t) \Delta \bar s_2(t)$. The result follows. 
\qed\\

In the next lemma, we obtain a semimartingale decomposition for $Y$.  

%with maximum total order 6, and the order of $y$ is 2. Since $|x_i(t)| \le 1$ and $y(t) \le 1$ for all $t$, thus $g(t,y)=F^y(\bfx(t),y)$ is Lipchitz on $y \in [0,1]$. 
\begin{Lemma}
	\label{lemma:decomp-s2-y}
The process $Y(t)=1/\bar s_2(t)$ is a $\{\clf_t\}_{0\le t < t_c}$ semimartingale of the form \eqref{eq:semimart} and

(i) With $F^y(\cdot)$  as defined in \eqref{eqn:def-fy}, 
\begin{equation}
	\bfd(Y)(t) = F^y(\bar \bfx(t), Y(t)) + O_{t_c}\left(\frac {I^2(t)Y(t)}{n}\right). \label{eqn:delta-y}
\end{equation}
(ii)  
\begin{equation*}
\bfv(Y)(t) = O_{t_c}\left(\frac{ I^2(t)Y^2(t)}{n} \right).
%+ \frac{ I^4(t)Y^4(t)}{n}\right).
\end{equation*}
\end{Lemma}

{\bf Proof:} Note that 
\be
\Delta Y(t) = \frac{1}{\bar s_2+\Delta \bar s_2}- \frac{1}{\bar s_2}= -\frac{\Delta \bar s_2}{\bar s_2^2}+ \frac{(\Delta \bar s_2)^2}{\bar s_2^2(\bar s_2 + \Delta \bar s_2)}  =-\frac{\Delta \bar s_2}{\bar s_2^2}+ O_{t_c}\left(\frac{(\Delta \bar s_2)^2}{\bar s_2^3}\right). \label{eqn:2591}
\ee
Thus by Lemma \ref{lemma:decompose-s2-jumps}(a), we have, 
\begin{align*}
	\bfd(Y)(t)
	=& -\frac{1} {(\bar s_2(t))^2} \bfd(\bar s_2)(t) + O_{t_c}\left( \frac{1}{(\bar s_2(t))^3} \bfv(\bar s_2)(t)\right)\\
	=& \left( -\frac{1} {(\bar s_2(t))^2}  \right) \left( F^s_2( \bar \bfx(t), \bar s_2(t)) +O_{t_c}\left(\frac{I^2(t)\bar s_2(t)}{n}\right) \right) + O_{t_c}\left( \frac{1}{(\bar s_2(s))^3} \cdot \frac{I^2(t)\bar s_2^2(t)}{n} \right) \\
	=& F^y(\bar \bfx(t), Y(t)) + O_{t_c}\left(\frac{I^2(t)Y(t)}{n}\right).
\end{align*}
This proves (i). For (ii), note that \eqref{eqn:2591} also implies
$ (\Delta Y(t))^2 \le \frac{(\Delta \bar s_2)^2}{\bar s_2^4}$
We then have
$$ \bfd(Y)(t)  \le  \frac{2\bfv(\bar s_2)(t)}{\bar s_2^4}  = O_{t_c}\left(\frac{I^2(t)Y^2(t)}{n}\right).$$
The result follows. \qed\\

We now give a semimartingale decomposition for $Z(t)=\bar s_3(t)/(\bar s_2(t))^3$. 
% Notice that 
% \begin{align*}
% 	z'(t) = 3y^2(t)s_3(t) y'(t) + y^3(t) s'_3(t) = 3z(t)\cdot \frac{1}{y(t)} F^y(\bfx(t), y(t)) + y^3(t) F^s_3(\bfx(t), 1/y(t), z(t)/y^3(t)), 
% \end{align*}

\begin{Lemma}
	\label{lemma:decomp-s3-z}
 The process $Z(t)=\bar s_3(t)/(\bar s_2(t))^3$ is a $\{\clf_t\}_{0\le t < t_c}$ semimartingale of the form \eqref{eq:semimart} and

(i) With $F^z(\cdot)$  as defined in \eqref{eqn:def-fy},
\begin{equation*}
	\bfd(Z)(t) =  F^z (\bar \bfx(t), Y(t), Z(t)) +  O_{t_c}\left(\frac{I^3(t)Y^2(t)}{n}\right).
\end{equation*}
(ii) 
\begin{equation*}
	 \bfv(Z)(t) =  O_{t_c}\left(\frac{ I^4(t)Y^4(t)}{n} + \frac{ I^6(t)Y^6(t)}{n}\right).
\end{equation*}
\end{Lemma}
\textbf{Proof: } Note that
\begin{equation*}
	\Delta Z = Y^3 \Delta \bar s_3 + 3 Y^2 \bar s_3 \Delta Y + R(\Delta Y, \Delta \bar s_3),
\end{equation*}
where $R(\Delta Y, \Delta \bar s_3)$ is the error term which,  using the observations
that $ \bar s_3 \le I \bar s_2$, $\Delta \bar s_3 \le 3 I \Delta \bar s_2$ and $|\Delta Y| \le Y^2 \Delta \bar s_2$,  can be bounded as follows. 
\begin{align*}
	|R(\Delta Y, \Delta \bar s_3)|
	\le& 3 Y^2 |\Delta Y||\Delta \bar s_3| + 3 Y \bar s_3 |\Delta Y|^2 \\
	\le& 3Y^2 \cdot {Y^2 \Delta \bar s_2} \cdot {3I\Delta \bar s_2} + 3 I \cdot \left( {Y^2 \Delta \bar s_2 }\right)^2
	= 12I Y^4 \cdot (\Delta \bar s_2)^2.
\end{align*}
From Lemma \ref{lemma:decompose-s2-jumps}(b), Lemma \ref{lemma:decomp-s2-y}(i) and Lemma \ref{lemma:decompose-s2-jumps}(c), we have
\begin{align*}
\bfd(Z)(t)
	=& Y^3(t) \bfd(\bar s_3)(t) + 3Y^2(t) \bar s_3(t) \bfd(Y)(t) + O_{t_c}\left(I(t) Y^4(t)\bfv(\bar s_2)(t)\right)\\
	=& Y^3(t)  \left(F^s_3(\bar \bfx(t),\bar s_2(t), \bar s_3(t))+ O_{t_c}\left(\frac{I^3(t) \bar s_2(t)}{n}\right)\right) \\
	&+3Y^2(t) \bar s_3(t) \left( F^y(\bar \bfx(t), Y(t)) + O_{t_c}\left(\frac{I^2(t) Y(t)}{n}\right)\right) + O_{t_c}\left( \frac{I^3(t)Y^2(t)}{n}\right)\\
	=& F^z(\bar \bfx(t), Y(t),Z(t)) + O_{t_c}\left(\frac{I^3(t) Y^2(t)}{n}\right).
\end{align*}
This proves (i). For  (ii), note that 
\bes
Y^3 |\Delta \bar s_3| + 3 Y^2\bar s_3 |\Delta Y| \le Y^3 \cdot 3I|\Delta \bar s_2|+ 3Y^2 \cdot I \bar s_2 \cdot Y^2|\Delta \bar s_2|=6Y^3I|\Delta \bar s_2|.
\ees
Thus, 
$$
|\Delta Z| \le 6Y^3I|\Delta \bar s_2| + 12 I Y^4 \cdot (\Delta \bar s_2)^2.$$
Applying Lemma \ref{lemma:decompose-s2-jumps}(c)  we now have,
\bes
 \bfv(Z)(t) = O_{t_c} \left( Y^6 I^2 \bfv(\bar s_2)(t)\right) + O_{t_c} \left(\frac{I^6 Y^6}{n}\right) =
 O_{t_c}\left( \frac{I^4Y^4}{n} + \frac{I^6 Y^6}{n} \right).
\ees
The result follows. \qed \\

\subsection{Proof of Theorem \ref{thm:suscept-limit}}
\label{sec:proof-conv-susceptibility}
We begin with  an upper bound on the size of the largest component at time $t \le t_n = t_c - n^{-\gamma}$ for $\gamma \in (0,1/4)$, which has been proved in \cite{bsr-2012}, and will play an important role in the proof of Theorem \ref{thm:suscept-limit}.
\begin{Theorem}[{\cite[Theorem 1.2]{bsr-2012}} {\bf Barely subcritical regime}]
\label{thm:subcrit-reg}
Fix  $\gamma \in (0,1/4)$.  Then there exists  $C_3  \in (0, \infty)$ such that, as $n \to \infty$,
\[ \prob\set{ I^{\sss (n)}(t) \le C_3 \frac{(\log n)^4}{(t_c-t)^2},~ \forall t <t_c-n^{-\gamma} } \to 1. \]
\end{Theorem}
The next lemma is an elementary consequence of Gronwall's inequality.
% In this section we prove the following abstract convergence theorem which captures the main ingredients of the proof of Theorem \ref{thm:suscept-limit}. For any sequence of time $\{t_n\}_n$, we shall use $\{Y^{\sss (n)}(t) : t \in [0,t_n] \}_{n} $ to denote a sequence of stochastic jump processes which are adapted to the filtrations $\{\FF^{\sss (n)}_t : t \in [ 0,t_n]\}_n$. We further assume 
% $$ \sup_n t_n \le T < \infty \mbox{ and }  \expec[(Y^{\sss (n)}(t))^2]<\infty. $$
% For a fixed $T$, let $g(t,y):[0,T]\times \Rbold\to \Rbold$ be a function satisfying a uniform Lipschitz condition: there is constant $C_1$, independent of $t$, such that for all $ t \in [0,T]$ and $y_1, y_2 \in \RRR$, we have
% $$ |g(t,y_1)-g(t,y_2)| \le C_1|y_1-y_2|. $$
% Let $\{y(t)\}_{t \in [0,T]}$ be the unique solution of the differential equation
% $$ y^\prime(t) = g(t,y(t)), \;\;\; y(0)=y_0. $$
% The following theorem establish the rate of convergence from $Y^{\sss (n)}(\cdot)$ to $y(\cdot)$ under certain conditions.
\begin{Lemma}
	\label{theo:conv-pure-jump}
Let $\{t_n\}$ be a sequence of positive reals such that $t_n \in [0, t_c)$ for all $n$.  Suppose that $U^{\sss(n)}$ is a semimartingale of the form \eqref{eq:semimart} with values in $\DDD \subset \RRR$.
Let $g:[0,t_c)\times \DDD \to \Rbold$ be such that, for some $C_4(g) \in (0, \infty)$,  
\begin{equation} \sup_{t \in [0, t_c)}|g(t,u_1)-g(t,u_2)| \le C_4(g)|u_1-u_2|, \; u_1, u_2 \in \DDD. \label{eq:eq2127}\end{equation}
Let $\{u(t)\}_{t \in [0,T]}$ be the unique solution of the differential equation
$$ u^\prime(t) = g(t,u(t)), \;\;\; u(0)=u_0. $$
Further suppose that there exist positive sequences:
\begin{enumerate}[(i)]
\item   $\{\theta_1(n)\}$  such that, whp, 
$|U^{\sss(n)}(0)-u_0|  \le \theta_1(n)$.

\item     $\{\theta_2(n)\}$  such that, whp,
$$  \int_0^{t_n}\left| \bfd(U^{\sss(n)})(t) - g(t,U^{\sss(n)}(t))\right|dt \le \theta_2(n). $$

\item    $\{\theta_3(n)\}$  such that, whp, 
$\langle\bfm(U^{\sss(n)}), \bfm(U^{\sss(n)})\rangle_{t_n} \le \theta_3(n)$.
\end{enumerate}
Then, whp, 
$$ \sup_{0\le t \le t_n}|U^{\sss (n)}(t)-u(t)| \le e^{C_4(g)T}(\theta_1(n) + \theta_2(n) + \theta_4(n)), $$
where $\theta_4=\theta_4(n)$ is any sequence satisfying $ \sqrt{\theta_3(n)} = o(\theta_4(n))$.
\end{Lemma}

{\bf Proof:} We suppress $n$ from the notation unless needed.  Using the Lipschitz property of $g$, we have, for all $t \in [0, t_n]$,
\begin{align*}
	|U(t)-u(t)|
	\le& |U(0)-u_0| + \int_0^t |\bfd(U)(s)-g(s,U(s))| ds + \int_0^t|g(s,U(s))-g(s,u(s))|ds + |\bfm(U)(t)|\\
	\le& |U(0)-u_0| + \int_0^t |\bfd(U)(s)-g(s,U(s))| ds + |\bfm(U)(t)| + C_4\int_0^t |U(s)-u(s)|ds.
\end{align*}
 Then by Gronwall's lemma 
\begin{equation}\label{eq:eq1844}
	\sup_{0 \le t \le t_n}| U(t)-u(t)| \le \left( |U(0)-u_0| + \int_0^{t_n} |\bfd(U)(s)-g(s,U(s))| ds + \sup_{0 \le t \le t_n}|\bfm(U)(t)| \right) e^{C_4 T}.
\end{equation}
Let  $\tau^{\sss (n)} = \inf\{ t\ge 0: \langle \bfm(U), \bfm(U) \rangle_t > \theta_3(n) \}$.  By Doob's inequality
$$ \E[\sup_{0\le t\le t_n}|\bfm(U)(t \wedge \tau)|^2] \le 4 \E[ |\bfm(U)(t_n \wedge \tau)|^2] = 4 \E \left[ \langle \bfm(U), \bfm(U) \rangle_{t_n \wedge \tau}\right]  \le 4 \theta_3(n). $$
Then for any $\theta_4(n)$ such that $ \theta_3=o((\theta_4)^2)$, we have
\begin{align*}
	\prob\{ \sup_{0 \le t \le t_n}|\bfm(U)(t)| > \theta_4(n) \} 
	\le& \prob\{ \tau^{\sss (n)} < t_n \}+ \prob\set{ \sup_{0\le t\le t_n}|\bfm(U)(t \wedge \tau)| > \theta_4(n) }\\
	\le& \prob\{ \langle \bfm(U), \bfm(U) \rangle_{t_n} > \theta_3(n) \} + 4 \theta_3(n)/\theta_4^2(n) \to 0.
\end{align*}
The result now follows on using the above observation in \eqref{eq:eq1844}.
 \qed \\

% In next section we will analyze the dynamic of $1/\bar s_2(t)$ and $\bar s_3(t)/\bar s_2^3(t)$, and give the explicit bounds on $\lambda_1, \lambda_2, \lambda_3$ that appear in the conditions of the above theorem.\\
% 
% 
% In this section we wrap up results in previous sections and prove Theorem \ref{thm:suscept-limit}.\\

{\bf Proof of Theorem \ref{thm:suscept-limit}:} Let $y$ and $z$ be as in the proof of Theorem  \ref{thm:suscept-funct}.  It suffices to show
\begin{align}
\sup_{0 \le t \le t_n} \left| Y(t) - y(t)\right|  n^{1/3} \convp 0 \label{eqn:conv-s2y} \\
\sup_{0 \leq t\le t_n} |Z(t)- z(t)| \convp 0. \label{eqn:conv-s3z}
\end{align} 
We begin by proving the following weaker result than \eqref{eqn:conv-s2y}:
\begin{equation}
	\sup_{0 \le t \le t_n}|Y(t)-y(t)| = O(n^{-1/5}), \; \mbox{whp}. \label{eqn:1142}
\end{equation}
Recalling from Theorem \ref{thm:suscept-funct} that $\bfx \mapsto F^y(\bfx, y)$ is Lipschitz, uniformly in $y$, we get for some $d_1 \in (0, \infty)$
\begin{equation*}
	\sup_{0 \le t \le T}|F^y(\bar \bfx(t), Y(t))-F^y(\bfx(t), Y(t))| \le d_1 \sup_{i \in \Omega_K} \sup_{ 0 \le t \le T}|\barx_i(t)-x_i(t)|.
\end{equation*}
From Lemma \ref{lemma:decomp-s2-y}(ii) and Lemma \ref{lemma:approx-xi}(a) we now get for some $d_2 \in (0, \infty)$, whp,
$$|\bfd(Y)(t) - F^y(\bfx(t), Y(t))|  \le  d_2\left(\frac {I^2(t)Y(t)}{n}  + n^{-2/5}\right), \mbox{ for all } t \in [0, t_n].$$
Thus, from Theorem \ref{thm:subcrit-reg} and recalling that $\gamma < 1/5$, we have whp,
\begin{align*}
	\int_0^{t_n} |\bfd(Y)(t) - F^y(\bfx(t), Y(t))|dt  =& O \left(\int_0^{t_n} \frac{(\log n)^8}{n(t_c-t)^4}dt + n^{-2/5}\right)  \\
	=& O((\log n)^8n^{3\gamma-1})+O(n^{-2/5})=O(n^{-2/5}).
\end{align*}

Next, by Lemma \ref{lemma:decomp-s2-y}(ii) and using the fact $Y(t) \le 1$ for all $t \in [0, t_c)$,
\begin{align}
	\langle \bfm(Y), \bfm(Y) \rangle_{t_n}
	 =& O\left( \int_0^{t_n} \frac{I^2(t) Y^2(t)}{n} dt\right)
	 = O\left( \int_0^{t_n} \frac{I^2(t)}{n} dt\right)\nonumber\\
	 =&
	O\left( \int_0^{t_n} \frac{(\log n)^8}{n(t_c-t)^4}dt \right)
	 = O((\log n)^8n^{3\gamma-1}). \label{eqn:2961}
\end{align}
The statement in \eqref{eqn:1142} now follows on observing that $((\log n)^8n^{3\gamma-1})^{1/2} = o(n^{-1/5})$ and applying Lemma \ref{theo:conv-pure-jump} with $\DDD:=[0,1]$, $g(t,y):=F^y(\bfx(t),y)$, $\theta_1=0$, $\theta_2 = n^{-2/5}$ and $\theta_3 = (\log n)^8n^{3\gamma-1}$.\\

We now strengthen the estimate in \eqref{eqn:1142} to prove \eqref{eqn:conv-s2y}.  From Theorem \ref{thm:suscept-funct} it follows that
$y(t_n) = \Theta(n^{-\gamma})$.  Since $\gamma < 1/5$, from \eqref{eqn:1142} we have, whp, $Y(t) \le 2 y(t)$ for all $t \le t_n$. 
Thus from the first equality in \eqref{eqn:2961} and Theorem \ref{thm:suscept-funct} we get, whp,
\bes
\langle \bfm(Y), \bfm(Y) \rangle_{t_n} = O\left( \int_0^{t_n} \frac{I^2(t) y^2(t)}{n} \right) = O\left( \int_0^{t_n} \frac{(\log n)^8}{n(t_c-t)^2}dt \right) = O((\log n)^8n^{\gamma-1}).
\ees
Since $((\log n)^8 n^{\gamma-1})^{1/2}=o(n^{-2/5})$, applying Lemma \ref{theo:conv-pure-jump} again gives
\be
	\sup_{0 \le t \le t_n}|Y(t)-y(t)| = O(n^{-2/5}), \; \mbox{whp}. \label{eqn:2971}
\ee
This proves \eqref{eqn:conv-s2y}.\\

We now prove \eqref{eqn:conv-s3z}.  We will apply Lemma \ref{theo:conv-pure-jump} to $\DDD := \RRR$ and
$
	g(t,z):= F^z(\bfx(t),y(t),z).
$
As noted in the proof of Theorem \ref{thm:suscept-funct}, $g$ defined as above satisfies \eqref{eq:eq2127}.

We now verify the three assumptions in Lemma \ref{theo:conv-pure-jump}.
Note that (i) is satisfied with $\theta_1 = 0$, since $Z(0)=z(0)=1$. Next,  by Lemma \ref{lemma:decomp-s3-z}(ii) and the fact $Y(t) \le 2y(t)$ for $t \le t_n$, whp, we have
\begin{align*}
\langle \bfm(Z), \bfm(Z) \rangle_{t_n} =& O \left(\int_0^{t_n} \left(\frac{I^4(t)Y^4(t)}{n} + \frac{I^6(t)Y^6(t)}{n}\right) dt \right)\\
 =& O\left(\int_0^{t_n}\left( \frac{ (\log n)^{16}}{n (t_c-t)^4} + \frac{ (\log n)^{24}}{n (t_c-t)^6}\right)  dt\right)
 = O((\log n)^{24}n^{5\gamma-1}). 
\end{align*}
Since $\gamma < 1/5$, we can find $\theta_4(n) \to 0$ such that $\sqrt{(\log n)^{24}n^{5\gamma-1} }= o(\theta_4(n))$
 Thus (iii) in  Lemma \ref{theo:conv-pure-jump} is satisfied. Next recall from the proof of Theorem \ref{thm:suscept-funct} that $g(t,z)$ is linear in $z$.  Also,  $Z(t) \le I(t)$.  Thus from Lemma \ref{lemma:approx-xi} and
\eqref{eqn:2971}, for some $d_3 \in (0, \infty)$ whp, for all $t \le t_n$
\begin{align*}
	&|F^z(\bar \bfx(t), Y(t), Z(t)) -g(t,Z(t))| \\
	\le&  d_3 (1 + Z(t)) \left(\sup_{1 \le i \le K} \sup_{0 \le t \le t_n}|\barx_i(t)-x_i(t)|
	 +  \sup_{0 \le t \le t_n}|Y(t)-y(t)|\right)
	=  I(t)O(n^{-2/5}).
\end{align*}
By  Lemma \ref{lemma:decomp-s3-z}(i) and the above bound,
\begin{align}
	\int_0^{t_n}|\bfd(Z)(t) - g(t,Z(t))|dt  
	=& O\left( \int_0^{t_n} n^{-2/5} I(t)dt \right)+O\left( \int_0^{t_n} \frac{y^2(t)I^3(t)}{n}  dt \right) \nonumber\\
	=& O((\log n)^4 n^{\gamma -2/5}) + O((\log n)^{12} n^{3\gamma -1}). \label{eqn:1297}
\end{align}
This verifies (ii) in  Lemma \ref{theo:conv-pure-jump} with $\theta_2(n) = O((\log n)^{12} n^{3\gamma -1})$.
From Lemma \ref{theo:conv-pure-jump} we now have
\bes
	\sup_{0\le t\le t_n}|Z(t)-z(t)|  \le \theta_1(n) + \theta_2(n) + \theta_4(n) =  o(1).
\ees
The result follows. \qed\\

\section{Coupling with the multiplicative coalescent}
\label{sec:main-coupling}

In this section we prove Theorem \ref{thm:crit-regime}.  Throughout this section we fix $\gamma \in (1/6,1/5)$.  The basic idea of the proof is as follows.  Recall $\alpha, \beta \in (0, \infty)$ from Theorem \ref{thm:suscept-funct} (see also Lemma \ref{lemma:alpha-beta}).  
We begin by approximating the BSR random graph process by a  process which until time $t_n := t_c - n^{-\gamma}$ is identical to
the BSR process and in the time interval $[t_n, t_c + \alpha \beta^{2/3}\frac{\lambda}{n^{1/3}}]$ evolves as an \erdos process, namely over this interval
edges between any pair of vertices appear at rate $1/\alpha n$, and self loops at any given vertex appear at rate  $1/2\alpha n$.
Asymptotic behavior of this random graph is analyzed using Theorem \ref{theo:aldous-full-gene}.  Theorems \ref{thm:suscept-funct}, \ref{thm:subcrit-reg}
and \ref{thm:suscept-limit} help in verifying the conditions \eqref{eqn:qsigma2}	and \eqref{eqn:additional-condition} in the statement of
Theorem \ref{theo:aldous-full-gene}.  We then complete the proof of  Theorem \ref{thm:crit-regime} by arguing  that the `difference' between the BSR process and the modified random graph process is asymptotically negligible.  

Let
\begin{equation}
\label{eqn:lambdan-defn}
t_n=t_c - \alpha \beta^{2/3}\frac{\lambda_n}{n^{1/3}} \mbox{ where }	\lambda_n = \frac{ n^{1/3-\gamma}}{\alpha \beta^{2/3}}.  
\end{equation}
Throughout this section, for $\lambda \in \RRR$, we denote $t^{\lambda} = t_c + \alpha \beta^{2/3}{\lambda}/{n^{1/3}}$.
Recall the random graph process $\BS^*(t)$ introduced in Section \ref{sec:main-bsr-susceptibility}.  Denote by
$(|\CC_{i}^{*}(t)|, \xi_{i}^{*}(t))_{i\ge 1}$ the vector of ordered component size and corresponding surplus in $\BS^*(t)$ (the components
are denoted by $\CC_{i}^{*}(t)$ ).
Let, for $\lambda \in \RRR$, 
$$ \bar \bfC^{{\sss (n)},*}(\lambda) =  \left(\frac{\beta^{1/3}}{n^{2/3}}\left|\CC_{i}^{*}(t^{\lambda})\right| :i\geq 1\right), \;\; \bar \bfY^{{\sss (n)},*}(\lambda) =  \left(\xi_{i}^{*}(t^{\lambda}) :i\geq 1\right).$$

For $i \ge 1$, denote $\bar \bfC_i^{{\sss (n)},*}(\lambda)$ and $\bar \bfY_i^{{\sss (n)},*}(\lambda)$ for the $i$-th coordinate of $\bar \bfC^{{\sss (n)},*}(\lambda)$ and $\bar \bfY^{{\sss (n)},*}(\lambda)$ respectively. 
Write
$\bar \bfY_i^{{\sss (n)},*} = \tilde \xi_{i}^{\sss (n)} + \hat \xi_{i}^{\sss (n)}$ where
$\tilde \xi_{i}^{\sss (n)}(\lambda)$ represents the surplus in $\BS^*(t^{\lambda})$ that is created before time $t_n$, namely
$$\tilde \xi_{i}^{\sss (n)}(\lambda) = \sum_{j: \CC_j^{*}(t_n) \subset \CC_{i}^{*}(t^{\lambda})} \bar \bfY_j^{{\sss (n)},*}(-\lambda_n).$$
In Section \ref{sec:old-surplus} we will show that the contribution from $\tilde \xi^{\sss (n)}(\lambda):= (\tilde \xi_i^{\sss (n)}(\lambda): i \ge 1)$ is asymptotically negligible.
First, in Section \ref{sec:two-coupling} below we analyze the contribution from the `new surplus', i.e. $\hat \xi^{\sss (n)}:= (\hat \xi_i^{\sss (n)} : i \ge 1)$.
\subsection{Surplus created after time $t_n$.}
\label{sec:two-coupling}
The main result of this section is as follows.  Recall the process $\bfZ(\lambda) = (\bfX(\lambda), \bfY(\lambda))$ introduced in Theorem \ref{thm:smc-surplus}.
\begin{Theorem}
	\label{theo:onedimc}
	For every $\lambda \in \RRR$, as $n\to \infty$, $(\bar \bfC^{{\sss (n)},*}(\lambda), \hat\xi^{\sss (n)}(\lambda))$ converges in distribution, in $\udown$, to
	$({\bfX}(\lambda), \bfY(\lambda))$.
\end{Theorem}
The basic idea in the proof of the above theorem is to argue that $\BS^*(t^{\lambda})$ `lies between' two \erdos random graph processes $\bfG^{\sss(n),-}(t^\lambda)$ and $\bfG^{\sss(n),+}(t^\lambda)$, whp, and then apply 
Theorem \ref{theo:aldous-full-gene} to each of these processes.  For a graph $\bfG$, denote by $|\CC_i(\bfG)|$ and $\xi_i(\bfG)$  the size and surplus, respectively, of the $i$-th largest component, $\CC_i(\bfG)$ of graph $\bfG$.
We begin with the following lemma. Recall $\lambda_n$ from \eqref{eqn:lambdan-defn}. 
% In this section we shall couple the process $\BS^*(\lambda):= \BS^*\left(t_c+\frac{\alpha\beta^{2/3}}{n^{1/3}}\lambda\right)$ with two 
% Erd\'{o}s-R\H{e}nyi process $\bfG^{\sss(n),-}(\lambda)$ and $\bfG^{\sss(n),+}(\lambda)$ such that $\BS^*(\lambda)$ is sandwiched between the two processes. Precisely, for any graph $\bfG$, write $\CC_i(\bfG)$ and $\xi_i(\bfG)$ for the size and surplus of the $i$-th largest component of graph $\bfG$, then define 
% \be
% \bar \bfC^{\sss(n),-}(\lambda) := \left(\frac{\beta^{1/3}}{n^{2/3}}\CC_i\left(\bfG^{\sss(n),-}(\lambda)\right) :i\geq 1 \right) 
% \mbox{ and }  \bar{\bfY}^{\sss (n),-}(\lambda) = \left (\xi_i\left( \bfG^{\sss(n),-}(\lambda)\right) :i\geq 1\right).
% \ee
% Also define $ \bar \bfC^{\sss(n),+}(\lambda)$ and $\bar{\bfY}^{\sss (n),+}(\lambda)$ for $\bfG^{\sss(n),+}(\lambda)$ similarly. Then we have the following lemma.
\begin{Lemma}
	\label{lemma:upper-lower-coupling}
	 There exists a construction of $\{\BS^*(t)\}_{t\ge 0}$ along with two other random graph processes $\{\bfG^{\sss(n),-}(t)\}_{t\ge 0}$ and $\{\bfG^{\sss(n),+}(t)\}_{t\ge 0}$
	such that:\\
%random graph processes $\bfG^{\sss(n),-}(\lambda)$ and $\bfG^{\sss(n),+}(\lambda)$ for $\lambda \in [-\lambda_n, \lambda_n]$ such that\\
	(i) With high probability,
	\begin{equation}
		\label{eq:eq1136}
		\bfG^{\sss(n),-}(t^\lambda) \subset \BS^*(t^\lambda)\subset \bfG^{\sss(n),+}(t^\lambda) \qquad \mbox{ for all } \lambda \in [-\lambda_n, \lambda_n]. \end{equation}
	(ii) Let for $i \ge 1$, $\bar \bfC^{\sss(n),\mp}_i(\lambda) = \frac{\beta^{1/3}}{n^{2/3}}|\CC_i\left(\bfG^{\sss(n),\mp}(t^\lambda)\right)|$
	and $$\bar{\bfY}^{\sss (n),\mp}_{i}(\lambda) = \xi_i\left( \bfG^{\sss(n),\mp}(t^\lambda)\right) - \sum_{j: \CC_j\left(\bfG^{\sss(n),\mp}(t_n)\right) \subset \CC_i\left(\bfG^{\sss(n),\mp}(t^{\lambda})\right)}
	 \xi_j\left( \bfG^{\sss(n),\mp}(t_n)\right).$$
	Then, for all $\lambda \in \RRR$
	\bes
	(\bar \bfC^{\sss(n),\bullet}(\lambda), \bar \bfY^{\sss(n),\bullet}(\lambda)) \convd (\bfX(\lambda), \bfY(\lambda)), \; \bullet = -, +,
	\ees
	where $\convd$ denotes weak convergence in $\udown$.
\end{Lemma}
We remark that $\bar{\bfY}^{\sss (n),\mp}(\lambda)$ represents the surplus in $\bfG^{\sss(n),\mp}(t^\lambda)$ created after time instant $t_n$.
Proof of the lemma relies on the following proposition which is an immediate consequence of Theorem \ref{thm:suscept-funct}, Theorem \ref{thm:suscept-limit} and Theorem \ref{thm:subcrit-reg}. 
\begin{Proposition}
	\label{prop:three-conditions} There exists a $\kappa \in (0, \frac{1}{3} - \gamma)$ such that
	$$ \frac{\bar s_3(t_n)}{(\bars_2(t_n))^3} \convp \beta, \;\; \frac{n^{1/3}}{\bars_2(t_n)}-\frac{n^{1/3-\gamma}}{\alpha} \convp 0, \;\; \frac{I^{\sss (n)}(t_n)}{n^{2\gamma + \kappa}} \convp 0. $$
\end{Proposition}
We now prove Lemma \ref{lemma:upper-lower-coupling}.

\textbf{Proof of Lemma \ref{lemma:upper-lower-coupling}:}\\
%**********************
We suppress $n$ in the notation for the random graph processes. Write $t_n^+:=t_c+n^{-\gamma}$.  Let
$\BS(t)$ for $t \in [0, t_n^+]$ be constructed as in Section \ref{sec:bsr} and define $\BS^*(t)$ for $t \in [0, t_n)$
as in Section \ref{sec:main-bsr-susceptibility}. Set
$$
\bfG^{\sss(n),-}(t) =  \bfG^{\sss(n),+}(t) = \BS^*(t), \mbox{ for } t \in [0, t_n).
$$
We now give the construction of these processes for
 $t \in [t_n, t_n^+]$. \\

The construction is done in two rounds. In the first round, we construct processes $\bfG^{I,-}(t)$, $\BS^{I,*}(t)$ and $\bfG^{I,+}(t)$ for $ t \in [t_n ,t_n^+]$ by using only the information of immigrations and attachments in $\BS(t)$, while the edge formation between large components is ignored.  We first construct the process $\{\BS^{ I}(t)\}_{t \in [t_n, t_n^+] }$ as follows. Let $\BS^{ I}(t_n) := \BS(t_n)$. For $t > t_n$, $\BS^{ I}(t)$ is constructed 
along with and same as $\BS(t)$, except for when
$$ c_{t-}(\vec{v}) \in \{ \vec{j} \in \Omega_K^4: \vec{j} \in F, j_1=j_2=\ompar \;\;\mbox{ or } \;\; \vec{j} \notin F, j_3=j_4=\ompar \},$$
in which case no edge is added to $\BS^{ I}(t)$.

Let $\bar x_i(t), a_i^*(t),b^*(t), c^*_i(t)$, $1 \le i \le K$, $t \in [t_n, t_n^+]$, be the processes determined from
$\{\BS(t)\}_{t\in [t_n, t_n^+] }$ as in Section \ref{sec:main-bsr-susceptibility}.  These processes will be used in the second round of the construction.

% Note that given $\{\BS^{ I}(t)\}_{t \in [t_n, t_n^+] }$, we have enough information to determine the density of vertices in small components $\bar x_i(t), i \in \Omega_K$ as well as the stochastic rate functions $(a_i^*(t),b^*(t), c^*_i(t))$ for $1 \le i \le K$, which will be used in the second round of exposure.

Now define $\BS^{I,*}(t)$ to be the subgraph that consists of all large components (components of size greater than $K$) in $\BS^{I}(t)$, and then define $\bfG^{I,-}(t)$ and $\bfG^{I,+}(t)$ for $t \in [t_n, t_n^+]$ as follows:
$$ \bfG^{I,-}(t) \equiv \BS^{I,*}(t_n), \mbox{ and } \bfG^{I,+}(t) \equiv \BS^{I,*}(t_n^+).$$
Then 
$$ \bfG^{I,-}(t) \subset \BS^{I,*}(t) \subset \bfG^{I,+}(t) \mbox{ for all } t \in [t_n ,t_n^+].$$

We now proceed to the second round of the construction. Let
$$ E_n =\set{ b(t_c) -n^{-1/6} < b^*(t) < b(t_c) + n^{-1/6}, \mbox{ for all } t \in [t_n ,t_n^+] }. $$
Note that Lemma \ref{lemma:approx-xi} and \eqref{eqn:def-fb-b} implies that with probability at least $1- C_1 e^{-C_2 n^{1/5}}$,
\begin{align*}
	\sup_{ t \in (t_n, t_n^+)}|b^*(t)-b(t_c)| 
	\le& \sup_{ t \in (t_n, t_n^+)}|b^*(t)-b(t)| + \sup_{ t \in (t_n, t_n^+)}|b(t)-b(t_c)|\\
	\le& d_1 n^{-2/5} + d_2 n^{-\gamma} = o(n^{-1/6}).
\end{align*}
Thus $\prob \set{E_n^c} \to 0$ as $n \to \infty$. 
Since we only need the coupling to be good with high probability, 
it suffices to construct the coupling of the three processes until the first time $t \in [t_n, t_n^+]$ when
$b^*(t) \not \in [b(t_c) -n^{-1/6}, b(t_c) + n^{-1/6}]$.  Equivalently, we can assume without loss of generality that
$b^*(t)  \in [b(t_c) -n^{-1/6}, b(t_c) + n^{-1/6}]$, for all $t \in [t_n, t_n^+]$, a.s.

%so we can define the three processes arbitrarily on $E_n^c$.\\

We will construct $\bfG^{ +}(t)$, $\BS^*(t)$ and $\bfG^{ -}(t)$ by adding new edges between  components in the three random graph processes  $\bfG^{I,-}(t)$, $\BS^{I,*}(t)$ and $\bfG^{I,+}(t)$ such that, at time $t \in [t_n ,t_n^+]$ edges are added
between each pair of vertices in  $\bfG^{I,-}(t)$, $\BS^{I,*}(t)$ and $\bfG^{I,+}(t)$,  at  rates $\frac{1}{n}(b(t_c)-n^{-1/6})$, $\frac{1}{n}b^*(t)$ and $\frac{1}{n}(b(t_c)+n^{-1/6})$, respectively.  The precise mechanism is as follows.

We first construct $\bfG^{+}(t)$ for $t \in (t_n,t_n^+]$ by adding edges between every pair of vertices in $\bfG^{I,+}(t)$
 at the rate $\frac{1}{n}(b(t_c)+n^{-1/6})$ and creating self-loops  at the rate $\frac{1}{2n}(b(t_c)+n^{-1/6})$ for each vertex in $\bfG^{I,+}(t)$. 
%Denote such an random graph process by $\{\bfG^{+}(t)\}_{t \in [t_n,t_n^+]}$. This will give a \erdos type random graph with multi-edges and self-loops.

Next, we  construct $\BS^*(t)$ and $\bfG^{-}(t)$ through successive thinning of $\bfG^{+}(t)$, thus obtaining the desired coupling. Let $(e_1,e_2, ...)$ be the sequence of edges that are added to $\bfG^{I,+}(t)$ to obtain $\bfG^{+}(t)$. 
Let $(u_1, u_2, ...)$ be i.i.d Uniform$[0,1]$ random variables that are also independent of the random variables used to construct
$\bfG^{I,-}, \BS^{I,*}, \bfG^{I,+}, \bfG^{+}$.
Suppose at time $t_k$, we have $\bfG^{+}(t_k)=\bfG^{+}(t_k-) \cup \{ e_k\}$, where  $e_k = \{v_1,v_2\}$.  We set $\BS^*(t_k) = \BS^*(t_k-) \cup \{ e_k\}$ if and only if
$$ v_1, v_2 \in \BS^{I,*}(t_k-) \mbox{ and } u_k \le \frac{b^*(t_k)}{b(t_c) + n^{-1/6}}, $$
otherwise let $\BS^*(t_k) = \BS^*(t_k-)$. 
This defines the process $\BS^*(t)$ (with the correct probability law) such that the second inclusion in \eqref{eq:eq1136} is satisfied.
Finally, construct $\bfG^{-}(t)$ by a thinning of $\BS^*(t)$
exactly as above by replacing $\frac{b^*(t_k)}{b(t_c) + n^{-1/6}}$ with $\frac{b(t_c) - n^{-1/6}}{b^*(t_k)}$.
Then $\bfG^{-}(t)$, for $t \in [t_n ,t_n^+]$ is
 an \erdos type processes and the first inclusion in \eqref{eq:eq1136} is satisfied.
This completes the proof of the first part of the lemma.

We now prove (ii). Consider first the case $\bullet = -$.  We will apply  Theorem \ref{theo:aldous-full-gene}.  With notation as in that theorem,
it follows  from the \erdos dynamics of $\bfG^{\sss(n),-}(t)$ that, the distribution of  $(\bar \bfC^{\sss(n),-}(\lambda), \bar \bfY^{\sss(n),-}(\lambda))$, conditioned on $\{\PP_{\vec{v}}(t),\, t \le t_n \; \vec{v} \in [n]^4\}$, for each
$\lambda \in [-\lambda_n, \lambda_n]$, is same as the distribution of $\bfZ(z^{\sss (n)}, q^{\sss (n)})$, where $z^{\sss (n)} = (\bar \bfC^{\sss(n),-}(-\lambda_n), \bf{0})$, $\bf{0}$ denotes the vector $(0, 0, \cdots)$
and  $q^{\sss (n)}$ is determined by the equality
% First we treat $\bfG^{\sss(n),-}(\lambda)$. In order to apply Theorem \ref{theo:aldous-full-gene}, we identify $q^{\sss (n)}$ in the Theorem as follows. The components in $\BS^*(t_n)$ corresponds to vertices in Theorem \ref{theo:aldous-full-gene} and the rescaled component sizes corresponds to the weights of vertices. Thus from the following equation and the fact $\alpha b(t_c)=1$ (Lemma \ref{lemma:alpha-beta}),
$$
q^{\sss (n)} \bar \bfC^{\sss(n),-}_i(-\lambda_n) \bar \bfC^{\sss(n),-}_j(-\lambda_n)
= \frac{\alpha\beta^{2/3}}{n^{1/3}} (\lambda + \lambda_n) \frac{(b(t_c) - n^{-1/6})}{n} |\CC_i(\bfG^{\sss(n),-}(t_n))| |\CC_j(\bfG^{\sss(n),-}(t_n))|,$$
for $i \neq j$. Recalling that $\alpha b(t_c)=1$ it then follows that 
$q^{\sss (n)} = \lambda + \frac{n^{1/3-\gamma}}{\alpha\beta^{2/3}} + O(n^{1/6-\gamma}).$
We now verify the conditions of Theorem \ref{theo:aldous-full-gene}.  Taking $x^{\sss (n)} = \bar \bfC^{\sss(n),-}(-\lambda_n)$ we see with,
$x^*, s_k$, $k = 1,2,3$ as in Theorem \ref{theo:aldous-full-gene},
% $$ \frac{\beta^{1/3}}{n^{2/3}} \CC_i^{\sss (n)}(t_n)\cdot \frac{\beta^{1/3}}{n^{2/3}} \CC_j^{\sss (n)}(t_n)\cdot q^{\sss (n)} = \CC_i^{\sss (n)} \CC_j^{\sss (n)}\cdot \left( \frac{\alpha \beta^{2/3}\lambda}{n^{1/3}} + \frac{1}{n^{\gamma}}  \right) \cdot \frac{1}{n}(b(t_c)-n^{-1/6}),$$
% we have $q^{\sss (n)}= \lambda + n^{1/3-\gamma}/\alpha\beta^{2/3} + O(n^{1/6-\gamma})$. For $s_k^{\sss (n)}$, $k=1,2,3$, and $x^{*}$ in Theorem \ref{theo:aldous-full-gene}, we have
$$ s^{\sss(n)}_1 \le \beta^{1/3} n^{1/3}, \;\; s^{\sss (n)}_2 = \frac{\beta^{2/3}}{n^{4/3}} \sum_{\CC \subset \BS^*(t_n)} |\CC|^2, \;\; s^{\sss (n)}_3 = \frac{\beta}{n^2} \sum_{\CC \subset \BS^*(t_n)} |\CC|^3.$$
Recall the definition of $\bar s_k$ and $\bar s_{k,\ompar}$ from \eqref{eqn:suscept-defn} and Section \ref{sec:main-bsr-susceptibility}.
%that $\bar s_{k,\ompar}(t) = \sum_{\CC \subset \BS^*(t_n)} \CC^k /n$ and $\bar s_{k}(t) = \sum_{\CC \subset \BS(t_n)} \CC^k /n$, for $k=2,3$. 
Then
$$ s_2^{\sss (n)}=\frac{\beta^{2/3} \bar s_{2,\ompar}(t_n)}{n^{1/3}},\;\; s_{3}^{\sss (n)} = \frac{\beta \bars_{3,\ompar}(t_n)}{n}, \;\; 
x^{*\sss(n)}= \beta^{1/3}\frac{I(t_n)}{n^{2/3}}. $$
%To verify the conditions in Theorem \ref{theo:aldous-full-gene} it suffices to show 
From the first two convergences in Proposition \ref{prop:three-conditions} and recalling that, for $k=1,2$,
$|\bar s_{k, \ompar}-\bar s_k| \le K^k$, we immediately get that the first two convergences in \eqref{eqn:qsigma2} hold.
Also,
$$\frac{x^*}{s_2} = \frac{I(t_n)}{\beta^{1/3}n^{1/3}\bar s_{2,\ompar}(t_n)} = \frac{I(t_n)}{\beta^{2/3}n^{ \gamma + 1/3}}O(1) \to 0 , \mbox{ in probability},$$
where the second equality is consequence of the second convergence in Proposition \ref{prop:three-conditions}, and the convergence of the last term follows from the third convergence in Proposition \ref{prop:three-conditions}. This proves the third convergence in \eqref{eqn:qsigma2}.  

Finally we note that the convergence in \eqref{eqn:additional-condition} holds with $\varsigma = \frac{1}{1 - 3(\gamma + \kappa)}$, where $\kappa$ is as in
Proposition \ref{prop:three-conditions}, since
$$
s_1 \left(\frac{x^*}{s_2}\right)^{\varsigma} \le  O(1) n^{1/3} \left (\frac{I(t_n)}{n^{\gamma + 1/3} }\right)^{\varsigma} = O(1) \left (\frac{I(t_n)}{n^{2\gamma + \kappa} }\right)^{\varsigma}
\to 0,$$
where the last equality follows from  our choice of $\varsigma$
and the convergence is a consequence of Proposition \ref{prop:three-conditions}. Thus we have verified all the  conditions in Theorem \ref{theo:aldous-full-gene}
and  therefore we have from this result that
	$(\bar \bfC^{\sss(n),-}(\lambda), \bar \bfY^{\sss(n),-}(\lambda))$
	converges in distribution, in $\udown$, to  $(\bfX^*(\lambda), \bfY^*(\lambda))$ proving part (ii) of the lemma for $\bullet = -$.

% 	$$ \frac{\bar s_3(t_n)}{(\bars_2(t_n))^3} \convp \beta, \;\; \frac{n^{1/3}}{\bars_2(t_c)}-\frac{n^{1/3-\gamma}}{\alpha} \convp 0, \;\; \frac{\CC_1^{\sss (n)}(t_n)}{n^{2\gamma + 0.01}} \convp 0, $$
% which has been proved in Section \ref{sec:main-bsr-susceptibility}.
%  Thus since Proposition \ref{prop:three-conditions} is true, we can apply Theorem \ref{theo:aldous-full-gene} to $\bfG^{\sss(n),-}(\lambda)$ and establish the convergence in d-metric:
% \bes
% (\bar \bfC^{\sss(n),-}(\lambda), \bar \bfY^{\sss(n),-}(\lambda)) \convd (\bfX^*(\lambda), \bfY^*(\lambda)).
% \ees

To prove part (ii) of the lemma for $\bullet = +$, one needs slightly more work. 
Once more we will apply Theorem   \ref{theo:aldous-full-gene}.  
As before, conditioned on $\{\bar \bfC^{\sss(n),+}(\lambda_0): \lambda_0 \le -\lambda_n\}$, for 
each
$\lambda \in [-\lambda_n, \lambda_n]$,  the distribution of  $(\bar \bfC^{\sss(n),+}(\lambda), \bar \bfY^{\sss(n),+}(\lambda))$  is same as the distribution of $\bfZ(\bar z^{\sss (n)}, \bar q^{\sss (n)})$, where $\bar z^{\sss (n)} = (\bar \bfC^{\sss(n),+}(-\lambda_n), \bf{0})$
and  
$\bar q^{\sss (n)} = \lambda + \frac{n^{1/3-\gamma}}{\alpha\beta^{2/3}} + O(n^{1/6-\gamma}).$
 Taking $x^{\sss (n)} = \bar \bfC^{\sss(n),+}(-\lambda_n)$ we see with,
$x^*, s_k$, $k = 1,2,3$ as in Theorem \ref{theo:aldous-full-gene},
% $$ \frac{\beta^{1/3}}{n^{2/3}} \CC_i^{\sss (n)}(t_n)\cdot \frac{\beta^{1/3}}{n^{2/3}} \CC_j^{\sss (n)}(t_n)\cdot q^{\sss (n)} = \CC_i^{\sss (n)} \CC_j^{\sss (n)}\cdot \left( \frac{\alpha \beta^{2/3}\lambda}{n^{1/3}} + \frac{1}{n^{\gamma}}  \right) \cdot \frac{1}{n}(b(t_c)-n^{-1/6}),$$
% we have $q^{\sss (n)}= \lambda + n^{1/3-\gamma}/\alpha\beta^{2/3} + O(n^{1/6-\gamma})$. For $s_k^{\sss (n)}$, $k=1,2,3$, and $x^{*}$ in Theorem \ref{theo:aldous-full-gene}, we have
$$ s^{\sss(n)}_1 \le \beta^{1/3} n^{1/3}, \;\; s^{\sss (n)}_2 = \frac{\beta^{2/3}}{n^{4/3}} \sum_{\CC \subset \BS^{I,*}(t_n^+)} |\CC|^2, \;\; s^{\sss (n)}_3 = \frac{\beta}{n^2} \sum_{\CC \subset \BS^{I,*}(t_n^+)} |\CC|^3.$$
Next note that for any component $\CC \subset \bfG^-(t_n)= \BS^{I,*}(t_n) $ there is a unique  component $ \CC^+ \subset \bfG^+(t_n) =\BS^{I,*}(t_n^+)$, 
such that $\CC \subset \CC^+$. Denote by $\CC_i$  the $i$-th largest component in $\BS^{I,*}(t_n)$, and let $\CC_i^+$
 be the corresponding component in $\BS^{I,*}(t_n^+)$ such that $\CC_i \subset \CC_i^+$. Denote by $N$  the number of immigrations that occur during $[t_n, t_n^+]$ in $\BS^{I,*}$, and denote by $\{\tilde \CC_i^+\}_{i=1}^N$  the components 
in $\BS^{I,*}(t_n^+)$
resulting from these immigrations.
Then
$$ s_2^{\sss (n)}=\frac{\beta^{2/3} \bar s_{2}^+}{n^1/3},\;\; s_{3}^{\sss (n)} = \frac{\beta \bars_{3}^+}{n}, \;\; 
x^{*\sss(n)}= \beta^{1/3}\frac{I^+}{n^{2/3}}, $$
where
\begin{align*}
	\bar s_2^+ :=& \frac{1}{n} \left( \sum_{i=1}^\infty |\CC^+_i|^2 + \sum_{i=1}^N |\tilde \CC_i^+|^2 \right),\\
	\bar s_3^+ :=& \frac{1}{n} \left( \sum_{i=1}^\infty |\CC^+_i|^3 + \sum_{i=1}^N |\tilde \CC_i^+|^3 \right),\\	
	I^+ :=& \max\set{ \max_{i} |\CC_i^+|, \max_{i} |\tilde \CC_i^+|}.
\end{align*}
To complete the proof it suffices to show that the statement in Proposition \ref{prop:three-conditions} holds with 
$(\bars_2(t_n), \bar s_3(t_n), I^{\sss (n)}(t_n))$ replaced with $(\bar s_2^+, \bar s_3^+, I^+)$.  This follows from Proposition
\ref{prop:excess-vertices} given below and hence completes the proof of the lemma. \qed \\

%In order to show that $\bar s_2^+, \bar s_3^+$ and $I^+$ satisfy Proposition \ref{prop:three-conditions}, one needs to show the following proposition.
\begin{Proposition}
	\label{prop:excess-vertices}
	With notation as in the proof of Lemma \ref{lemma:upper-lower-coupling}, as $n \to \infty$, we have
	$$ I^+ = O(I), \;\; \frac{\bar s_2^+}{\bar s_2(t_n)} \convp 1, \;\; \frac{\bar s_3^+}{\bar s_3(t_n)} \convp 1, 
	\;\; \frac{n^{1/3}}{\bar s_2(t_n)} - \frac{n^{1/3}}{\bar s_2^+} \convp 0.$$
\end{Proposition}

{\bf Proof.} The proof is similar to that of  Proposition 8.1 in \cite{bhamidi-budhiraja-wang2011} thus we only give a sketch. 

Observe that the total rate of attachments is $\sum_{i=1}^K c^*_i(t) \le 1$ and each attachment has size no bigger than $K$. 
Recall that  $\CC_i$  denotes the $i$-th largest component in $\BS^{I,*}(t_n)$.  Denote 
by $V_i(t)$, $t \in [t_n, t_n^+]$, the stochastic process defining the size of the component containing $\CC_i$
in $\BS^{I,*}(t)$.  Note that
  $V_i(t_n) = |\CC_i|$ and $V_i(t_n^+)=|\CC_i^+|$. Then ${V_i(t)}/{K}$ can be stochastically dominated by a Yule process starting with $\lceil |\CC_i|/K\rceil $ particles and birth rate $K$. 
Using this and an argument similar to \cite{bhamidi-budhiraja-wang2011}, it follows that, 
$$ |\CC_i^+| - |\CC_i| \le_d K \cdot \mbox{Negative-Binomial}(\lceil |\CC_i|/K\rceil, e^{-2Kn^{-\gamma}}). $$
Next, note that the immigrations are of size no bigger than $2K$, and thus for the same reason, we have the  bound,
$$ |\tilde \CC_i^+| \le_d 2K + K \cdot \mbox{Negative-Binomial}(2, e^{-2Kn^{-\gamma}}).$$
Since the total number of vertices is $n$,  the number of immigrations $N$ can be bounded by $n/K$.

With the above three bounds the proof of the proposition follows exactly as the proof of Proposition 8.1 in \cite{bhamidi-budhiraja-wang2011} with obvious
changes  needed  due to the constant $K$ that appears in the above bounds. Details are omitted. \qed\\

% In this section, we always assume the connection between the two time scales $t$ and $\lambda$ as $t=t_c + \alpha \beta^{2/3}\lambda/n^{1/3}$. Suppose $\lambda_n$ is such that $t_n=t_c-n^{-\gamma}=t_c - \alpha \beta^{2/3}\lambda_n/n^{1/3}$, then
% \begin{equation}
% \label{eqn:lamdan-def}
% 	-\lambda_n = -\frac{ n^{1/3-\gamma}}{\alpha \beta^{2/3}}.
% \end{equation}
% Also, in this section we will mainly work with $\BS^*(t)$, instead of $\BS(t)$. We can also define the rescaled size and surplus vector $\bar \bfC^{\sss (n)}(\lambda)$ for $\BS^*(t)$ similar to the one defined in Theorem \ref{thm:crit-regime}. In this section we also use $(\bar \bfC^{\sss (n)}(\lambda), \bar \bfY^{\sss (n)}(\lambda))$ for the analogue parts defined for $\BS^*(t)$.
% 
% \subsubsection{Proof of Theorem \ref{thm:crit-regime}: the sandwich method}
% \label{sec:generalized-aldous}
We will now use Lemma \ref{lemma:upper-lower-coupling} to complete the proof of Theorem \ref{theo:onedimc}.
We begin with the following elementary lemma.
\begin{Lemma}
	\label{lemma:l1cgce}
	Let $\{x^{\sss (n)}_i, y^{\sss (n)}_i, x_i,  y_i, i \ge 1, n\ge 1\} $ be a collection of non-negative numbers such that, for each fixed $i$, as $n\to \infty$, 
	$x^{\sss (n)}_i \to   x_i$ and  $y^{\sss (n)}_i\to y_i$.
	Also suppose that 
	$\sum_i x^{\sss (n)}_iy^{\sss (n)}_i \to \sum_i x_iy_i < \infty$ as $n \to \infty$.
	Then $\sum_i|x^{\sss (n)}_iy^{\sss (n)}_i -  x_i y_i| \to 0$ as $n\to \infty$.
\end{Lemma}
	{\bf Proof.}  Proof is immediate on applying Fatou's lemma, indeed
	\begin{align*}
		2\sum_{i} x_iy_i \le& \liminf_{n\to \infty} \sum_i (x^{\sss (n)}_iy^{\sss (n)}_i +  x_i y_i - |x^{\sss (n)}_iy^{\sss (n)}_i -  x_i y_i|)\\
	=& 2 \sum_i x_iy_i - \limsup_{n\to \infty} \sum_i |x^{\sss (n)}_iy^{\sss (n)}_i -  x_i y_i|.\end{align*}
\qed\\

The next proposition says that the inclusion in \eqref{eq:eq1136} can be strengthened to component-wise inclusion.
\begin{Proposition}
	\label{prop:surplus} Fix $\lambda\in \Rbold$ and $i_0\geq 1$.  Then, as $n\to\infty$,
	\[\prob\set{\CC_i(\bfG^{\sss(n), -}(t^\lambda))\subset \CC_i(\BS^{*}(t^\lambda)) \subset \CC_i(\bfG^{\sss(n), +}(t^\lambda))~~~ \forall~ 1\leq i\leq i_0} \to 1.\]
\end{Proposition}
{\bf Proof: }
From Lemma \ref{lemma:upper-lower-coupling}  and Lemma 15 in \cite{aldous2000random} (see also Section 8.2 of \cite{bhamidi-budhiraja-wang2011} for a similar
argument), we have, as $n\to\infty$,
\begin{equation}
	\label{eq:eq1154}(\bar{\bfC}^{\sss(n),-}(\lambda),\bar{\bfC}^{{\sss(n)},*}(\lambda), \bar{\bfC}^{\sss(n),+}(\lambda) ) \convd (\bfX(\lambda), \bfX(\lambda), \bfX(\lambda)),\end{equation}
in $\ldown \times \ldown \times \ldown$, where $\bfX$ is as in Theorem \ref{thm:smc-surplus}.
 Define events $E_n, F_n$ as
$$ E_n=\set{ \bar \bfC_i^{\sss (n),-}(\lambda) > \bar \bfC_{i+1}^{\sss (n),+}(\lambda) : 1\le i \le i_0 }, F_n=\set{ \bfG^{\sss (n),-}(\lambda) \subset \BS^*(\lambda)  \subset \bfG^{\sss (n),+}(\lambda) }. $$
Then on the set $E_n \cap F_n$ 
$$ \CC_i(\bfG^{\sss(n), -}(\lambda))\subset \CC_i(\BS^{*}(\lambda)) \subset \CC_i(\bfG^{\sss(n), +}(\lambda)), ~~\forall~ 1\leq i\leq i_0.$$
From Lemma \ref{lemma:upper-lower-coupling} (i) $\prob\{ F_n^c\} \to 1$. Also
\begin{align*}
	\limsup_{n} \prob( E_n^c) 
	\le& \limsup_n \prob \set{\bar \bfC_i^{\sss (n),-}(\lambda) \le \bar \bfC_{i+1}^{\sss (n),+}(\lambda)  \mbox{ for some } 1\le i \le i_0}\\
	\le& \prob \set{ \bfX_i(\lambda) \le \bfX_{i+1}(\lambda) \mbox{ for some } 1 \le i \le i_0 }=0.
\end{align*}
This shows that $\prob(E_n \cap F_n) \to 1$ as $n \to \infty$.  The result follows.
 \qed\\

We will also need the following elementary lemma.  Proof is omitted.
\begin{Lemma}
	\label{lemma:sandwich-1}
	Let $ \eta ^{\sss (n),-}, \eta^{\sss (n), +}, \eta^*$ be real random variables such that $\eta^{\sss (n),-} \le \eta^{\sss (n),+}$ with high probability. Further assume $\eta^{\sss (n),-} \convd \eta^*$ and $\eta^{\sss (n),+} \convd \eta^*$. Then 
	$ \eta^{\sss (n),+} - \eta^{\sss (n),-} \convp 0. $
	Furthermore, if $\eta^{\sss (n)}$ are random variables such that $\eta^{\sss (n),-} \le \eta^{\sss (n)} \le \eta^{\sss (n),+}$ with high probability, then 
	$ \eta^{\sss (n)} \convd \eta^*$  and  $\eta^{\sss (n)} - \eta^{\sss (n),-} \convp 0.$
\end{Lemma}

We now complete the proof of Theorem \ref{theo:onedimc}.\\

\textbf{Proof of Theorem \ref{theo:onedimc}:} 
From Lemma \ref{lemma:upper-lower-coupling} (ii) we have that
% Denote $\bar\bfY^{\sss (n),-}(\lambda)$, $bar\bfY^{\sss (n)}(\lambda)$ and $bar\bfY^{\sss (n),+}(\lambda)$ for the surplus vector of the corresponding random graphs. Then Theorem \ref{theo:aldous-full-gene} shows that 
% $$ (\bar \bfC^{\sss (n),-}(\lambda), \bar \bfY^{\sss (n),-}(\lambda)) \convd (\bfX^*(\lambda), \bfY^*(\lambda)), \;\; (\bar \bfC^{\sss (n),+}(\lambda), \bar \bfY^{\sss (n,+)}(\lambda)) \convd (\bfX^*(\lambda), \bfY^*(\lambda)),$$
% in the space $\udown$. This in particular shows that
\begin{equation} \label{eq:eq1206}
\left(\bar {\bfC}^{\sss (n),-}(\lambda), \bar \bfY^{\sss (n),-}(\lambda) , \sum_{i=1}^\infty \bar {\bfC}^{\sss (n),-}_i(\lambda)\bfY_i^{\sss (n),-}(\lambda)\right)
\convd \left( {\bfX}(\lambda), \bfY(\lambda), \sum_{i=1}^\infty  {\bfX}_i(\lambda) \bfY_i(\lambda)\right),\end{equation}
in $\ldown \times \NNN^{\infty} \times \RRR$, where on $\NNN^{\infty}$ we consider the product topology.

In order to prove the theorem it suffices, in view of Lemma \ref{lemma:l1cgce},  to show that
\be
\left(\bar {\bfC}^{\sss (n)}_*(\lambda), \hat \xi^{\sss (n)}(\lambda) , \sum_{i=1}^\infty \bar {\bfC}_i^{{\sss (n)},*}(\lambda)\hat \xi_i^{\sss (n)}(\lambda)\right)
\convd \left( {\bfX}(\lambda), \bfY(\lambda), \sum_{i=1}^\infty  {\bfX}_i(\lambda) \bfY_i(\lambda)\right), \label{eqn:2576}
\ee
in $\ldown \times \NNN^{\infty} \times \RRR$. 
% Firstly note that 
% \begin{align}
% &\left(\bar {\bfC}^{\sss (n),-}(\lambda),  \bar {\bfC}^{\sss (n)}(\lambda), \bar \bfY^{\sss (n),-}(\lambda) , \sum_{i=1}^\infty \bar {\bfC}^{\sss (n),-}_i(\lambda)\xi_i^{\sss (n),-}(\lambda)\right)\nonumber\\
% \convd& \left( {\bfX}^*(\lambda), {\bfX}^*(\lambda), \bfY^{*}(\lambda), \sum_{i=1}^\infty  {\bfX}^*_i(\lambda) \bfY_i^{*}(\lambda)\right), \label{eqn:2580}
% \end{align}
% with the topology of $\ldown \times \ldown \times \NNN \times \RRR$. 
From Proposition \ref{prop:surplus}, we have for any $i_0 \in \NNN$, with high probability
$$  \bar \bfY_i^{\sss (n),-}(\lambda) \le \hat \xi_i^{\sss (n)}(\lambda) \le \bar \bfY_i^{\sss (n),+} \mbox{ for } 1 \le i \le i_0.$$
Also, from Lemma \ref{lemma:upper-lower-coupling} (i),
whp,
$$  \sum_{i=1}^\infty \bar {\bfC}^{\sss (n),-}_i(\lambda)\bar\bfY_i^{\sss (n),-}(\lambda) \le \sum_{i=1}^\infty \bar {\bfC}^{\sss (n)}_i(\lambda)\bar\bfY_i^{\sss (n)}(\lambda) \le \sum_{i=1}^\infty \bar {\bfC}^{\sss (n),+}_i(\lambda)\bar\bfY_i^{\sss (n),+}(\lambda). $$
From Lemma \ref{lemma:sandwich-1} and Lemma \ref{lemma:upper-lower-coupling} (ii), we then have
$$ \left( \left |\hat \xi^{\sss (n)}(\lambda) - \bar \bfY^{\sss (n),-}(\lambda)\right |,\;\; \sum_{i=1}^\infty \bar {\bfC}_i^{{\sss (n)},*}(\lambda)\hat \xi_i^{\sss (n)}(\lambda) -\sum_{i=1}^\infty \bar {\bfC}^{\sss (n),-}_i(\lambda)\bar \bfY_i^{\sss (n),-}(\lambda) \right) \convp 0, $$
 in $\NNN^{\infty} \times \RRR$, where for $y = (y_1, y_2, \cdots ) \in \ZZZ^{\infty}$, $|y| = (|y_1|, |y_2|, \cdots)$.
The convergence in \eqref{eqn:2576} now follows on combining \eqref{eq:eq1206} and  \eqref{eq:eq1154}.  The result follows. \qed\\

\subsection{Proof of Theorem \ref{thm:crit-regime}.}
\label{sec:old-surplus}
As a first step towards the proof we show the following convergence result for one dimensional distributions.
\begin{Theorem}
	\label{theo:onedimcnew}
	For every $\lambda \in \RRR$, as $n\to \infty$, $(\bar {\bfC}^{\sss (n)}(\lambda), \bar \bfY^{\sss (n)}(\lambda))$ converges in distribution, in $\udown$, to
	$({\bfX}(\lambda), \bfY(\lambda))$.
\end{Theorem}
\textbf{Proof.}
Fix $\lambda \in \RRR$.  We first argue that
\begin{equation} \label{eq:eq1536}
(\bar \bfC^{{\sss (n)},*}(\lambda), \bar \bfY^{{\sss (n)},*}(\lambda))	\convd ({\bfX}(\lambda), \bfY(\lambda)), \mbox{ in }\udown.
\end{equation}
For this, it suffices to show that
\begin{equation}\label{eq:eq1539}
\sum_{i=1}^\infty \tilde \xi_i^{\sss (n)}(\lambda) \bar \bfC_{i}^{{\sss (n)},*} (\lambda) \convp 0.	
	\end{equation}
	Define 
	$$E_n=\set{ I(s) \le C_3\frac{(\log n)^4}{(t_c-s)^2} \mbox{ for } s \le t_c-n^{-\gamma} }.$$
	By Theorem \ref{thm:subcrit-reg}, $\prob\{E_n^c\} \to 0$ and $E_n \in \tilde \FF(\lambda) = \sigma \{ |\CC_i(s)|: i\ge 1,  s \le t^{\lambda}\}$ for all $\lambda \ge -\lambda_n$.
We begin by showing that there exists $d_1 \in (0, \infty)$ such that, for all $i \in \NNN$,
\be
\label{eq:eq1601}
\E \left[\tilde \xi_i^{\sss (n)}(\lambda)\mid \tilde{\FF}_\lambda\right ] 1_{E_n} \le d_1 \bar {\bfC}_i^{{\sss (n)},*}(\lambda) n^{\gamma -1/3}(\log n)^4. 
\ee
Note that at any time $s < t^{\lambda}$, for a component of size $\CC \subset \BS^{\sss (n)}(s)$, there are at most $2 |\CC|^2 n^2$ quadruples of vertices which may provide a surplus edge within $\CC$. Since edges are formed at rate $2/n^3$, we have that
% Since the realization of these point processes (recall $\PP _{\vec v}$ are independent point processes with rate  as defined above \eqref{eqn:f-rule-def-cts}) associate with these quadruples after time $s$ are independent of $\tilde \FF_\lambda$, then after conditioning, the rate of surplus in this component is bounded by 
% $$ \frac{1}{2n^3} \cdot 2 |\CC|^2 n^2= |\CC|^2/n.$$
\begin{align*}
	\E \left[\tilde \xi_i^{n}(\lambda)\mid \tilde{\FF}_\lambda\right ] 
	\le& \int_0^{t_n} \left[\sum_{j: \CC_j(\BS^{\sss (n)}(s)) \subset \CC_i(\BS^{*}(t^{\lambda}))}
	\frac{1}{2n^3}2 n^2 |\CC_j(\BS^{\sss (n)}(s))|^2 \right] ds\\
	\le & \frac{1}{n} |\CC_i(\BS^{*}(t^{\lambda}))| \int_0^{t_n}  I(s) ds.
\end{align*}	
Thus, for some $d_0, d_1 \in (0, \infty)$
%using Theorem \ref{thm:subcrit-reg}	
	% \\
	% =& \frac{\CC_1^{\sss (n)}(\lambda)}{n} \int_0^{t_c-n^{-\gamma}}  I(s) 1_{E_n} ds
	% \le \frac{\CC_1^{\sss (n)}(\lambda)}{n} \int_0^{t_c-n^{-\gamma}}  C_3 \frac{(\log n)^4}{(t_c-s)^2} ds\\
	% \le& C_5 \bar {\bfC}^{\sss (n)}_1(\lambda) n^{2/3+\gamma -1} (\log n )^4.
\begin{align*}
\E \left[\tilde \xi_i^{\sss (n)}(\lambda)\mid \tilde{\FF}(\lambda)\right ] 1_{E_n} 
\le d_0\frac{\bar {\bfC}_i^{{\sss (n)},*}(\lambda)}{n^{1/3}} \int_0^{t_c-n^{-\gamma}}  \frac{(\log n)^4}{(t_c-s)^2} ds
\le d_1 
\bar {\bfC}_i^{{\sss (n)},*}(\lambda) n^{\gamma -1/3
} (\log n )^4.
\end{align*}	
This proves \eqref{eq:eq1601}.
As an immediate consequence of this inequality we have that
\begin{align*}
	\E \left [ \sum_{i}\tilde \xi_i^{\sss (n)}(\lambda)\bar {\bfC}_i^{{\sss (n)},*}(\lambda) \mid \tilde\FF(\lambda)\right]1_{E_n}
	=& \sum_{i}\bar {\bfC}_i^{{\sss (n)},*}(\lambda)1_{E_n}\E \left [ \tilde \xi_i^{\sss (n)}(\lambda) \mid \tilde\FF(\lambda)\right] \\
	\le&  d_1  n^{\gamma -1/3} (\log n )^4\sum_{i} \left(\bar {\bfC}_i^{{\sss (n)},*}(\lambda) \right)^2.
\end{align*}
Observing that $\gamma - 1/3 < 0$ and, from Theorem \ref{theo:onedimc}, that $\sum_{i}\left(\bar {\bfC}_i^{{\sss (n)},*}(\lambda) \right)^2$  converges in distribution,
we have that
$$\E \left [ \sum_{i}\tilde \xi_i^{\sss (n)}(\lambda)\bar {\bfC}_i^{{\sss (n)},*}(\lambda) \mid \tilde\FF(\lambda)\right]1_{E_n} \convp 0.$$
Since $\prob(E_n) \to 1$, letting $\eta^{\sss (n)} = \sum_{i}\tilde \xi_i^{\sss (n)}(\lambda)\bar {\bfC}_i^{{\sss (n)},*}(\lambda)$, we have that
$\E(\eta^{\sss (n)} \mid \tilde\FF(\lambda)) \to 0$ in probability.  Convergence in \eqref{eq:eq1539} now follows
on noting that, as $n\to \infty$,
 $$ \E[\eta^{\sss (n)} \wedge 1] = E \left [ E[\eta^{\sss (n)} \wedge 1 \mid \tilde \clf(\lambda)]\right] \le \E \left [ \E [\eta^{\sss (n)}  \mid \tilde \clf(\lambda)]\wedge 1\right] \to 0.$$
This proves \eqref{eq:eq1536}.  Next note that
\begin{equation}
	\label{eq:eq1619}
	\sum_{i=1}^{\infty} |\bar {\bfC}^{\sss (n)}_{i}(\lambda) - \bar {\bfC}_i^{{\sss (n)},*}(\lambda)|^2 \le \frac{n}{n^{4/3}}O(1) \to 0, \mbox{ as }
	n \to \infty .
	\end{equation}
	Also,
	\begin{align*}
		\E \left[\bar \bfY_i^{\sss (n)}(\lambda)\mid \tilde{\FF}(\lambda)\right ] 1_{\{|\CC_i(t^{\lambda})| \le K\}}
		\le& \left[\int_0^{t^{\lambda}} \sum_{j: \CC_j(s) \subset \CC_i(t^{\lambda})}
		\frac{1}{2n^3}2 n^2 |\CC_j(s)|^2 ds\right]  1_{\{|\CC_i(t^{\lambda})| \le K\}}\\
		\le & \frac{K^2}{n}.
	\end{align*}	
Thus, as $n \to \infty$,
$$
\E \sum_{i=1}^{\infty} |\bar {\bfC}^{\sss (n)}_{i}(\lambda)\bar \bfY_i^{\sss (n)}(\lambda) - \bar {\bfC}_i^{{\sss (n)},*}(\lambda)\bar \bfY_{*,i}^{\sss (n)}(\lambda)|
\le \frac{O(1)}{n} \E \left[\sum_{i=1}^{\infty}\bar {\bfC}^{\sss (n)}_{i}(\lambda) \right] = O(n^{-2/3}) \to 0.$$
The result now follows on combining the above convergence with \eqref{eq:eq1619} and \eqref{eq:eq1536}. \qed \\

\begin{Remark}
	\label{rem:rem1905}
	The proofs of Theorems \ref{theo:onedimc} and \ref{theo:onedimcnew} in fact establish the following stronger statement:
	For all $\lambda \in \RRR$,
	\begin{align*}
	&\left (|\bar \bfY^{\sss (n),-}(\lambda) - \bar \bfY^{\sss (n)}(\lambda)|,
	\sum_{i=1}^{\infty}|\bar {\bfC}^{\sss (n),-}_i(\lambda) - \bar {\bfC}^{\sss (n)}_i(\lambda)|^2,\right .\\
	& \left .\;\;\;\;\; \sum_{i=1}^{\infty} |\bar {\bfC}^{\sss (n),-}_i(\lambda)
	 \bar \bfY^{\sss (n),-}_i(\lambda) - \bar {\bfC}^{\sss (n)}_i(\lambda)
	\bar \bfY^{\sss (n)}_i(\lambda)|\right) 
	\to  ({\bf 0},0,0),\end{align*}
	 in probability, in $\NNN^{\infty}\times \RRR \times \RRR$. 
\end{Remark}

%Finally the multi-dimensional convergence in Theorem \ref{thm:crit-regime} is just one step away from us.\\
 
{\bf Proof of Theorem \ref{thm:crit-regime}:}
For simplicity we present the proof for the case $m=2$.  The general case can be treated similarly.
Fix $-\infty < \lambda_1 < \lambda_2  < \infty$. 
Denote, for $\lambda \in \RRR$, $\bar \bfZ^{\sss (n),-}(\lambda) = (\bar {\bfC}^{\sss (n),-}(\lambda), \bar \bfY^{\sss (n),-}(\lambda))$.
In view of Remark \ref{rem:rem1905} it suffices to show that, as $n \to \infty$,
$$(\bar \bfZ^{\sss (n),-}(\lambda_1),\bar \bfZ^{\sss (n),-}(\lambda_2)) \convd (\bfZ(\lambda_1), \bfZ(\lambda_2)),$$
for which it is enough to show that for all $f_1, f_2 \in C_b(\udown^0)$
\begin{equation}
	\label{eq:eq1943}
	\E\left [ f_1(\bar \bfZ^{\sss (n),-}(\lambda_1))f_2(\bar \bfZ^{\sss (n),-}(\lambda_2))\right]
	\to \E\left [ f_1( \bfZ(\lambda_1))f_2( \bfZ(\lambda_2))\right].
\end{equation}
Note that the left side of \eqref{eq:eq1943} equals
$$	\E\left [ f_1(\bar \bfZ^{\sss (n),-}(\lambda_1))\clt_{\lambda_2-\lambda_1}f_2(\bar \bfZ^{\sss (n),-}(\lambda_1))\right],$$
which using Theorem \ref{thm:smc-surplus} (2), Lemma \ref{lemma:upper-lower-coupling} (ii) and the fact that
$\bfX(\lambda) \in \udown^1$ a.s., converges to
$$
\E\left [ f_1( \bfZ(\lambda_1))\clt_{\lambda_2-\lambda_1}f_2( \bfZ(\lambda_1))\right]
= \E\left [ f_1( \bfZ(\lambda_1))f_2( \bfZ(\lambda_2))\right],$$
where the last equality follows from Theorem \ref{thm:smc-surplus} (3).  This proves \eqref{eq:eq1943} and the result follows. \qed \\

% For the lower bound process $\bfG^{\sss (n),-}(\lambda)$, 
% convergence of $\left ((\bar {\bfC}^{\sss (n),-}(\lambda_k), \bar \bfY^{\sss (n),-}(\lambda_k)), k = 1, \cdots m \right)$, in $d$ metric, to
% $\left (({\bfX}^*(\lambda_k), \bfY^{*}(\lambda_k)), k = 1, \cdots m \right)$ is an immediate consequence of the nearly Feller property.
% Same is true with $-$ replaced by $+$.
% To complete the proof it suffices to show that
% $$\sum_i | \bar {\bfC}^{\sss (n),-}_i(\lambda_k)\xi_i^{\sss (n),-}(\lambda_k) -  \bar {\bfC}^{\sss (n)}_i(\lambda_k)\xi_i^{\sss (n)}(\lambda_k)| \convp 0$$
% as $n\to \infty$, for each fixed $k$.
% However, this is immediate from Lemma \ref{lemma:l1cgce} and \eqref{eqn:2576}.\qed\\
% 
% 

%\subsubsection{Analyzing surplus in the critical regime}

%\subsection{Barely supercritical regime}
%\label{sec:super-crit}
%\todo[inline]{Can prove a weaker result using \cite{janson2010phase}. Might make sense just to give idea of the proof and identify $\gamma=\rho^\prime(t_c)$ in terms of $\alpha, \beta$.}

% \section{Appendix}
% \begin{Lemma}
% 	\label{augfeller}
% 	For every $t \ge 0$, $\clt_t(C_b(\udown)) \subset C_b(\udown)$.
% 	
% \end{Lemma}
% \begin{Lemma}
% 	\label{augconsist}
% 	For all $\lambda \in \Rbold$ and $t\ge 0$, $\nu_{\lambda}\clt_t = \nu_{\lambda +t}$.
% 	
% \end{Lemma}
% \todo[inline]{Xuan:  See if you can fill in the proofs of the above lemmas.}
\section*{Acknowledgements}
AB and XW has been supported in part by the National Science Foundation (DMS-1004418, DMS-1016441), the Army Research Office (W911NF-0-1-0080, W911NF-10-1-0158) and the US-Israel Binational Science Foundation (2008466). SB and XW have been supported in part NSF-DMS grant 1105581.

% 
% \newpage
% 
% %\listoftodos
% 
% 
% \newpage


% \bibliographystyle{plain}
% \bibliography{crit-bib}
\end{document}